\documentclass[11pt]{article}
\usepackage{amsmath,amsthm,amssymb,amsfonts}

\newfont{\lie}{eufm10 at 12pt}
\newfont{\liepequenos}{eufm10 at 10pt}

\newfont{\corpos}{msbm10 at 12pt}
\newfont{\corpospequenos}{msbm10 at 10pt}

\newcommand{\CC}{\mbox{\corpos \symbol{67}}}         
\newcommand{\RR}{\mbox{\corpos \symbol{82}}}         
\newcommand{\ZZ}{\mbox{\corpos \symbol{90}}}         

\newcommand{\blu}{\mbox{~}}
\newcommand{\ti}[1]{\mbox{\tiny $#1$}}
\newcommand{\scr}[1]{\mbox{\scriptsize $#1$}}

\newcommand{\sm}[1]{\mbox{\small $#1$}}
\newcommand{\la}[1]{\mbox{\large $#1$}}
\newcommand{\La}[1]{\mbox{\Large $#1$}}
\newcommand{\LA}[1]{\mbox{\LARGE $#1$}}

\newcommand{\lnab}[1]{\la{\nabla}_{\!#1}}

\newcommand{\hlnab}[1]{\hat{\la{\nabla}}_{\!\!#1}}
\newcommand{\lnabp}[1]{\la{\nabla'}_{\!\!#1}}
\newcommand{\lnabo}[1]{\la{\nabla}^{\bot}_{\!\!#1}}
\newcommand{\lnabe}[1]{\la{\nabla}^{E}_{\!\!#1}}

\newcommand{\ra}{\rightarrow}
\newcommand{\non}{\nonumber}
\newcommand{\ha}{\scr{\frac{1}{2}}}

\newcommand{\al}{\alpha}
\newcommand{\be}{\beta}
\newcommand{\ga}{\gamma}
\newcommand{\bal}{\bar{\alpha}}
\newcommand{\bbe}{\bar{\beta}}
\newcommand{\bga}{\bar{\gamma}}
\newcommand{\bmu}{\bar{\mu}}
\newcommand{\brho}{\bar{\rho}}

\newcommand{\gdf}[3]{g(\mbox{\large ${\nabla}$} dF({#1},{#2}),JdF({#3}))}
\newcommand{\gf}[3]{g_{#1}{#2}{#3}}
\newcommand{\Jw}{J_{\omega}}

\newcommand{\Fw}{F^{*}\omega}
\newcommand{\Fws}{(F^{*}\omega)^{\sharp}}
\newcommand{\wt}{\omega^{\bot}}
\newcommand{\wts}{(\omega^{\bot})^{\sharp}}
\newcommand{\Jt}{J^{\bot}}
\newcommand{\lh}{\varphi}
\newcommand{\tUpsilon}{\tilde{\Upsilon}}
\newcommand{\tPsi}{\tilde{\Psi}}

\def\cyclic{\mathop{\kern0.9ex{{+}
\kern-2.2ex\raise-.25ex\hbox{\Large\hbox{$\circlearrowright$}}}}\limits}
\newtheorem{Lm}{Lemma}
\newtheorem{Pp}{Proposition}
\newtheorem{Th}{Theorem}
\newtheorem{Cr}{Corollary}

\begin{document}
\setcounter{page}{1}
\title{Cayley submanifolds of Calabi-Yau 4-folds.}
\author{Isabel M.C.\ Salavessa$^{1*}$ and Ana Pereira do Vale$^{2**}$\\[2mm]}
\protect\footnotetext{$\!\!\!\!\!\!\!\!\!\!\!\!\!$
{\bf MSC 2000:} Primary:53C42,53C55,53C25,53C38; Secondary:57R20,57R45.\\
{\bf ~~Key Words:}
Minimal submanifold, K\"{a}hler angles, Cayley submanifolds,
 K\"{a}hler-Einstein manifold, Residue.\\
$^*$ and $^{**}$ Partially supported by Funda\c{c}\~{a}o Ci\^{e}ncia e 
Tecnologia through POCI/MAT/60671/2004.\\
$^*$Partially supported  by Funda\c{c}\~{a}o Ci\^{e}ncia e Tecnologia through
 Plurianual of CFIF.}
\date{ \em  Dedicated to Jim Eells ~~~~~~~~~~~~~~~~~~~~~~~~~~~~
~~~~~~~~~~~~~~~
~~~~~~~~~~~~~~~~~~~~~~~~~~~~~\em}
\maketitle ~~~\\[-9mm]
{\footnotesize $^1$ Centro de F\'{\i}sica das
Interac\c{c}\~{o}es Fundamentais, Instituto Superior T\'{e}cnico,
Edif\'{\i}cio Ci\^{e}ncia,\\[-1mm] Piso 3,
1049-001 LISBOA, Portugal;~~
e-mail: isabel@cartan.ist.utl.pt}\\
{\footnotesize $^2$ Centro de Matem\'{a}tica,
Universidade do Minho,
Campus de Gualtar\\[-1mm]
4710-057 BRAGA, Portugal;~~
e-mail: avale@math.uminho.pt}\\[5mm]
{\small {\bf Abstract:}
Our main  results are: (1)
The complex  and Lagrangian points  of a non-complex and non-Lagrangian
$2n$-dimensional submanifold  $F\!:\! M\!\ra\! N $, immersed  
 with parallel mean curvature and with equal K\"{a}hler angles into a  
K\"{a}hler-Einstein manifold $(N,J,g)$ of complex dimension $2n$, are 
zeros of finite order of  $\sin^2\theta$ and $\cos^2\theta$
respectively, where $\theta$ is  the common $J$-K\"{a}hler angle.
(2) If $M$ is a  Cayley submanifold of a Calabi-Yau (CY) manifold $N$ 
of complex dimension 4, then $\bigwedge^2_+NM$ is naturally isomorphic
to $\bigwedge^2_+TM$. (3) If $N$ is Ricci-flat (not necessarily CY) 
and $M$ is a  Cayley submanifold,  
then $p_1(\bigwedge^2_+NM)=p_1(\bigwedge^2_+TM)$ still holds,
but $p_1(\bigwedge^2_-NM)-p_1(\bigwedge^2_-TM)$ may describe
a residue on the $J$-complex points, 
in the sense of Harvey and
Lawson. We describe this residue by a PDE on
a natural morphism $\Phi:TM\ra NM$, $\Phi(X)=(JX)^{\bot}$, with
singularities at the complex points.  We give an explicit formula
of this residue in a particular case. When
 $(N,I,J,K,g)$ is a hyper-K\"{a}hler manifold
and $M$ is an $I$-complex closed 4-submanifold,
the first Weyl curvature invariant of $M$ may be
described as a residue  on the $J$-K\"{a}hler angle
at the  $J$-Lagrangian points by a Lelong-Poincar\'{e} type formula.  
We study the almost complex structure $\Jw$ on $M$ induced by $F$.}
\markright{\sl\hfill  Salavessa--Pereira do Vale \hfill}
\section{Introduction}
\setcounter{Th}{0}
\setcounter{Pp}{0}
\setcounter{Cr} {0}
\setcounter{Lm} {0}
\setcounter{equation} {0}
\baselineskip .5cm

The role of the complex and anti-complex points on the topology-geometry
of closed non-complex  minimal surfaces immersed into
complex K\"{a}hler surfaces has been studied
in \cite{[We1]}, \cite{[E-S]}, \cite{[E-G-T]}, and \cite{[Wo]}. 
In these papers, it is proved that
the set ${\cal C}={\cal C}^+\cup {\cal C}^-$ of complex and
anti-complex points  is a set of  isolated points, and
each of such points is of  finite order.  The order  of the complex and
anti-complex points is defined  as a  multiplicity of
a zero of $(1\pm \cos\theta)$, where $\theta$ is the K\"{a}hler angle,
and  adjunction formulas were obtained in \cite{[We1]},
\cite{[E-G-T]} and  \cite{[Wo]}:\\
\begin{equation}
 -\sum_{p\in {\cal C}^-} order(p)-\sum_{p\in {\cal C}^+} order(p)
= {\cal X}(M) + {\cal X}(NM)
\end{equation}
\begin{equation}
  \sum_{p\in {\cal C}^-} order(p) -\sum_{p\in {\cal C}^+} order(p)
= F^{*}c_I(N)[M].
\end{equation}
The proofs of these formulas come, respectively, from the
following  PDEs of second order
on the cosine of the K\"{a}hler angle, with
singularities at complex and anti-complex points:
\begin{equation}
\ha\Delta\log\sin^2\theta =  (K^M+K^{\bot})
\end{equation}
\begin{equation}
\ha \Delta \log\left(\frac{1+\cos\theta}{1-\cos\theta}\right)=
- Ricci^N(e_1,e_2),
\end{equation}
where $K^M$ and $K^{\bot}$ are respectively the Gaussian curvature of $M$
and the curvature of the normal bundle $NM$,
 and $e_1, e_2$ is a direct orthonormal frame of $M$. Therefore, 
(1.1) and (1.2)  are formulas that describe some polynomials  
of topological invariants of the immersed surface, normal bundle and ambient 
space, as residue formulas of  certain functions that have singularities at 
those special points.\\[-3mm] 

In higher dimensions, the papers \cite{[L]}, \cite{[We2]}, \cite{[We3]}
show how Pontrjagin classes and Euler classes of a closed (generic)
submanifold $M$ of a complex manifold $(N,J)$ are carried by subsets of
CR-singular points, that is, points with sufficiently many complex
directions. The investigation of complex tangents on a $m$-dimensional
submanifold $M$ embedded into a K\"{a}hler manifold $N$ of complex dimension
$m$ is very much justified, by the well known embedding theorem of Whitney.
More generally, if $N$ has a calibration $\Omega$ of rank $m$ (see definitions
in \cite{[H-L1]}) and
$M$ is not $\Omega$-calibrated 
we may expect that
$\Omega$-calibrated points may have a similar role (\cite{[S2]}). Minimality
of $M$ should guarantee the order of such points to be finite.
\\[-3mm]
 
In \cite{[H-L2]}, \cite{[H-L3]} a general 
framework is shown to obtain this sort of geometric residues, inspired
by the above examples. Given two Riemannian vector
bundles $(E, g_E)$, $(F, g_F)$ over $M$, of the same rank $m$,
  with Riemannian connections $\lnab{}^E$ and $\lnab{}^F$,
  and a bundle map $\Phi:F\ra E$, degenerated at a set of points $\Sigma$,
we may compare a
$m$-characteristic classe $Ch$ of $E$ and the one of $F$, 
describing these invariants
using the curvature tensors with
respect to $\lnab{}^{F}$ and $\lnab{}^E$, via Chern-Weil theory.
$\Phi$ induces on
$F$ a singular connection $\lnabp{}={\Phi^{-1}}^*\lnabe{}$, Riemannian
for a degenerated metric, and that makes
$\Phi$ a parallel isometric bundle map, but $R'$ and
$Ch(R')$ can be smoothly extended to $\Sigma$ by the identities
$ R'(X,Y,Z,W)= g_E(R^E(X,Y)\Phi(Z),\Phi(W))$, $Ch(R')=Ch(R^E)$.
The difference $Ch(R')-Ch(R^F)$ is of the form $dT$ where $T$ is a 
transgression form with singularities along $\Sigma$.
If  $\Sigma$ is  sufficiently small and regular, the Stokes theorem
reads $\int_{M\sim V_{\epsilon}(\Sigma)} dT=-\int_{\partial
V_{\epsilon}(\Sigma)}T$,
where $V_{\epsilon}(\Sigma)$ is a tubular neighbourhood of $\Sigma$
of radius $\epsilon$, and letting $\epsilon \ra 0$ may describe
$Ch(E)\!-\!Ch(F)$ as a residue of $T$
 along $\Sigma$ and expressed in terms of the zeros of $\Phi$.  

Inspired in this framework, 
the present paper shows some formulas of the type (1.3)-(1.4)
for 4-dimensional submanifolds of certain K\"{a}hler manifolds.
As we will see, to workout such formulas in dimension $>2$ 
is considerably more difficult then in the surface case.
\\[-3mm]

We study the set ${\cal C}$ of complex points and the  set ${\cal L}$ 
of Lagrangian points of  non-holomorphic and non-Lagrangian immersed
submanifolds $F:M\ra N$ of real dimension $2n$ of a K\"{a}hler-Einstein 
 (KE) manifold of complex dimension $2n$, namely if $F$ is immersed
with equal K\"{a}hler angles (e.k.a.$\scr{s}$). A natural bundle map 
$\Phi:TM\ra NM$, $\Phi(X)=(JX)^{\bot}$, is defined
and was first used by Webster. $\Phi$ is
degenerated at points with complex directions, 
 and has maximum norm at Lagrangian points, where it is
an isometry.  If $F$ has e.k.a.$\scr{s}$, 
$\Phi$ is conformal with $\|\Phi(X)\|^2
=\sin^2\theta\|X\|^2$ where $\theta$ is the common K\"{a}hler angle, and
away from  ${\cal L}$, one can define  smooth
almost complex structures $\Jw$ on $M$ and $\Jt$ on the normal  bundle
$NM$ that are naturally inherited from the ambient space, and they
coincide with  the induced complex structure at  complex points.
These almost complex structures, with the K\"{a}hler angle, will be
fundamental for our formulas. 
In section 3, if $n=2$ we study  $\Jw$. 
\\[-3mm]

If $n=2$ and $(N,J,g)$ is  Ricci-flat  KE,
$M$ is a Cayley submanifold if it is minimal and
with equal K\"{a}hler angles. 
 If $N$ is Calabi-Yau these Cayley submanifolds
are calibrated by one of the $S^1$-family of Cayley calibrations
(\cite{[H-L1]},\cite{[J]}).  
The Cayley calibrations $\Omega$ do 
not specify the complex or the Lagrangian points, but induce a natural
isomorphism  $\Omega^{M}:\bigwedge^2_+TM\ra
\bigwedge^2_+NM $, $\langle\Omega^{M}(X\wedge Y), U\wedge V\rangle
=  \Omega(X,Y,U,V)$ (see Prop.3.2)\\[-3mm]

In Section 4
we  prove  that complex and Lagrangian points of a $n$-submanifold with
parallel mean curvature are zeros of a
system  of complex-valued  functions that satisfy
a second-order partial differential system of inequalities of the
Aronszajn type, and so, if the submanifold is not  complex
or Lagrangian, they are zeros of finite order.
These inequalities
are obtained from some estimates on the Laplacian of the pull-back
of the K\"{a}hler form of $N$ by $F$, and on the Laplacian of
 $\Phi:TM\ra NM$. Furthermore, the sets ${\cal C}$ and ${\cal L}$
have  Hausdorff codimension at least 1, and if $M$ is
closed and $n=2$, ${\cal L}$ is a set of Hausdorff codimension at 
least 2.\\[-3mm]

  In Section 5
we prove the following residue-type formula, in the same spirit as 
formulas (1.3) and (1.4):
\begin{Th} If $F:M\ra N$ is  a non-$J$-holomorphic Cayley submanifold 
immersed into a 4-fold Ricci-flat K\"{a}hler manifold (not necessarily 
Calabi-Yau), 
the following equalities hold, for some representatives in the  cohomology
classes of $M$:
\begin{eqnarray}
p_1(\sm{\bigwedge^2_+}NM)&=&   p_1(\sm{\bigwedge^2_+}TM)\\
p_1(\sm{\bigwedge^2_-}NM) &=& 
 p_1(\sm{\bigwedge^2_-}TM)+ \frac{1}{\pi^2}d\eta
\end{eqnarray}
where $\eta=\eta(\Phi)$ is a $3$-form, defined away from the complex 
points, which is given by
\begin{eqnarray}
\eta(\Phi)
&=&-\frac{1}{4}\La{\langle}\, \Phi^{-1}\lnab{}\Phi\, \wedge\,  
 \La{(}\Phi(R^{\bot}) + R^M +\frac{1}{3}
[\Phi^{-1}\lnab{}\Phi,\Phi^{-1}\lnab{}\Phi]
\La{)}\, \La{\rangle}~~~~
\end{eqnarray}
where 
$\Phi(R^{\bot}):\bigwedge^2TM\ra\bigwedge^2TM$ is given by
${\Phi(R^{\bot})(X,Y)(Z)=\Phi^{-1}R^{\bot}(X,Y)\Phi(Z)}$.
Furthermore:\\[1mm]
$(A)$ If $F$ has no complex points $TM$ and $NM$ have the same
Pontrjagin and Euler classes.\\[1mm]
$(B)$ If $d\Phi=0$, or if $g(\lnab{X}\Phi(Y),\Phi(Z))$ is skew symmetric 
on $(Y,Z)$,  then $\theta$ is constant and $\Phi:TM \ra NM$ is a
parallel homothetic diffeomorphism.\\[1mm]
$(C)$ If $g(\lnab{X}\Phi(Y),\Phi(Z))$ is symmetric on $(Y,Z)$, or if
$\bar{R}(X,Y)\Phi=0$, where $\bar{R}$ is the curvature tensor of
$TM^*\otimes NM$,    then
$p_1(\bigwedge^2_-TM)=p_1(\bigwedge^2_-NM)$ holds, or
equivalently $M$ and $NM$ have the same Pontrjagin and Euler classes.
\end{Th}
We will say that $F$ has \em regular homogeneous complex points \em, if
${\cal C}=\bigcup_i \Sigma_i$ is a
disjoint finite union of closed submanifolds $\Sigma_i$
of dimension $d_i\leq 3$, 
and for each $i$,  on a neighbourhood $V$ of $\Sigma_i$
in $M$,
 $\sin\theta=  f_i^{r_i}$ where 
$r_i$ is the common  order of the zeros of $\Phi$ (and of $\|\Phi\|$)
along $\Sigma_i$, and
 $f_i$ is a nonnegative continuous function,
smooth on $V\sim {\cal C}$, such that  $\|\nabla f_i\|$ exists as a positive  
$C^{\mu}$ function on all $V$, with $\mu\geq r_i+2$ and the flow of 
$X_{f_i}=\frac{\nabla f_i}{\|\nabla f_i\|^2}$
can be extended to ${G}_{t_0}=\{(p,w)\in N{\Sigma_i}: \|w\|<t_0\}$
as a $C^{\mu+1}$ diffeomorphism $\xi:G_{t_0}\ra V$. That is, $X_{f_i}$ 
is a multivalued vector field at points $p\in \Sigma_i$, with sublimits
spanning all $T_p\Sigma^{\bot}_i$ and for each $u$ unit vector of 
$T_p\Sigma_i^{\bot}$ it is defined an integral curve 
$\gamma_{(p,u)}(t)=\xi(p, tu)$ with
$\gamma_{(p,u)}(0)=p$ and $\gamma'_{(p,u)}(0)=\frac{u}{c(p)}$,
where $c(p)=\|\nabla f_i\|(p)$. 
This flow map $\xi$ defines for each sufficiently small 
$\epsilon>0$ a diffeomorphism
from $C_{\epsilon}=\{(p,w)\in N\Sigma_i: \|w\|=\epsilon\}$
onto $f^{-1}_i(\epsilon)$. Furthermore, 
for each sufficiently small coordinate chart $y$ of $\Sigma$
we have a Farmi-type coordinate chart $x$ of $V$ of class $C^{\mu +1}$,
extending $y$ and
satisfying $f_i
=\sqrt{x_{d_i+1}^2+ \ldots + x_4^2}$ (see Prop.\ 5.8).
Examples of such functions $f_i$
are the distance function $\sigma$ to a submanifold $\Sigma_i$.
Let  $\pi:N^1\Sigma_i:=C_1\ra \Sigma_i$, $\pi(p,u)=p$, and
$S(p,1)$ the unit sphere of $T_p\Sigma_i^{\bot}\subset T_pM$.
For
$u\in S(p,1)$ and $X\in T_pM$, set $X^{\bot_u}=X-g(X,u)u$, and
define $\varsigma(u)(X)=(\lnab{u}\tilde{X})^{\bot}\in T_pM$ 
where $\tilde{X}$ is 
any smooth
section of $TM$ with $\tilde{X}_p=X$.
\begin{Cr} Assume $M$ is compact and 
$F$ in Theorem 1.1 has  regular homogeneous complex points of
order $r_i$ on $\Sigma_i$.
Let $\tilde{\Phi}=\frac{\Phi}{\|\Phi\|}$
and set for each $(p,u)\in N^1{\Sigma_i}$ and $X$
 vector field on $M$,
\begin{equation}\begin{array}{l}
\Upsilon_i(p,u):=\lnab{u^{r_i}}^{r_i}\Phi(p),~~~~~~~
\Psi_i(p,u)(X_p):= \lnab{u^{r_i-1},X}^{r_i}\Phi(p),\\[1mm]
G_i(p,u)(X_p):=\lnab{u^{r_i},X}^{(r_i +1)}\Phi(p)+r_i\Psi_i(p,u)(
\varsigma(u)(X_p))
\end{array}
\end{equation}
defining smooth sections $\Upsilon_i$ of $\pi^{-1}(TM^*\otimes NM)$ and
$\Psi_i, G_i$ of 
$\pi^{-1}(TM^*\otimes (TM^*\otimes NM))$.
Then 
there exist the following limits
\begin{eqnarray}
\lim_{\epsilon\ra 0}\tilde{\Phi}(\xi(p,\epsilon u))
&=&\frac{1}{{r_i}!c(p)^{r_i}}\Upsilon_i(p,u),\\
\lim_{\epsilon\ra 0}\epsilon\lnab{X}
\tilde{\Phi}(\xi(p,\epsilon u))&=&\frac{1}{{(r_i-1)}!c(p)^{r_i-1}}
\Psi_i(p,u)(X^{\bot_u})~~~~~~~~~\\
\mbox{if~}X\bot \nabla f_i \mbox{~near~}p, 
~~~~\lim_{\epsilon\ra 0}\lnab{X}
\tilde{\Phi}(\xi(p,\epsilon u))&=&\frac{1}{{r_i}!c(p)^{r_i}}G_i(p,u)(X_p),
\end{eqnarray}
with $\frac{1}{{r_i}!c(p)^{r_i}}\Upsilon_i(p,u): T_p M\ra NM_p$ an isometry.
Furthermore, set 
\[T_i^{(1)}(p,u)(X)=\Upsilon_i(p,u)^{-1}\circ \Psi_i(p,u)(X)~~~~~~~~
T_i^{(0)}(p,u)(X)=\Upsilon_i(p,u)^{-1}\circ G_i(p,u)(X).\]
Then
\begin{eqnarray*}
\lefteqn{p_1(\sm{\bigwedge^2_-}NM)[M]-p_1(\sm{\bigwedge^2_-}TM)[M]=
}~~~~~~~~~~~~~~~~~~~~~~~\\[3mm]
=\sum_{i:d_i=2}
-\frac{r_i}{4}\int_{\Sigma_i}
\LA{(}\int_{S(p,1)}
&&
\!\!\!\!\!\!\!\!\!\!\!\!\!\!\!\!\!\!\! {c(p)^{-1}}
\langle T_i^{(1)}(p,u) \wedge \La{(}
\Upsilon(p,u)(R^{\bot})+R^M\La{)}\rangle (*u)
d_{S(p,1)}(u)\LA{)}d_{\Sigma_i}(p)\non\\[-1mm]
+\sum_{k=0}^3\sum_{{i:d_i=k}} \sum_{
{\,\al\!+\!\be\! +\! \ga\! =\!3\!-\!k}
}&&\!\!\!\!\!\!\!\!\!\!\!\!\!\!\!\!\!\!
-\frac{r_i^{3-k}}{12}\int_{\Sigma_i}\!\!\LA{(}\int_{S(p,1)}\!\!\!\!
\!\!\!\!\!\!\!\!c(p)^{-1}\la{\langle}T_i^{(\al)}(p,u)\wedge
 [T_i^{(\be)}(p,u),T_i^{(\ga)}(p,u) ]\la{\rangle}(*u)
d_{S(p,1)}(u)\LA{)}d_{\Sigma_i}(p).~~~~~~~~~~~~~\non\\[-2mm]\non
\end{eqnarray*}
\end{Cr}

In section 6 we prove that
if $M$ is a $J$-complex submanifold and $N$ is Ricci-flat, then (1.5)
still holds. Moreover, if $c_1(M)=0$, then $\bigwedge^2_{+}TM$
and $\bigwedge^2_{+}NM$ are both flat, and $\bigwedge^2_{-}TM$
and $\bigwedge^2_{-}NM$ are both anti-self-dual.\\[-2mm]

 If $N= (N, I,J,K, g)$ is hyper-K\"{a}hler (HK) of complex
dimension 4, and $M$ is an $I$-complex submanifold of complex dimension 2, 
then, considering on $N$ the complex structure 
$J$, $M$ is a Cayley submanifold with a
$J$-K\"{a}hler angle $\theta$ that can assume any value. 
Furthermore, $M$ has a hyper-Hermitian
structure $(M, I,J_{\omega_J}, J_{\omega_K}, g)$,
 defined away from totally complex points.
More generally, if $M$ is an "$I$-K\"{a}hler" Cayley submanifold 
of a Ricci-flat K\"{a}hler  4-fold  $(N,J,g)$, i.e., locally
 on a open dense set of $M\sim {\cal L}$, 
a smooth K\"{a}hler  structure $I$ 
 exists and that anti-commutes with
$\Jw$,  then we conclude (in subsection 3.2) that the 
$J$-K\"{a}hler angle $\theta$ also satisfies the PDE
\begin{equation}
\Delta\log \cos^2\theta = s^M
\end{equation}
where $s^M$ is the scalar curvature of $M$. If $M$ is
closed, this is a residue-type formula
for the first curvature invariant of Weyl of $M$,
$\kappa_2(M) =\ha \int_M s^M \mbox{Vol}_M$, 
 in terms of the zero set $\Sigma$ of  $\cos\theta$, which is the set of the
 the $J$-Lagrangian  points ${\cal L}$ of $M$.
 We prove in section 7:
\begin{Pp} If $N$ is HK and $M$ is a non-totally complex 
closed $I$-K\"{a}hler
submanifold with $I$-K\"{a}hler form $\omega_I$, then there
exist a 
locally finite union of irreducible analytic subvarieties
of complex codimension 1 (i.e analytic surfaces) $\Sigma_i$
and integers $a_i$ such that $\Sigma=\bigcup_i \Sigma_i$ where $\cos\theta$
vanish at homogeneous order $a_i$ along $\Sigma_i$ and
a formula of Lelong-Poincar\'{e} type in terms of characteristic divisors
exist: 
$~\frac{1}{\pi}\kappa_2(M)=-\sum_i a_i \int_{\Sigma_i}\omega_I.$\\[-2mm]
\end{Pp}
\noindent
If $I$ does not exist globally on $M$, we still can obtain a residue formula
under some conditions, and a removable high rank singularity theorem
(see  Proposition 7.1 and Corollary 7.1).

In section 8 we give some examples of complete 
 non-linear Cayley submanifolds of $(\RR^{8}, J_0,g_0)$,
with no complex $J_0$-points,    with
only one complex point,  with a 2-plane set of complex points,
or with a 2-plane set of Lagrangian points. They are all
holomorphic for some other  complex structure of $\RR^{8}$.
We also observe that all submanifolds, and in particular
coassociative ones,  that are graphs of
maps $f:\RR^4 \ra \RR^3$, do not have $J_0$-complex points.
\section{The K\"{a}hler angles}
\setcounter{Th}{0}
\setcounter{Pp}{0}
\setcounter{Cr} {0}
\setcounter{Lm} {0}
\setcounter{equation} {0}
We recall the notion of  K\"{a}hler angles introduced in 
\cite{[S-V0]}, \cite{[S-V1]} and \cite{[T1]}
for an immersed  $2m$-submanifold $F:M\ra N$ of a K\"{a}hler manifold
of complex dimension $2n$ where  $m\leq n$.
We denote by $J$ and  $g$  the complex and Hermitian structure of $N$
and  $\omega(X,Y)=
g(JX,Y)$ its K\"{a}hler form.
Let $NM=(dF(TM))^{\bot}$ be the normal bundle and denote by $(~)^{\bot}$ the
orthogonal projection of $F^{-1}TN$ onto $NM$.
The pullback
2-form $\Fw$ defines at each point  $p\in M$ the K\"{a}hler angles
$\theta_1, \ldots, \theta_m$ of $M$, 
$\theta_{\al}\in [0,\frac{\pi}{2}]$, such that
$\cos\theta_1\geq\ldots\geq \cos\theta_m\geq 0$ and
$\{\pm i\cos\theta_{\al}\}_{\al=1,\ldots,m}$
are the  eigenvalues of the complex extension $\Fw$ to $T^c_pM$.
Polar decomposition of the endomorphism $\Fws=|\Fws|\Jw$ 
where $\sharp$ is the usual musical isomorphism,
defines a partial isometry $\Jw:TM\ra TM$
with the same kernel ${\cal K}_{\omega}$ as $\Fw$. Then
$M=\bigcup_{k=0}^m{\cal L}_{m-k}$ where
a point $p\in M$ is in ${\cal L}_{m-k}$
iff $Rank (\Fw)_p= 2k$. At each  $p\in {\cal L}_{m-k}$,
 we may take
an o.n.\ basis of $T_pM$ of the form $\{X_1,Y_1=\Jw X_1,\ldots,
X_k, Y_k=\Jw X_k, X_{k+1},Y_{k+1}, \ldots, X_m,Y_m\}$
where $\{X_{k+1},Y_{k+1}, \ldots, X_m,Y_m\}$ is any o.n.\ basis
of ${\cal K}_{\omega}$. If $O$ is an open set of $M$ lying in
${\cal L}_{m-k}$, then $X_{\al},Y_{\al}$ can be chosen  smoothly
on a neighbourhood of each point of $O$. 
The complex frame
\begin{equation}
\al=Z_{\al}:= \frac{X_{\al}-iY_{\al}}{2}~~~~~~~~\bal=\bar{Z}_{\al}
~~~~~~~~~~~\al=1, \ldots, m
\end{equation}
diagonalizes $\Fw$, $\Fws(Z_{\al})=i\cos\theta_{\al}Z_{\al}$, and
for $\al\leq k$,  $Z_{\al}\in T^{(1,0)}M$ w.r.t. $\Jw$.
We will use the greek letters $\al,\be,\mu,\ldots$ and their conjugates
to denote both the integer in $\{1,2,\ldots,m\}$ it represents or the
corresponding  complex vector of $T^cM$ above defined in (2.1).
If $M$ is orientable 
then   ${\cal C}^+$ and  ${\cal C}^-$,  are  the set of points such that
$\cos\theta_{\al}(p)=1$ $\forall \al$, and respectively, $\Jw$ defines
the same or the opposite orientation of $M$.
The eigenvalues
$\cos\theta_{\al}$ are only locally Lipschitz on $M$,
while  the product $2\epsilon(p)\cos\theta_1\ldots\cos\theta_m
=\langle \Fw^m, Vol_M\rangle $   is smooth
everywhere, where $\epsilon(p)$ is the orientation of $\Jw$ for
$p\in{\cal L}_0$. For $E$ subspace of $T_pM$ set $E^J=E\cap JE$.\\[5mm]
Let $\wt=\omega_{| NM}$ be the restriction of the K\"{a}hler form
$\omega$ of $N$ to the  normal  bundle $NM$, and
$\wts= |\wts|\Jt$
be its polar decomposition. 
We define the following morphisms
\[\begin{array}{rcccccccc}
\Phi :& TM & \ra  & NM &~~~~~~~~~~~~~~~~~& \Xi :& NM & \ra & TM\\
& X &  \ra  & (JX)^{\bot}&~~~~~~~~~~~~~~& &U & \ra & (JU)^{\top}
\end{array}\]
Note that $JX= \Fws(X)+\Phi(X)$, and $\Phi(X)=0$ iff
$\{X,JX\}$ is a complex direction of $F$.
Similarly for $\wt$ and $\Xi$. 
 $\Phi:(TM^J)^{\bot}\ra NM$, $\Xi: (NM^J )^{\bot}\ra TM$ are 1-1.
Set $2s= dim(TM^J)$, $2t=dim(NM^J)$.
The o.n.\ basis $\{U_{A},V_{A}\}=$
$\{U_{1},JU_1,\ldots,U_t,JU_t, \Phi(\frac{Y_{\al}}{\sin\theta_{\al}}),~
\Phi(\frac{X_{\al}}{\sin\theta_{\al}})\}$, where $\al$ are s.t. 
$\sin\theta_{\al}\neq 0$, diagonalize $\wt$, and  so  $2n=2m+t-s$,
and  for $A=\al+t-s$,  $\sigma_A=\theta_{\al}$  are the the non-zero 
K\"{a}hler angles of $NM$. 
That is, $TM$ and $NM$ have the same nonzero K\"{a}hler
angles, and they have the same multiplicity. Only the eigenvalues $\pm i$
of $\Fws$ and of $\wt$
may or not exist and may appear with different multiplicity, 
$t$ and $s$, respectively. Set $E_{\al}=span\{
X_{\al}, Y_{\al}\}$, $F_{A}=span\{U_A, V_{A}\}$, and $P_{E_{\al}}$, $P_{F_{A}}$
the corresponding orthonormal projections of $TM$ and $NM$.
We use the Hilbert-Schmidt inner products on tensors and forms.
We have
\begin{eqnarray}
&&
\|\Fw\|^2=\ha\|\Fws\|^2=\sm{\sum_{\al}}\cos^2\theta_{\al}=\|\wt\|^2+(s-t)=
\|\wt\|^2-2(n-m)\non\\[-1mm]
&&g(\Phi(X), \Phi(Y)) =(1-\cos\theta_{\al}\cos\theta_{\be})
g(X, Y)~~~~\mbox{for}~X\in E_{\al}, Y\in E_{\be}\\[-1mm]
&&g(\Xi(U), \Xi(V)) =(1-\cos\theta_{\al}\cos\theta_{\be})
g( U, V)~~~~\mbox{for}~ U\in F_{\al}, V\in F_{\be}\non\\[-1mm]
&&\|\Phi\|^2=2\sm{\sum_{\al}}\sin^2\theta_{\al}=\|\Xi\|^2-4(n-m)
~~~~~\non\\[-1mm]
&&\wt\circ \Phi = -\Phi\circ \Fws ~~~~~~~~~~~~
\Fws\circ \Xi = -\Xi\circ \wt\non\\[-1mm]
&&\Jt\circ \Phi = -\Phi\circ \Jw ~~~~~~~~~~~~~~~~
\Jw\circ \Xi = -\Xi\circ \Jt~~~\mbox{on}~ {\cal L}_0\\[-1mm]
&&-\Xi\circ \Phi=\sm{\sum_{\al}}\sin^2\theta_{\al}P_{E_{\al}}
~~~~~~~~~-\Phi\circ\Xi=\sm{\sum_{\al}}\sin^2\theta_{\al}P_{F_{\al}}.
\end{eqnarray}
If $X\in T_pM$ and $U\in NM_p$, then
$\omega(U,\Phi(X))=\omega(\Xi(U),X)$,
$\omega( U, \Jw X)=\omega( \Jt U,X)$, 
$g( U,\Phi(X)\rangle)=-g( \Xi(U),X)$. 

  We denote by $\lnab{}$
both Levi-Civita connections of $M$ and $N$ or $F^{-1}TN$,
if no confusion exists,
otherwise we explicit them by $\lnab{}^M$ and $\lnab{}^N$. We take on
$NM$ the usual connection $\lnabo{}$, given by $~\lnabo{X}U=
(\lnab{X}U)^{\bot}$, for $X$ and $U$ smooth sections of  $TM$ and
$NM\subset F^{-1}TN$, respectively. 
We denote the
corresponding curvature tensors by $R^M$, $R^N$ and $R^{\bot}$.
The sign convention we
choose for the curvature tensors is
$R(X,Y)Z=-\lnab{X}\lnab{Y}Z+\lnab{Y}\lnab{Z}+\lnab{[X,Y]}Z$.
The  second fundamental form of $F$, $\lnab{X}dF(Y)=\lnab{}dF(X,Y)$
is a symmetric 2-tensor on $M$ that takes values on the normal bundle.
Its covariant derivative $\lnabo{}\lnab{}dF$ is defined  considering
$\lnab{}dF$ with values on $NM$.
We denote by $i_{\ti{NM}}:NM\ra F^{-1}TN$ the inclusion bundle map,
 and  its covariant derivative $\lnab{X}i_{\ti{NM}}$  is a morphism
from $NM$ into $TM$. Then
$\forall X,Y\in T_pM$, $U\in NM_p$, $p\in M$,
\begin{equation}
\sm{g( \lnab{X}i_{\ti{NM}}(U),Y) =
g((\lnab{X}^NU)^{\top},Y)=
-g(U, \lnab{X}dF(Y))
 =-g(A^U(X),Y) =-g( U, (\lnab{X}^NY)^{\bot})}
\end{equation}
where $A:NM_p\ra L(T_pM;T_pM)$ is the shape operator.
Let $H=\frac{1}{dim(M)}trace_{g_{M}}\lnab{}dF$
denote the mean curvature of $F$.
$F$  is \em minimal \em  (resp.\ with parallel mean curvature)
if $H=0$ (resp.\ $\lnabo{}H=0$).
$F$ is  $\Jw$-pluriminimal in ${\cal L}_0$
if $(\lnab{}dF)^{(1,1)}(X,Y)=\ha(
\lnab{}dF(X,Y) +\lnab{}dF(\Jw X,\Jw Y))=0$. In this case $F$
is minimal on ${\cal L}_0$.
For  $p\in M$, $X,Y,Z\in T_{p}M$, $U,V\in NM_p$,
\begin{eqnarray}
\lnab{Z}F^{*}\omega (X,Y) &=&-g(\lnab{Z}dF(X),\Phi(Y))+
g(\lnab{Z}dF(Y),\Phi(X))\\
\lnab{Z}\wt (U,V) &=& -g(\lnab{Z}i_{\ti{NM}}(U),\Xi(V)) +
g(\lnab{Z}i_{\ti{NM}}(V),\Xi(U)).
\end{eqnarray}
If $(E,g_E)$ is a Riemannian vector bundle and $T,S: TM\ra E$ are  vector 
bundle maps, we define a  2-form  $\langle T\wedge S\rangle$ by
\[\langle T\wedge S\rangle (X,Y) =g_E( T(X), S(Y)) -
g_E(T(Y), S(X)). \]
From the symmetry of $\lnab{}dF$,
${ \langle \lnab{Z}dF\wedge \lnab{W}dF\rangle (X,Y)
= \langle \lnab{X}dF\wedge
\lnab{Y}dF \rangle (Z,W)}.$ 
Recall the Gauss, Ricci and Coddazzi equations: For $X,Y,Z\in C^{\infty}(TM)$,
and $U,V\in C^{\infty}(NM)$
\begin{eqnarray}
R^M(X,Y,Z,W) &=& R^N(X,Y,Z,W) + \langle \lnab{Z}dF\wedge  \lnab{W}dF\rangle
(X,Y)\\
R^{\bot}(X,Y,U,V)&=& R^N(X,Y,U,V) +\langle A^U\wedge A^V\rangle (X,Y)\\
-R^N(X,Y,Z,U) &=& g( ~\lnabo{X}\lnab{}dF(Y,Z)-
\lnabo{Y}\lnab{}dF(X,Z)~,~U~).
\end{eqnarray}
\subsection{$\Delta\Phi, \Delta \Fw$}
\begin{Lm}
Let $F:M\ra N$ be a $2m$-dimensional immersed  submanifold.
For any $X,Y\in T_pM$, $U,V\in NM_p$, and any local o.n.\ frame 
$e_i$ of $M$,
\\[2mm]
$(i)$~~~$\lnab{X}\Phi(Y)=\wt(\lnab{X}dF(Y))-\lnab{X}dF(\Fws(Y))$.\\[1mm]
$(ii)$~~$d\Phi(X,Y)=
-\lnab{X}dF((F^{*}\omega)^{\sharp}(Y)) +
\lnab{Y}dF((F^{*}\omega)^{\sharp}(X)).$\\[1mm]
$(iii)$~~$\delta\Phi =-2m(JH)^{\bot} $.\\[1mm]
$(iv)$\\[-14mm]
\begin{eqnarray*}
\lefteqn{~~\Delta \Phi(X) 
=  2m \, \LA{(} \lnabo{\Fws(X)} H-\lnabo{X}\wt (H)- \wt(\lnabo{X}H)\LA{)}
 +\lnab{X}dF(\delta\Fws)+}\\[-1mm]
&&+\sm{\sum_{i}}\lnab{}dF(\lnab{e_i}\Fws(X),e_i)
+\sm{\sum_i }\La{(} R^N(e_i,X, \Fws(e_i))-R^N(e_i,\Fws(X),e_i)\La{)}^{\bot}.
\end{eqnarray*}
$(v)$\\[-14mm]
\begin{eqnarray*}
\Delta \Fw (X,Y) &=&  2m \La{(} g( \lnabo{X}H,\Phi(Y))
-g(\lnabo{Y}H,\Phi(X)) \La{)} + 2m \,g( H,d\Phi(X,Y))\\
&& +Trace_M R^N(X,Y, dF(\cdot),\Phi(\cdot)) +\langle
\lnab{X}dF,\lnab{Y}\Phi\rangle -\langle
\lnab{Y}dF,\lnab{X}\Phi\rangle.
\end{eqnarray*}
$(vi)$ \\[-14mm]
\begin{eqnarray*}
g(\lnab{X}\Xi(U),Y) &=& g(\Fws \La{(}\lnab{X}i_{\ti{NM}}(U)\La{)}
-\lnab{X}i_{\ti{NM}}(\wt(U)),Y)=-g(\lnab{X}\Phi(Y), U).
\end{eqnarray*}
\end{Lm}
\noindent
\em Proof.  \em We take smooth vector fields
$X,Y$ of $M$ such that at a given point $p_{0}$, $\lnab{}Y(p_{0})=
\lnab{}X(p_{0})=0$, and assume also that $\lnab{}e_i=0$ at $p_0$.
 Then at $p_0$
\begin{eqnarray*}
\lnab{X}\Phi(Y) &=& \lnabo{X}(\Phi(Y))
=\La{(}\lnab{X}(JdF(Y))^{\bot}\La{)}^{\bot}=
\La{(}\lnab{X}(JdF(Y)-\Fws(Y))\La{)}^{\bot}\\
&=&\La{(}J\lnab{X}dF(Y)\La{)}^{\bot}-\lnab{X}dF(\Fws(Y))
=\wt(\lnab{X}dF(Y))-\lnab{X}dF(\Fws(Y))
\end{eqnarray*}
and we get $(i)$.
$(ii)$ follows from $d\Phi(X,Y)=\lnab{X}\Phi(Y)-\lnab{Y}\Phi(X)$,
and the symmetry of $\lnab{}dF$. It follows that  $
\delta\Phi =-\sum_{i} \lnab{e_{i}}\Phi (e_{i})= -2m(J H)^{\bot} +
\sum_{i}\lnab{e_{i}}dF (\Fws(e_{i})).$ 
Note that $\sum_{i}\lnab{e_{i}}dF (\Fws(e_{i}))=0$ because $\lnab{}dF(X,Y)$
is symmetric and $\Fws$ is skew symmetric. Then $(iii)$ is proved.
Now, 
\begin{eqnarray*}
d\delta\Phi(X) &=& -2m\, d((JH)^{\bot})(X)=
-2m\lnabo{X}(\wt(H))=-2m\lnabo{X}\wt (H) -2m\, \wt(\lnabo{X}H).\\[1mm]
\delta d\Phi(X)
&=& -\sm{\sum_i} \lnab{e_i}d\Phi (e_i,X)
=-\sm{\sum_i} \lnabo{e_i}(d\Phi (e_i,X) )\\
&=& -\sm{\sum_i} \lnabo{e_i}\La{(} -\lnab{e_i}dF(\Fws(X))+\lnab{X}dF(
\Fws(e_i))\La{)}\\
&=&\sm{\sum_i} \lnabo{e_i}\La{(}\lnab{\Fws(X)}dF(e_i)-\lnab{X}dF(
\Fws (e_i))\La{)}\\
&=&  {\sum_i}~~\lnabo{e_i}\lnab{}dF (\Fws(X), e_i) +
\lnab{}dF(\lnab{e_i}\Fws (X), e_i)\\[-3mm]
&&~~~~~~-\lnabo{e_i}\lnab{}dF(X, \Fws(e_i)) 
- \lnab{X}dF(\lnab{e_i}(\Fws) (e_i))\\
 &=&\sum_i ~~~\lnabo{\Fws(X)}\lnab{}dF(e_i,e_i)-
\la{(}R^N(e_i, \Fws(X)) e_i\la{)}^{\bot}
-\lnabo{X}\lnab{}dF(e_i,\Fws(e_i)) \\[-3mm]
&&~~~~~~~+\la{(}R^N(e_i, X) \Fws(e_i)\la{)}^{\bot}
+\lnab{}dF(\lnab{e_i}\Fws (X), e_i) + \lnab{X}dF(\delta (\Fws))
\end{eqnarray*}
where we applied Coddazzi's equation (2.10) in the last equality.
Since $\lnabo{X}\lnab{}dF$ is symmetric and
$\Fw$ is skew symmetric
$\sum_i \lnabo{X}\lnab{}dF(e_i, \Fws(e_i))=0$.
Thus,
\begin{eqnarray*}
\delta d\Phi(X)
&=& \lnabo{\Fws(X)} (2m H) +
\sm{\sum_i}\La{(}R^N(e_i, X) \Fws(e_i)-R^N(e_i,\Fws(X))e_i \La{)}^{\bot}\\
&&+\sm{\sum_i}\lnab{}dF(\lnab{e_i}\Fws (X), e_i) + \lnab{X}dF(\delta (\Fws)).
\end{eqnarray*}
From
$\Delta \Phi=(d\delta+\delta d )\Phi (X)$, we get the expression
in $(iv)$.
Since $\Fw$ is closed, then, for $Y$ vector field
with $\lnab{}Y(p_0)=0$, and using (2.10) and (2.6)
\begin{eqnarray*}
\lefteqn{\Delta\Fw(X,Y)=(d\delta+\delta d )\Fw (X,Y)
= d(\delta\Fw)(X,Y)=}\\
&=&\lnab{X}(\delta\Fw)(Y) -\lnab{Y}(\delta\Fw)(X)=
\sm{\sum_i} d(-\lnab{e_i}\Fw(e_i,Y))(X)
 +d(\lnab{e_i}\Fw(e_i,X))(Y)\\
&=&\sm{\sum_i} d\LA{(} g( \lnab{}dF(e_i,e_i), \Phi(Y))
-g(\lnab{}dF(e_i,Y), \Phi(e_i))\LA{)}(X)\\[-2mm]
&&+\sm{\sum_i} d\LA{(} -g( \lnab{}dF(e_i,e_i), \Phi(X))
+g(\lnab{}dF(e_i,X), \Phi(e_i))\LA{)}(X)\\[1mm]
&=&d(g( 2mH,\Phi(Y)))(X) -d(g( 2mH,\Phi(X)))(Y)
-\sm{\sum_i} g(\lnabo{X}\lnab{}dF(e_i, Y),\Phi(e_i))\\
&& -\sm{\sum_i} g(\lnab{}dF(e_i,Y),\lnab{X}\Phi(e_i))
+g(\lnabo{Y}\lnab{}dF(e_i, X),\Phi(e_i))
+g(\lnab{}dF(e_i,X),\lnab{Y}\Phi(e_i))\\[1mm]
&=&\sm{2m( g(\lnabo{X}H,\Phi(Y)) -g( \lnabo{Y}H,\Phi(X)))
+ 2m\, g( H,\lnab{X}\Phi(Y))-2m\, g(H,\lnab{Y}\Phi(X))}\\
&&+\sm{\sum_i}~g( -\lnabo{X}\lnab{}dF(Y,e_i) +\lnabo{Y}\lnab{}dF(X,e_i)~,~
\Phi(e_i))\\
&&+\sm{\sum_i}-g( \lnab{}dF(e_i,Y),\lnab{X}\Phi(e_i))
+g(\lnab{}dF(e_i,X),\lnab{Y}\Phi(e_i))\\[1mm]
&=&\sm{2m\la{(}g(\lnabo{X}H,\Phi(Y))-g(\lnabo{Y}H,\Phi(X))
+g(H, d\Phi(X,Y))\la{)}
+\sum_i R^N(X,Y,e_i,\Phi(e_i)\rangle}\\
&&+\sm{\sum_i} -g( \lnab{}dF(e_i,Y),\lnab{X}\Phi(e_i))
+g(\lnab{}dF(e_i,X),\lnab{Y}\Phi(e_i))
\end{eqnarray*}
obtaining the expression $\Delta \Fw$ of $(v)$. Finally
we prove $(vi)$.
From $\Xi(U) + \wt(U)=JU$,  and assuming at a point $p_0$~
$\lnabo{}U(p_0)=0$ ( and so $\lnab{X}^NU =\lnab{X}i_{\ti{NM}}(U)$)
 we obtain, at $p_0$
\begin{eqnarray*}
\lnab{X}\Xi(U) +\lnab{X}i_{\ti{NM}} (\wt(U)) &=&
\la{(}\lnab{X}^{N}(\Xi(U) + \wt(U))\la{)}^{\top}=(J\lnab{X}^NU)^{\top}\\[-1mm]
&=&\la{(}J\La{(}\lnab{X}i_{\ti{NM}}(U)\La{)}\la{)}^{\top}=\Fws\La{(}\lnab{X}
i_{\ti{NM}}(U)\La{)}.\\[-9mm]
\end{eqnarray*}
Therefore, using $(2.5)$,  and $(i)$,
\begin{eqnarray*}
g( \lnab{X}\Xi (U), Y) &=&
-g(\lnab{X}i_{\ti{NM}}(U), \Fws(Y))
 + g(\wt(U),\lnab{X}dF(Y))\\
&=&g( \lnab{X}dF(\Fws(Y)),U) -g( \wt(\lnab{X}dF(Y)),U)
=-g(\lnab{X}\Phi(Y),U).~~~ \qed
\end{eqnarray*}
\section{Cayley submanifolds}
\setcounter{Th}{0}
\setcounter{Pp}{0}
\setcounter{Cr} {0}
\setcounter{Lm} {0}
\setcounter{equation} {0}
We assume $m=n$ and  that $F:M^{2n}\ra N^{2n}$ has equal K\"{a}hler angles 
(e.k.a.$\scr{s}$), that is $\theta_{\al}=\theta$ $\forall \al$ and 
we denote by ${\cal L}={\cal L}_n$. In this case
$\Phi$ and $\Xi$ are conformal bundle maps and  
\begin{eqnarray}
\Fws = \cos\theta \Jw~~~~&&~~~~~~~\wt=\cos\theta\Jt\\
-\Xi\circ \Phi =
\sin^{2}\theta Id_{TM}~~~~&&~~~
-\Phi\circ \Xi =\sin^{2}\theta Id_{NM}. ~~~\mbox{on}~M.
\end{eqnarray}
On $M\sim {\cal L},$~  $\Jw$ is $g_{M}$-orthogonal. Thus $\forall \al,\be$,
\begin{equation}
g( \lnab{Z}\Jw({\al}),{\be}) =2i g(\lnab{Z}{\al},
{\be})= -g(
{\al},\lnab{Z}\Jw({\be})),~~~~~~~
g(\lnab{Z}\Jw({\al}),{\bbe}) =0.
\end{equation}
The Ricci tensor of $N$ can be expressed in terms of the
frame (2.1) as (see \cite{[S-V1]})
\begin{equation}
\sin^2\theta\, Ricci^N(U,V)= \sm{\sum_{\al}} 4R^N(U,JV, \al, \Phi(\bal))=
Trace_M R^N(U,JV,dF(\cdot), \Phi(\cdot))
\end{equation}
valid at all points $p\in M$, and $U,V\in T_{F(p)}N$.
We have
\begin{Pp} Assume $(N,J,g)$ is KE with $Ricci=R g$,
and $F:M\ra N$ is a  $2n$-dimensional immersed  submanifold with
e.k.a.$\scr{s}$. Then\\[3mm]
$(1)$~ $d\Phi(X,Y)=
2\cos\theta (\lnab{}dF)^{(1,1)}(\Jw X,Y)$ and
~ $ Trace_{\,\CC ,\Jw}d\Phi 
=\sum_{\al}d\Phi(X_{\al},Y_{\al})=2n\cos\theta ~H$.
\\[1mm]
$(2)$~ $d\sin^2\theta(X)\cdot g(Y,Z)
=g( \lnab{X}\Phi(Y),\Phi(Z))
+g(\lnab{X}\Phi(Z),\Phi(Y))$.\\[1mm]
$(3)$~ $d\Phi=0$ iff $d\Phi(X,\Jw X)=0$ iff
$F$ is $\Jw$-pluriminimal or Lagrangian.
Furthermore:
$(a)$ if $R\neq 0$, then  $d\Phi=0$ iff $F$ is complex or Lagrangian;~
$(b)$ if $R=0$, $d\Phi=0$  iff $F$ has constant K\"{a}hler angle
and $\Phi:(TM,\lnab{},g_M)\ra (NM, \lnabo{}, g)$ is a parallel
homothetic morphism. \\[2mm]
$(4)$~ $\delta \Phi =0$ iff $H$ is a Lagrangian direction of $NM$,
 iff $F$ is minimal away from ${\cal L}$.
Consequently,   $\Phi:(TM,\lnab{},g_M)\ra (NM, \lnabo{}) $ is
 closed and co-closed
1-form (and so harmonic)  iff  $\Phi:(TM,\lnab{},g_M)\ra (NM, \lnabo{})$ is
parallel iff $F$ is Lagrangian or $\Jw$- pluriminimal.\\[2mm]
$(5)$~ If $F$  has parallel mean curvature
then
\begin{eqnarray}
\Delta \Phi(X) &=& -2n\lnabo{X}\wt (H)+
\lnab{X}dF(\delta
\Fws)+\sm{\sum_{i}}\lnab{}dF(\lnab{e_i}\Fws(X),e_i)~~~~~~\\[-1mm]
&&+\sm{\sum_i} \La{(} R^N(e_i,X,\Fws(e_i))-R^N(e_i,\Fws(X),e_i)\La{)}^{\bot}\\[1mm]
\Delta \Fw (X,Y) &=& 2n\, g( ~H\, ,\, -\lnab{X}dF(\Fws (Y))+\lnab{Y}
dF(\Fws (X))~) \non\\
&&+\sin^2\theta R\Fw(X,Y) +\sm{\sum_i}
2\wt(\lnab{e_i}dF(Y),\lnab{e_i}dF(X))\non \\
&&+\langle \lnab{Y}dF, \lnab{X}dF\circ \Fws\rangle -
\langle \lnab{X}dF, \lnab{Y}dF\circ \Fws\rangle. \non
\end{eqnarray}
\end{Pp}
\noindent
\em Proof. \em $(1)$ follows from Lemma 2.1, (2) from differentiation of (2.2),
(3) and (4) are consequence of $\cite{[S-V2]}$
and Lemma 2.1 $(ii)$ and $(iii)$, and  (5) follows directly from Lemma 2.1.,
(3.4)  and the $J$-invariance of
$Ricci^N$. \qed\\[5mm]
Now  we assume $n=2$.  Four dimensional submanifolds of any 
K\"{a}hler manifold of complex dimension $\geq 4$, immersed
with equal K\"{a}hler angles,
are just the same as submanifolds  satisfying
\[ *\Fw=\pm \Fw.\]
Since pointwise  $\Fw$ is  self-dual or anti-self-dual, and is a closed
2-form, then it is co-closed as well. In particular  it is
an harmonic 2-form. This is not the case of $n\neq 2$, unless if
$\theta=constant $ (see \cite{[S-V2]}, or next lemma 3.1(2)). 
In case that $N$ is a  KE 
manifold of zero Ricci tensor and of real  dimension $8$,
a Cayley submanifold is a minimal 4-dimensional  submanifold
with equal K\"{a}hler angles $\theta_1=\theta_2=\theta$. If
$N$ is a Calabi-Yau 4-fold, that is a K\"{a}hler manifold with
a complex volume form $\rho\in \sm{\bigwedge}^{(4,0)}_{\CC}M$ (this
condition implies Ricci-flat, and the converse also  holds
in case $N$ is simply connected),
these submanifolds are characterised by being calibrated by one of  the
Cayley calibrations $\Omega=\ha\omega\wedge \omega + Re(\rho)$,
 where $\rho $ is one of the
$S^1$-family of parallel holomorphic volumes of $N$ (\cite{[H-L1]}).
Calabi-Yau 4-folds are Spin(7) manifolds. So,
locally on $N$ there is a section $\{e_1, \ldots, e_8\}$ of the
principal $Spin(7)$-bundle
of frames of $N$ defined on a open set $U$ of $N$, and that at each point
$p\in U$, defines a
isometry  of $T_pN$ onto $\RR^{8}$ such that $\Omega$ looks like
(see \cite{[J]})
\begin{equation}
\begin{array}{ccl}
\Omega&=& dx_{1234}+dx_{5678}+
\la{(}dx_{12}+dx_{34}\la{)}\wedge\la{(}dx_{56}+dx_{78}\la{)}\\
& &+\la{(}dx_{13}-dx_{24}\la{)}\wedge\la{(}dx_{57}-dx_{68}\la{)}
-\la{(}dx_{14}+dx_{23}\la{)}\wedge\la{(}dx_{58}+dx_{67}\la{)}.
\end{array}
\end{equation}
From this equation we see that
the subspaces spanned by $e_1,\ldots,e_4$ and $e_5, \ldots, e_8$
are Cayley subspaces.
We note that we use the  opposite orientation on Cayley subspaces
that Harvey and Lawson do in \cite{[H-L1]}
, and the calibration they use is  given by
$\Omega'=-\ha \omega^2+ Re(\rho)$
that is
\begin{equation}
\begin{array}{ccl}
\Omega'&=& dx_{1234}+dx_{5678}+
\la{(}dx_{12}-dx_{34}\la{)}\wedge\la{(}dx_{56}-dx_{78}\la{)}\\
&&+\la{(}dx_{13}+dx_{24}\la{)}\wedge\la{(}dx_{57}+dx_{68}\la{)}
+\la{(}dx_{14}-dx_{23}\la{)}\wedge\la{(}dx_{58}-dx_{67}\la{)}
\end{array}
\end{equation}
and so $-\Omega'$ and $\Omega$ differ on the chosen parallel holomorphic
volume, giving opposite fase on the special Lagrangian
calibration. In \cite{[H-L1]} it is proved that
$Spin(7)$ acts transitively  on the grassmannian $G(\Omega)$
of Cayley 4-planes of $\RR^{8}$ and the isotropic subgroup
of a Cayley subspace $E$ is  $K\equiv SU(2)\times SU(2)\times SU(2)/{\ZZ}_2$
(see in (1.39) of \cite{[H-L1]} how  $K$ is embedded in Spin(7)). 
Thus, we can assume that
$B=\{e_1,e_2,e_3,e_4\}$ and  $B^{\bot}=\{e_5,e_6,e_7,e_8\}$ are direct
o.n. basis of $T_pM$
and $NM_p$ respectively.
We identify isometrically in the usual way
bivectors with 2-forms. So $J_1^B=e_1\wedge e_2 + e_3\wedge e_4$,
$J_2^B=e_1\wedge e_3-e_2\wedge e_4$, and $J_3^B=e_1\wedge e_4+
e_2\wedge e_3$ defines a direct o.n. basis  (of norm $\sqrt{2}$) of
$\bigwedge^2_+T_pM$. Similar for $\bigwedge^2_+NM_p$.
We define a bilinear map:
\[\begin{array}{lc}
\Omega^{\triangle}: &\sm{\bigwedge^2}T_pM \times 
\sm{\bigwedge^2} NM_p\ra \RR^{}\\
 \blu & \Omega^{\triangle}(X\wedge Y, U\wedge V)=\Omega(X,Y,U,V)
\end{array}
\]
\begin{Pp} $\Omega^{\triangle}$ defines a natural 
orientation reversing
isometric bundle isomorphism
between $\sm{\bigwedge_+^2} TM$ and $\sm{\bigwedge_+^2} NM$.
\end{Pp}
\noindent
\em Proof. \em
Identifying isometrically in the canonic way (via musical isomorphisms
with respect to the induced metrics) the
 bilinear map $\Omega^{\triangle}$  restricted to 
$\bigwedge_+^2 T_pM\times {\bigwedge_+^2} NM_p$
with a linear map
$\Omega^{\triangle}:\bigwedge_+^2 T_pM\ra \bigwedge_+^2 NM_p
$, and  using the frame $e_i$ adapted to $M$, from (3.7)
we see that $\Omega^{\triangle}$ applies $J_1^B$ to $J_1^{B^{\bot}}$,
$J_2^B$ to $J_2^{B^{\bot}}$ and $J_3^B$ to $-J_3^{B^{\bot}}$,
 and so it gives a global orientation reversing isometry bundle map between
 the vector bundles $\bigwedge_+^2 TM$ and $ \bigwedge_+^2 NM $.\qed\\[5mm]
In particular, $p_1(\bigwedge_+^2 TM)=p_1(\bigwedge_+^2 NM)$.
We will see in section 5 that this equality still holds
for the case of $M$ Cayley but $N$ only a Ricci-flat K\"{a}hler
manifold. In this case, we cannot  guarantee the existence of
a global isomorphism between the two bundles. 
If $M$ were calibrated for $\Omega'$ we would obtain a
global orientation preserving isometry  bundle map between
 the vector bundles $\bigwedge_-^2 TM$ and $ \bigwedge_-^2 NM$.
\\[3mm]
For simplicity of notation we denote by
\begin{equation}
\gf{Z}{X}{Y}=\gdf{Z}{X}{Y}=g(\lnab{Z}dF(X),\Phi(Y))
\end{equation}
and define
$\omega_{M}(X,Y)=g_{M}(\Jw X,Y)$. We have
\begin{eqnarray}
d\omega_{M}(X,Y,Z) 
&=& g(\lnab{X}\Jw (Y), Z) - g(\lnab{Y}\Jw (X), Z)
+g( \lnab{Z}\Jw (X), Y)\\
&=&g( d\Jw (X,Y),Z) +g( \lnab{Z}\Jw (X), Y).
\end{eqnarray}
Recall the Weitzenb\"{o}ck operator of $\bigwedge^2T^*M$ applied to 
$\Fw$ is given by
\[S\Fw(X,Y)=\sm{\sum_{i}}-\overline{R}(e_i,X)\Fw(e_i,Y)
+\overline{R}(e_i,Y)\Fw(e_i,X)\]
where $e_i$ is an o.n.b.\ of $T_pM$, and 
$\overline{R}$ is the curvature operator on $\bigwedge^2T^*M$:
$\forall X,Y,u,v\in T_pM$, $\phi\in \bigwedge^2T_p^*M$, 
$(\overline{R}(X,Y)\phi)~(u,v)$ $=-\phi(R^M(X,Y)u,v)-\phi(u,R^M(X,Y)v).$ 
Let $s^M= trace~ Ricci^M=\sum_{\mu}4Ricci^M(\mu,\bmu)$
be the scalar curvature of $M$.
\begin{Lm} For an immersion with e.k.a.$\scr{s}$, $\forall X\in T_pM$, 
$\forall\al,\be$:\\[2mm]
$(1)$~~$\|\lnab{}\Fw\|^2= n\|\nabla \cos\theta\|^2 +
\ha\cos^2\theta\|\lnab{}\Jw\|^2$.\\[1mm]
$(2)$~~$\delta (\Fws) = (n-2)\Jw(\nabla \cos\theta)$.\\[1mm]
$(3)$~~$\cos\theta(\delta\Jw)=(n-1)\Jw(\nabla \cos\theta)$.\\[1mm]
$(4)$~~$\delta(\Fw)(X)=\sum_{\mu}(-2g_X\mu\bar{\mu}-2g_X\bar{\mu}\mu)
+2n g( H, JdF(X)).$\\[1mm]
$(5)$~~$g_X\be\al=g_X\al\be+\cos\theta \, g(\lnab{X}\Jw(\al),\be).$\\[1mm]
$(6)$~~$g_X\bbe\al=g_X\al\bbe +\frac{i}{2}d\cos\theta (X).$\\[1mm]
$(7)$~~$\frac{i}{2}d\cos\theta(X)=-g_{X}\be\bbe +g_X\bbe\be$
(~no summation on $\be$).\\[1mm]
$(8)$~~For $n=2$, $\delta\Fw=0$, and if $H=0$
$~\frac{i}{2}d\cos\theta(X)=\sum_{\mu}-g_X\mu\bmu=
\sum_{\mu}g_X\bmu\mu $.\\[1mm]
$(9)$~$S^M=\sum_{\mu\rho}
8 R^M(\mu,\rho,\bmu,\brho)+8 R^M(\mu,\brho,\bmu,\rho)=\sum_{\mu\rho}
16 R^M(\mu,\rho,\bmu,\brho)-8 R^M(\mu,\bmu,\rho,\brho)$.\\[1mm]
$(10)$~$\langle S\Fw,\Fw\rangle
=16\cos^2\theta\sum_{\mu\rho}R^M(\mu,\rho,\bmu,\bar{\rho})=
\cos^2\theta s^M+\sum_{\mu\rho}8\cos^2\theta
R^M(\mu,\bmu,\rho,\bar{\rho}) $.\\[1mm]
$(11)$~For $n=2$, $H=0$, and $N$ Ricci-flat,
$~\Delta\cos^2\theta =\langle
S\Fw, \Fw\rangle + \|\lnab{}\Fw\|^2$.
\end{Lm}
\noindent
\em Proof. \em  All formulas are somewhere proved in \cite{[S-V2]}. 
We only need to
check  (4) and part of (8). Since $2nH=\sum_{\mu}2\lnab{\mu}dF\bmu
+\sum_{\mu}2\lnab{\bmu}dF\mu$ and $\lnab{}dF$ is symmetric, by (2.6)
\begin{eqnarray*}
\delta\Fw(X)&=&\sm{\sum_{\mu}} -2\lnab{\mu}\Fw(\bmu,X)-2\lnab{\bmu}\Fw(\mu,X)=
\sm{\sum_{\mu}}2g_{\mu}\bmu X-2g_{\mu}X\bmu 
+2g_{\bmu}\mu X-2g_{\bmu}X\mu\\
&=& g(2nH, JdF(X))- \sm{\sum_{\mu}}(2g_{\mu}X\bmu +2g_{\bmu}X\mu).
\end{eqnarray*}
For $n=2$, from (2),
$\delta\Fw=0$ and so if $F$ is minimal by (4), $\sum_{\mu}g_{X}\mu\bmu=
-\sum_{\mu} g_{X}\bmu\mu$. Finally, by (7)
$id\cos\theta(X)= \sm{\sum_{\mu}}-g_X\mu\bmu +g_X\bmu\mu =\sm{\sum_{\mu}}
-2g_X\mu\bmu =\sm{\sum_{\mu}}2 g_X\bmu\mu.$\qed\\[5mm]
Consequently
\begin{Pp} If $n=2$ and
$F:M\ra N$ is a  submanifold with e.k.a.$\scr{s}$, then \\[2mm]
$(1)$~~$ d\omega_M=-d\log \cos\theta \wedge \omega_M,$\\[1mm]
$(2)$~~$d\cos\theta(1)=-2i \cos\theta\, g( \lnab{\bar{2}}\Jw(2),
1)= 4\cos\theta\, g( \lnab{\bar{2}}2,1)$,
\\$~~~~~~d\cos\theta(2) =-2i \cos\theta g(
\lnab{\bar{1}}\Jw(1), 2)=
4\cos\theta\, g(\lnab{\bar{1}}1,2) $,\\[2mm]
$(3)$~~and if $F$ is a Cayley submanifold, then
~$g_X2\bar{2}=-g_X\bar{1}1.$
\end{Pp}
\noindent
\em Proof. \em (1) From $0=d\Fw=d(\cos\theta\omega_M)
=d\cos\theta\wedge \omega_M +\cos\theta d\omega_M $ we obtain (1).
From (1), 
$\frac{i}{2}d\cos\theta(1) =d\cos\theta\wedge \omega_M(1,2,\bar{2})
=-\cos\theta d\omega_M(1,2, \bar{2}).$
Equations  (3.3) and (3.10) give
$ d\omega_M(1,2, \bar{2})= g( \lnab{\bar{2}}\Jw(1),2 )=
-g( \lnab{\bar{2}}\Jw(2),1)=-2ig(\lnab{\bar{2}}2,1)$
obtaining the first equality in (2). Similar for $d\cos\theta(2)$.
Lemma 3.1(7)(8), and  minimality of $F$ imply
$-2g_X1\bar{1}+2g_X\bar{1}1=id\cos\theta(X)
=\sum_{\mu}-2g_X\mu\bmu=-2g_X1\bar{1}-2g_X2\bar{2}$
and we get (3).\qed
\subsection{The almost complex structure  $\Jw$}
\begin{Pp} If $n=2$, and $F$  has e.k.a.$\scr{s}$, on $~M\sim{\cal L}$, 
we have
\begin{eqnarray*}
(1)&&d\Jw=0\mbox{~iff~}  \lnab{}\Jw=0.
\mbox{~~In~this~case~}d\omega_M=0.\\
(2)&& \cos\theta = const. ~~\mbox{iff}~~ d\omega_M=0 ~
~ \mbox{iff}~~ \lnab{\bar{\ga}}\Jw(\al)=0~~\mbox{iff}~~\delta\Jw=0.\\
(3)&&\mbox{If}~\lnab{}\Jw=0\mbox{~then~}\cos\theta = constant\mbox{
~and~}\Fw\mbox{~parallel}.\\
(4)&&\Jw\mbox{~is~integrable}~~~\mbox{iff}~~~d\Jw(1,2)=0
~~~\mbox{iff}~~~\lnab{\ga}\Jw(\al)=0.\\
(5)&&\mbox{If~} F \mbox{~is}~\Jw\mbox{-pluriminimal}~ \mbox{then} 
~\cos\theta =~constant.~~
\end{eqnarray*}
Consequently we have:\\[2mm]
$(A)$ If  $\Jw$ is integrable, then  $\cos\theta=constant$ iff
$\Jw$ is K\"{a}hler.\\[1mm]
$(B)$  If $\Jw$ is an almost complex structure of the Gray list \cite{[G1]}
that is, K\"{a}hler, or almost, quasi, nearly or semi-K\"{a}hler, then
$\cos\theta$ is constant. Therefore, immersions  with non-constant 
e.k.a.$\scr{s}$
produces submanifolds with an almost complex structure on $M\sim {\cal L}$, 
that might be integrable  but not one of the Gray list. 
\end{Pp}
\noindent
\em Proof. \em  From (3.10)  $d\Jw=0$ implies\\[-4mm]
\begin{equation}
d\omega_M(X,Y,Z)=g(\lnab{Z}\Jw(X),Y)
\end{equation}
with $(Z,X)\ra \lnab{Z}\Jw(X)$ symmetric. But then $d\omega_M=0$.
By (3.12) we get $\lnab{}\Jw=0$, and we have proved (1).
Since $\lnab{}\Jw( T^{1,0}M)\subset T^{0,1}M$ (see (3.3)), then, ~~
$ d\omega_{M}(\al,\be,\bga)=g( \lnab{\bga}\Jw(\al),\be),$
and for $n=2$, $\al,\be,\ga$ must have repeated vectors, so, 
$0 =d\omega_{M}(\al,\be,\ga)=g( d\Jw(\al,\be),\ga) +
g( \lnab{\ga}\Jw(\al), \be)$.
Using Prop.3.3(2) and (3.3)
\begin{eqnarray}
\|d\omega_{M}\|^{2}&=&\sum_{\ga, \al<\be} 16 | d\omega_{M}(\al,\be,\bga)|^2
=\sum_{\ga, \al<\be}16|\langle \lnab{\bga}\Jw(\al),\be\rangle|^{2}=
\sum_{\ga,\al}8\|\lnab{\bga}\Jw(\al)\|^2\non \\[-1mm]
&=&\frac{4}{\cos^2\theta}(|d\cos\theta(1)|^2+|d\cos\theta(2)|^2)=
2\|\nabla\log(\cos\theta)\|^2
\end{eqnarray}
and we have proved the first 3 equivalences of (2). The  last one comes
from lemma 3.1(3).
If $\lnab{}\Jw=0$ then $\delta \Jw=0$. Thus, by (2) $\cos\theta $
is constant, and by Lemma 3.1(1) $\lnab{}\Fw=0$ and (3) is proved.
The integrability
of $\Jw$ is  equivalently to the vanishing
of the tensor:
$ N_{\Jw}(X,Y)=[\Jw X,\Jw Y]-[X,Y] -\Jw [X,\Jw Y] -\Jw [\Jw X, Y].$
Using the connection on $M$ we have
$N_{\Jw}(X,Y)=-\Jw(d\Jw (X,Y)) + (\lnab{\Jw X}\Jw) (Y)
-(\lnab{\Jw Y}\Jw)(X).$
From $\Jw\circ \Jw=-Id$ we have $\lnab{X}\Jw(\Jw Y)=-\Jw(\lnab{X}\Jw(Y))$.
Thus,
$N_{\Jw}=0$ iff $d\Jw (X,Y)=d\Jw (\Jw X, \Jw Y)$, iff
$d\Jw (\al,\be)=0$, iff $\lnab{\ga}\Jw(\al)=0$.
Recall  that pluriminimality implies $\cos\theta= constant$ (\cite{[S-V2]}).
Now we prove (A). If $\Jw$ is integrable and $\cos\theta$ is constant,
by (2) and (4) $\lnab{}\Jw=0$, i.e. $\Jw$ is K\"{a}hler. 
Now we prove the last remark (B). If $\Jw$ is almost-K\"{a}hler, that is
$d\omega_M=0$, from (2) $\cos\theta=constant$.
If $\Jw$ is nearly-K\"{a}hler, that is $\lnab{X}\Jw(X)=0$, then $\delta
\Jw=0$, and so it is semi-K\"{a}hler. The later implies by (2) that
$\theta$ is constant.
If $\Jw$ is quasi-K\"{a}hler, then
$\lnab{X}\Jw(Y)=-\lnab{\Jw X}\Jw(\Jw Y)$, and so $\lnab{1}\Jw(\bar{1})=
\lnab{2}\Jw(\bar{2})=0$.
By Proposition 3.3 this implies $\cos\theta=constant$.
$\qed$\\[5mm]
\em Remark. \em 
As an observation, we conclude that if
 $R\neq 0$ and $M$  is compact immersed with e.k.as and with  almost
K\"{a}hler $\Jw$
(i.e. $d\omega_M=0$), then $M$ is K\"{a}hler, confirming the 
Goldberg conjecture in this case. In fact, from previous proposition
we have $\theta$ constant, and so theorem 1.2 of \cite{[S1]} concludes
that $\Jw$ is  K\"{a}hler.
\\[3mm]
Let $B=\{X_1,Y_1,X_2,Y_2\}$ be a diagonalising o.n. basis of $\Fw$
at the point $p\in M$, and let $\epsilon(p)\in \{-1,+1\}$
be the sign of this basis (well defined for $p\notin {\cal L}$),
and if $p\in M\sim  {\cal C}$,
we denote by   $\epsilon'(p)$ the sign of the basis $B'=\Phi(B)$ of $NM_p$.
Then $\frac{\Phi(Y_1)}{\sin\theta}, \frac{\Phi(X_1)}{\sin\theta},
\frac{\Phi(Y_2)}{\sin\theta},\frac{\Phi(X_2)}{\sin\theta}$ diagonalizes
$\wt$, and has the same orientation as $B'$.
\begin{Lm} Let $n=2$ and $F:M\ra N$ be an immersed oriented  submanifold
with e.k.a.$\scr{s}$. Then $\Phi:TM\ra NM$ is, along 
$M\sim {\cal L}\cup {\cal C}$, an orientation
preserving bundle morphism, with respect to the
orientations defined by $\Jw$ and $\Jt$ respectively. Moreover,
for $p\in M\sim  {\cal C}$, $\epsilon(p)\epsilon'(p)=+1$
always hold.
\end{Lm}
\noindent
\em Proof. \em From (2.3)
$\Phi$ is a $\Jw-\Jt$-anti-holomorphic morphism, and
so it preserves the orientations. If $p\in {\cal L}$, for any
chosen basis $B$, $\epsilon(p)\epsilon'(p)=+1$ still holds. Now assume
$p\notin {\cal L}\cup {\cal C}$. Since
$\Phi(X_{\al})=JX_{\al}-\cos\theta Y_{\al}$,
$\Phi(Y_{\al})=JY_{\al}+\cos\theta X_{\al}$,
then
\begin{eqnarray*}
\lefteqn{0 \neq A_{\Phi}(T_pM):=
\mbox{Vol}_{N}\left(X_1,Y_1,X_2,Y_2,\Phi(X_1),
\Phi(Y_1), \Phi(X_2),\Phi(Y_2)\right)=\epsilon(p)\epsilon'(p)\sin^4\theta}\\
&=&\sm{\mbox{Vol}_{N}(X_1,Y_1,X_2,Y_2,JX_{1},JY_{1},
JX_{2}, JY_{2})=
\mbox{Vol}_N(X_1,JX_1,Y_1,JY_1,X_2, JX_{2},Y_2, JY_{2}).}
\end{eqnarray*}
Note  that if $T_pM$ is a Lagrangian subspace, then $A_{\Phi}=1$
for any o.n.b.\ $X_1,Y_1, X_2, Y_2$ of $T_pM$ we choose.
Now we prove that $A_{\Phi}>0$ holds for any $T_pM$ with e.k.a.$\scr{s}$ with
$\cos\theta\neq 1$. We consider the strictly decreasing continuous
curve,  $a:[0,\frac{\pi}{2}]\ra [0,1]$, defined  by
$ a(t)= \frac{1- \sin(t)}{\cos(t)}$, $a(0)=1$, and $a( \frac{\pi}{2})=0 $.
We may identify $T_pN$ with $V=\RR^{4}\times \RR^{4}$, with
complex structure
$J_0(X,Y)= (-Y,X)$.
We consider the family of maps, for $t\in [0,\frac{\pi}{2}]$,
$~ \Gamma_t: \RR^{4}\ra V,~\Gamma_t(X)= (X, a(\frac{\pi}{2}- t)\Jw(X))$, 
where $\Jw$ is a fixed $g_0$-orthogonal complex structure
on $\RR^{4}$. For each $t$, $\Gamma_t$
is an isomorphism of  $\RR^{4}$
onto a subspace $E_t=\Gamma_t(\RR^{4})$ with e.k.a
$\theta=(\frac{\pi}{2}-t) $
and complex structure $\Jw$ (see \cite{[S1]}),
giving a continuous curve of subspaces starting
from a Lagrangian subspace $E_0=\RR^{4}$ of $V$.
Then $t\in [0, \frac{\pi}{2}[\ra A_{\Phi_t}(E_t)$
is a continuous curve of
nonzero numbers, starting with value 1 at $t=0$. Thus  it remains positive
on $[0, \frac{\pi}{2}[$.
Now, all subspaces $E$ of dimension 4 of $V$ that have the same K\"{a}hler
angles
are the same up to a unitary transformation of $V$, ( such
transformation maps  a diagonalising basis of $E$ and of $E^{\bot}$
into the corresponding ones of $E'$ and $E'^{\bot}$, see \cite{[S1]} 
or \cite{[T1]}).
This proves
$\epsilon (p)\epsilon' (p)=\frac{A_{\Phi}(T_pM)}{\sin^4\theta}>0$.\qed\\[4mm]
Recall that a hyper-K\"{a}hler manifold is a Riemannian manifold $(N,g)$
endowed with two K\"{a}hler structures $I,J$ that anti-commute, $IJ=-JI$.
Then $I,J,K:=IJ$ define a family of K\"{a}hler structures indexed on
$S^2$,  $(J_x)_{x\in S^2}$, $J_x=aI+bJ+cK$, $x=(a,b,c)$.
\begin{Pp} Let $n=2$, and
 $F:M\ra N$ be an immersed connected oriented submanifold
with e.k.as. Then we have:\\[2mm]
$(1)$~ With respect to the  given orientation of $M$,
 $\Phi: TM\ra NM$ is, away from ${\cal C}$,
an orientation preserving morphism. Furthermore, 
we may assume that the orientation of $M$ is
such that $\Fw$ is self-dual on all $M$, and so  $\Jw$ and $\Jt$ define
the orientation of $M$ and $NM$ resp..\\[1mm] 
$(2)$~ If $(N,(J_x)_{x\in S^2}, g)$ is hyper-K\"{a}hler, and $F:M\ra N$
is a $J_x$-complex submanifold, then $\forall y\in S^2$,
$\Fw_y$ is selfdual , where $\omega_y$ is the K\"{a}hler form of
$(N,J_y,g)$. Moreover $M$ is a Cayley submanifold of $(N,J_y,g)$ with
k.a $\cos\theta_y(p)=\|(J_yX)^{\top}\|$ where $X\in T_pM$ is any
unit vector. 
\end{Pp}
\noindent
\em Proof. \em (1) follows immediately from Lemma 3.2.
Each of the 2-forms  $\eta_{\pm}=\Fw  \pm *\Fw $
is  harmonic and so, if not identically to zero, its zero set
has empty interior. 
This implies that one of the $\eta_{\pm}$ must vanish
identically. Thus, we may choose the orientation of $M$ s.t.
 $\Fw$ is self-dual  on all $M$, and so 
 $\Jw$ defines the orientation of $M$.
From Lemma 3.2 $\Jt$ defines the orientation of $NM$.
 Now we prove (2).
First we recall  that
for any $x,y\in S^2$
$J_xJ_y=J_{x\cdot y}=-\langle x,y\rangle Id + J_{x\times y}$.
Let
$B=\{e_1,e_2,e_3,e_4\}=\{X,J_xX,Z, J_xZ\}$ be an o.n.b. of $T_pM$.
We have~
$g( J_yJ_xX,X) =-\langle y,x\rangle =g( J_yJ_xZ,Z)$
, $g(J_yJ_xZ,X)
=-g( J_yX, J_xZ)$,
$ g(J_yZ, J_xX)=-g( J_{y\times x}X, Z ),$
and  $g( J_yJ_xZ,J_xX)
=g(J_yX, Z)$.
A basis for the self-dual 2-forms on $M$ is given by\\[-3mm]
\begin{equation}
J^B_i=e_1\wedge e_2+ e_3\wedge e_4,~~~~~J^B_j=e_1\wedge e_3- e_2\wedge e_4,
~~~~~J^B_k=e_1\wedge e_4+ e_2\wedge e_3.
\end{equation}
Then we see that $ \Fw_y=\cos\theta_y(p)J^B_u$
where
$ u=\frac{1}{t}(<y,x>i+g( J_yX,Z)j+g( J_{y\times x}X,Z)k)$
($i,j,k$ is the usual basis of $\in \RR^{3}$)
 and\\[-3mm]
\[\cos\theta_y(p)=t=\sqrt{(<y,x>^2+g( J_yX,Z)^2
+g(J_{y\times x}X,Z)^2)}=\|(J_yX)^{\top}\|,\]
proving that $\Fw_y$ is self-dual. In particular $M$ is
a Cayley submanifold (see also  \cite{[S2]}).
\qed
\subsection{A particular case}
Let us first assume that $N$ is an hyper-K\"{a}hler (HK) manifold
$(N,I,J,K,g)$ of complex dimension 4, 
where $I$ and $J$ are $g$-orthogonal K\"{a}hler
structures on $N$ that anti-commute and  $K=IJ$.  
If $M$ is an $I$-complex submanifold of complex dimension 2, 
then from Prop.3.5 $M$ is a 
$J$-Cayley submanifold of $N$. Let $\omega_J$ be the 
K\"{a}hler form of $(N,J,g)$.
Then ${{\Fw}_J}^{\sharp}=\cos\theta_J J_{\omega_J}$.
Similar for $K$.
We are going to describe an o.n. basis
that diagonalizes ${\Fw}_J$ and ${\omega_J^{\bot}}$ that we use in 
\cite{[A-PV-S]} and \cite{[S-PV]}.  
Let $p\in M$ and $X\in T_pM$ be a unit vector and $H_X=
span \{X,IX,JX,KX\}$. Since $T_pN$ is a vector space of dimension 8
and $X\in T_pM$, there exist
 $U\in H_X^{\bot}\cap NM_p$ unit vector
and  o.n. basis $B$ 
of $T_pM$ and $B^{\bot}$ of $NM_p$ of the form
\[\begin{array}{l}
B=\{W_1,W_2=IW_1,W_3,W_4=IW_3\}
=\{X,IX,J(cX+sU),K(cX+sU)\}\\
B^{\bot}=\{U_1,U_2=IU_1,U_3,U_4= IU_3\}
=\{-U,-IU,J(cU-sX), K(cU-sX)\}
\end{array}\]
where $c^2+s^2=1$.
The basis $\{W_1,J_{\omega_J}W_1, W_2,J_{\omega_J}W_2\}$ and
$\{W_1,J_{\omega_K}W_1, W_2,J_{\omega_K}W_2\}$ 
diagonalize ${\Fw}_J$ and  ${\Fw}_K$ respectively 
with $\cos\theta_J=|c|=\cos\theta_K$. Moreover, 
$I=J_i^B$, $J_{{\omega}_J}=\epsilon J^B_j$,
$J_{{\omega}_K}=\epsilon J^B_k$, 
where $\epsilon= sign~c$ (if $c=0$ take $\epsilon =1$) (see (3.14)).
Consequently, the hyper-Hermitian structure on $TM$ defined by the orientation
determined by $I$
is given by $\{ I, J_{{\omega}_J}, J_{{\omega}_K} \}$, defined on ${\cal L}_0$,
that is, away from totally complex points, i.e.\ points with
$c\!=\!0$, (see next remark).\\[4mm] 
\em Remark. \em  
A convenient multiple of the $m$-power of the 
fundamental 4-form $\Omega$ gives  a calibration (\cite{[T2]}). Then we may
define a quaternionic angle for any $4m$-dimensional submanifold.
  The  $\Omega$-angle  of a complex 4-submanifold 
(in the sense of quaternionic-Khaler geometry) is given by
$\theta(p)$ s.t. $\cos\theta(p)=\Omega(X_1,X_2,X_3,X_4)$ where
$X_i$ is a d.o.n.b. of $T_pM$, and  has values
between $\frac{1}{3}$ and $1$. The first extreme value corresponds to
a \em totally complex point, \em  that is a point s.t. $T_pM$  is
$I$-complex and $J$-Lagrangian, for a local almost hyper-Hermitian structure
$I,J,K=IJ$ of $N$. The other extreme value corresponds to a
\em quaternionic point, \em  a point s.t. $T_pM$ is a quaternionic
subspace of $T_pN$, or equivalently $T_pM$  is
$I$ and $J$ -complex.\\

Now we return to the general case of 
$(N,J,g)$ being a Ricci-flat KE 4-fold and $F:M\ra N$ an immersed
4-submanifold with e.k.a.$\scr{s}$.
For each non $J$-Lagrangian point $p$ and local almost
complex structure $I$ on $M$ orthogonal to $\Jw$ we can
find a local basis  ${Z_{\al}}_{1\leq \al\leq 2}$  defined by (2.1) satisfying
\begin{equation}
I(1)=\bar{2}~~~~~I(2)=-\bar{1}~~~~~I(\bar{1})=2~~~~~~I(\bar{2})=-1. 
\end{equation}
In this subsection we are going to assume that  for each 
$p_0\in M\!\sim\! {\cal L}$  a K\"{a}hler complex structure $I$
orthogonal to $\Jw$ exists on a  open dense set $O$ of a 
neighbourhood of $p_0$, and we will say that
$M$ is \em  $I$-K\"{a}hler. \em  
Then (see e.g.\ \cite{[E-S]}, \cite{[LeB]}) we have the following orthogonal 
decomposition w.r.t. $I$ on $O$:
\begin{equation}
\sm{\bigwedge^2_{+}}TM=\RR{}\omega_I\oplus (T^{(2,0)}M\oplus T^{(0,2)}M)
\end{equation}
where $\omega_I$
denotes the $I$-K\"{a}hler form on $O$.
Moreover any real self-dual harmonic 2-form $\zeta$ orthogonal to 
$\omega_I$ (for the Hilber-Schmidt inner product) is of the form
$\zeta=\varphi + \bar{\varphi}$ where $\varphi$ is a holomorphic 
$(2,0)$-form over 
$O$. The $I$-Ricci form $\rho^M$, $\rho^M(X,Y)=Ricci^M(IX,Y)$ 
and the scalar curvature $s^M$ of $M$ is given by
\begin{eqnarray}
\rho^M &=& -i \partial\bar{\partial}\log \|\varphi\|^2=
 -i \partial\bar{\partial}\log \|\zeta \|^2\\
s^M &=& \Delta^+ \log \|\varphi\|^2=\Delta^+ \log \|\zeta \|^2
\end{eqnarray}
away from the zero set of $\zeta$. 
If we take $\zeta={\Fw}_J$ we conclude:\\[-4mm]
\begin{Th}  on a open set  $O$ of $M$ where $M$ is $I$-K\"{a}hler
\begin{eqnarray}
Ricci^M(IX,Y)&=&-\ha d(d\log\cos^2\theta \circ I)(X,Y)\\
s^M &=& \Delta\log\cos^2\theta.
\end{eqnarray}
\end{Th}
\begin{Pp} If $M$ is $I$-K\"{a}hler,
then  $\Jw$ is integrable iff $\theta$ is constant.
In the particular case that $N$ is hyper-K\"{a}hler, and
$M$ is $I$-complex submanifold, then $\Jw$ is integrable iff
$\theta$ is constant, iff
$\Jw$ is K\"{a}hler, iff $(M,I, {\Jw}, I{\Jw})$ is hyper-K\"{a}hler.
\end{Pp}
\noindent
\em proof. \em 
Since $\Jw$ and $I$ anti-commute and
$I$ is parallel on $M$, then
$(\lnab{}\Jw) \circ I=\lnab{}(\Jw \circ I) =-\lnab{}(I\circ \Jw )=
-I\circ \lnab{}\Jw.$
Therefore,
$g( \lnab{Z}J_{\omega} (\bar{1}), \bar{2}) =
g( \lnab{Z}J_{\omega} (-I(2)), I(1)) =
g(I\lnab{Z}J_{\omega} (2), I(1))=
g(\lnab{Z}J_{\omega} (2), 1)$.
By Prop.3.3(2) and (3.15)
$d\cos\theta(1)$ $=2i\cos\theta\, 
g( \lnab{I(1)}\Jw (1),2)$, $
d\cos\theta(2)$ $=2i\cos\theta\, 
g( \lnab{I(2)}\Jw (1),2)$. Thus,
 $d\cos\theta(Z)= 
2i \cos\theta g( \lnab{IZ}\Jw (1), 2)$, or equivalently
\begin{equation}
\lnab{Z}\Jw(Y)=-d\log\cos\theta(IZ)\Jw I(Y).
\end{equation}
From Prop.3.4(4) and the above formula $(3.21)$ we conclude that 
$\Jw$ is integrable iff
$0=g( d\Jw(1,2),1) =\ha d\log\cos\theta(1)$
and $0=g( d\Jw(1,2),2) =-\ha d\log\cos\theta(2)$,
that is, iff $\cos\theta$ is constant. This condition turns out
to be equivalent to $\Jw$ to be K\"{a}hler, by Prop.3.4(A). \qed.
\section{Complex and Lagrangian points}
\setcounter{Th}{0}
\setcounter{Pp}{0}
\setcounter{Cr} {0}
\setcounter{Lm} {0}
\setcounter{equation} {0}
In this section we will study the nature of 
the complex and the Lagrangian points 
of a submanifold  $F:M^{2n}\ra N^{2n}$ immersed with
parallel mean curvature and with e.k.as.
We  introduce some natural complex vector subbundles ${\cal F}^{+}$
and ${\cal F}^{-}$ of $F^{-1}TN$, over $M\sim {\cal L}$,
that generalizes to higher dimensions  the
special complex vector subbundles defined for immersed real surfaces
into K\"{a}hler surfaces given in \cite{[E-G-T]}.
Namely, for each  point $p\in M\sim {\cal L}$, we define the  $J$-complex
vector subspaces of $T_{F(p)}N$
\begin{equation}
{\cal F}^{+}_p= \{X-J\Jw X : X\in T_{p}M\},~~~~~~~~~~
{\cal F}^{-}_p =({\cal F}^{+})^{\bot}
\end{equation}
and the linear morphisms, over $M\sim {\cal L}$, $\Psi^{\pm}:
TM\ra {\cal F}^{\pm}$, $\Psi^{\pm}(X)=\ha(X\pm J\Jw X)$,
being $\Psi^{+}$ a complex  morphism,
while $\Psi^{-}$ is an anti-complex one, both conformal:
\[
 \Psi^{\pm}\circ \Jw = \pm J\circ \Psi^{\pm}, 
~~~~~~~~~~g( \Psi^{\pm}(X), \Psi^{\pm}(Y)) =
\sm{\frac{(1\pm\cos\theta)}{2}}
g(X,Y)~~~~~~ \forall X, Y.
\]
In particular $\Psi^+$ is an isomorphism  over $M\sim{\cal L}$, 
what implies
${\cal F}^{+}$ to be smooth of real  rank $2n$. Thus, the same holds for
${\cal F}^{-}$. 
Denoting the decompositions
$TM^c=T^{1,0}M\oplus T^{0,1}M$  and $TN^c=T^{1,0}N\oplus T^{0,1}N$,
with respect to $\Jw$ and  $J$ respectively,  we have
$\Psi^{+}:T^{1,0}M\ra T^{1,0}N$, $\Psi^{-}:T^{1,0}M\ra T^{0,1}N$.
At $p\in M\sim {\cal L}$,  $\forall X\in T_{p}M$~
$X=\Psi^{+}(X)+\Psi^{-}(X)$, $J\Jw X= \Psi^{-}(X)-\Psi^{+}(X)$.
Note that, w.r.t the complex structure $J$,
$T_p^{1,0}N=({\cal F}^+_p)^{1,0}\oplus ({\cal F}^-_p)^{1,0}$. Then
we may take a local unitary o.n. frame
$(\sqrt{2}W_{\al}, \sqrt{2}K_{\al})_{1\leq
\al\leq n}$ of $T^{(1,0)}N$, along $M\sim {\cal L}$ s.t.
\begin{eqnarray}
W_{\al} \in T^{1,0}N\cap ({\cal F}^{+})^{c}
~~~~~~~~~~
K_{\al} \in T^{1,0}N\cap
({\cal F}^{-})^{c}.
\end{eqnarray}
Let $p$ a non Lagrangian point, and $X_{\al},Y_{\al}=\Jw X_{\al}$ a
diagonalising o.n. local frame of $\Fw$, on a neighbourhood of $p$ and
let $Z_{\al}=\al$ as in $(2.1)$. 
Note that $Z_{\al}\in T^{1,0}M$ with respect to $\Jw$. Define
 some local complex maps $u_{\al\be}$, $v_{\al\be}$ on 
$M\sim {\cal L}$ by
\begin{eqnarray}
\Psi^+(\al)=(Z_{\al})^{1,0}=\sm{\sum_{\be}}u_{\al\be}W_{\be},~~~~~~~~~~~
\Psi^-(\al)=(Z_{\al})^{0,1}=\sm{\sum_{\be}}v_{\al\be}K_{\bbe}.~~~~~~~~~~~~~
\end{eqnarray}
Consider the $n\times n$ complex matrices
$u=[u_{\al\be}]_{1\leq \al, \be \leq n},$
~$ v= [v_{\al\be}]_{1\leq \al, \be \leq n}$. 
\begin{Lm} 
$u\cdot \bar{u}^t=\bar{u}^t\cdot u =  
 \ha (1+\cos\theta)Id,$ ~
$v\cdot \bar{v}^t=\bar{v}^t\cdot v =
\ha (1-\cos\theta)Id$.
In particular, for each $\al,\mu$, 
$|u_{\al\mu}|^{2}\leq \frac{1+\cos\theta}{2}$, 
$|v_{\al\mu}|^{2}\leq \frac{1-\cos\theta}{2}$.
\end{Lm}
\noindent
\em Proof. \em We have
$\sum_{\ga}\ha u_{\al\ga}\overline{u_{\be\ga}}=\sum_{\ga\rho}g(
u_{\al\ga}W_{\ga}, \overline{{u}_{\be\rho}}W_{\bar{\rho}})=
g( \Psi^{+}(\al), \overline{\Psi^{+}(\be)})=
\ha(1+\cos\theta)\frac{\delta_{\al\be}}{2},$
and similar for $v$.
Recall  that
for matrices,
$A\bar{A}^t=D$, where $D$ is a real diagonal matrix,
 implies $\bar{A}^t A= D$.
\qed\\[3mm]
Now we obtain some estimates:
\begin{Lm} On a neighbourhood of a point $p\in  M\sim {\cal L}$,
there exists a constant $C>0$ s.t.
\begin{eqnarray}
\forall \be,\mu \mbox{~and~}\forall A,B\in C^{\infty}(T^cN)
 ~~~~~|R^N(\be,\mu,A, B)| &\leq &  C \|\Phi\|~~~~~~~~~~~~~~~~\\
\|\lnab{}\Jw\| &\leq &   C\|\Phi\|.~~~~~~~~~~~~~~~~
\end{eqnarray}
\end{Lm}
\noindent
\em Proof. \em 
Since $(N,J,g)$ is K\"{a}hler,  by (4.3)
\[
R^N(\be,\mu,A,B)=
\sm{\sum_{\al,\rho}}u_{\be\al}v_{\mu\rho}
R^{N}(W_{\al},K_{\brho},A,B)+
v_{\be\al}u_{\mu\rho}
R^{N}(K_{\al},W_{\rho},A,B)
\]
Thus, the estimate of $|v_{\mu\rho}|$ in Lemma 4.1 and that
$\|\Phi\|^2=4\sin^2\theta=4(1-\cos\theta)(1+\cos\theta)$ implies (4.4). 
By (3.3), to estimate $\|\lnab{}\Jw\|$ we need only to estimate
$|\langle \lnab{Z}\Jw (\mu), \rho\rangle|$. We have
\[\lnab{Z} \Fw (\mu,\rho) = \langle \lnab{Z} \Fws (\mu),\rho\rangle
= \langle d\cos\theta(Z)\Jw(\mu) + \cos\theta\lnab{Z}\Jw (\mu),
\rho\rangle
=  \cos\theta \langle \lnab{Z}\Jw(\mu),\rho\rangle\]
and from (2.6) we obtain (4.5).
\qed \\[1mm]
\begin{Pp} Assume  $F$  with parallel mean curvature and e.k.as.\\[2mm]
$(1)$ Locally  there exist a constant $C>0$ such that
$\|\Delta \Phi\|\leq C\|\Phi\|.$\\[1mm]
$(2)$ If $N$ is KE, locally  
there exist a constant $C>0$ such that
$\|\Delta \Fw\|\leq  C\|\Fw\|.$\\[-4mm]
\end{Pp}
\noindent
\em Proof. \em 
Note first that
by  the expression of $\lnab{Z}\Fw$ in (2.6), we have
$ \|\lnab{}\Fw\|, ~\|\delta \Fw\|\leq C\|\Phi\|.$
The term $\lnabo{X}\wt (H)$ of $\Delta \Phi (X)$ in Prop.3.1(5)
can be estimate using (2.7):
$\|\lnabo{X}\wt(H)\|\leq C \|\Xi\|\leq C \|\Phi\|.$
Now let $ L(X)=\sum_i \La{(} R^N(e_i,X,\Fws(e_i)) -
R^N(e_i,\Fws(X), e_i)~\La{)}^{\bot}.$
$L$ vanish at Lagrangian points. On $M\sim {\cal L}$ we take
a local unitary complex frame $\sqrt{2}U_{\al}, \sqrt{2}U_{\bal}$
of the complexified normal bundle. Then
\[
\sm{L(\mu)}=
\sm{\sum_{\al} -4i\cos\theta (R^N(\al,\mu, \bal))^{\bot}
= \sum_{\al, \ga}-8i\cos\theta R^N(\al,\mu, \bal, U_{\ga})U_{\bar{\gamma}}
-8i\cos\theta R^N(\al,\mu, \bal, U_{\bar{\ga}})U_{\ga}.}
\]
Therefore, by Lemma 4.2 and Lemma 4.1,
we conclude that
$ \|L(\mu)\|\leq C\|\Phi\|$, and 
we have proved that
$\|L(X)\|\leq C\|\Phi\|.$ 
By  Prop.3.1(5), away from ${\cal L}$,
 $\|\Delta\Phi\|\leq  C\|\Phi\|$. 
At a point $p\in {\cal L}$, $\Phi$ is an isometry and is smooth
so the inequality also holds. Again by  Prop.3.1(5),
$\|\Delta \Fw\|\leq C (\|\Fw\|+\|\wt\|)\leq C |\cos\theta|\leq C\|\Fw\|.$
\qed \\[5mm]
In order to conclude from  Proposition 4.1 that complex points and 
Lagrangian points 
are zeros of finite order of $\Phi$ and $\Fw$ respectively, 
we need to translate some  inequalities of 
Aronszajn-type for vector bundle maps to similar inequalities for 
the components.
If $\psi$ is a $r$-form on $M$ with values on a Riemannian vector bundle
$E$ over $M$, the Weitzenb\"{o}ck formula reads
\[\Delta \psi=(d\delta +\delta d)\psi =\sm{\sum_i}
-\lnab{e_i}(\lnab{e_i}\psi) +\lnab{\lnab{e_i}e_i}
\psi_A  ~~+ S(\psi)
\]
where $e_i$ is any o.n.\ frame of $M$ and $S(\psi)$ is the
Weitzenb\"{o}ck operator on $\bigwedge^{r}T^*M\otimes E$.
Assume $M$ is  a connected Riemannian manifold, and $E_A$ is a finite family
of Riemannian vector bundles over $M$, and $\forall A$, 
$\psi_A\in C^{\infty}(\bigwedge^{r_A}T^*M\otimes E_A)$ 
is $E_A$-valued $r_A$-form on $M$.
We need the following Aronszajn-type theorem:
\begin{Lm} Assume $M$ is connected and  there is a constant $C>0$ s.t.
$\|\Delta \psi_A\|\leq {\sum_{B}} C\La{(} \|\psi_B\| 
+\|\lnab{}\psi_B\|\La{)}~\forall A.$
If $\{\psi_A\}$ have a common zero of infinite order, then all $\psi_A
\equiv 0$ on all $M$.
\end{Lm}
\noindent
\em Proof. \em 
Let $e_i$ and $w_{A,\al}$ be a local o.n. frames of $M$ and of $E_A$
respectively. For each $\sigma =\{i_1<\ldots < i_{r_A}\}$,
 $e^{\sigma}_{A,\al}=e_*^{i_1}\wedge \ldots \wedge e_*^{i_{r}}\otimes
w_{A,\al}$,  defines an o.n. frame of $\bigwedge^{r_A}T^*M\otimes
E_A$ for the Hilbert-Schmidt inner product. Let $a^{\al}_{A,\sigma}$
be the local components of $\psi_A$ w.r.t  $e^{\sigma}_{A,\al}$,
$\psi_A = \sum_{\sigma, \al} a^{\al}_{A,\sigma}e^{\sigma}_{A,\al}$.
Then $\lnab{X}\psi_A=\sum_{\al,\sigma}
da^{\al}_{A,\sigma}(X)e^{\sigma}_{A,\al}+
a^{\al}_{A,\sigma}\lnab{X}e^{\sigma}_{A,\al}$ and applying
Weitzenb\"{o}ck formula to $\psi_A$,
\[\Delta a^{\al}_{A,\sigma} = -\langle \Delta \psi_A,
e^{\sigma}_{A,\al}\rangle -\sm{\sum_{i \be \rho}} 2da^{\be}_{A,\rho}(e_i)
\langle\lnab{e_i}e^{\rho}_{A,\be}, e^{\sigma}_{A,\al}\rangle
+\sm{\sum_{\be \rho}}a^{\be}_{A,\rho}\langle\Delta e^{\rho}_{A,\be},
e^{\sigma}_{A,\al}\rangle.\]
Consequently,  there exists  constant
$C',C>0$ s.t.
$|\Delta a^{\al}_{A,\sigma}|\leq  \|\Delta \psi_A\| +\sm{\sum_{\be,\rho}} C'
(|a^{\be}_{A,\rho}|+\|\nabla a^{\be}_{A,\rho}\|)
$ $\leq \sm{\sum_{B,\be \rho}}C\La{(}
|a^{\be}_{B,\rho}| +\|\nabla a^{\be}_{B,\rho}\|\La{)}.$
A common zero of infinite order of $\{\psi_A\}$ is a common zero of 
infinite order of the family $\{a^{\al}_{B,\sigma}\}$, and the lemma 
follows from last remark of \cite{[Ar]} .~~\qed.\\[-3mm]

\begin{Pp} Let $m$ be the dimension of  $M$  and  $Z$ be the
 set of common zeros of  $\psi_A$. If for
each $p\in Z$ there exist an $A$ s.t. $p$ is a zero of finite order
of $\psi_A$, then $Z$ is a countably $(m-1)$-rectifiable set, and in
particular has Hausdorff codimension at least 1. 
\end{Pp}
\noindent
\em Proof. \em See, for example,  a proof in \cite{[Ba]}.\qed\\[3mm]
Since the zeros of infinite order of $\Phi$ (resp. $\Fw$) are  zeros of 
infinite order of $\sin^2\theta$ (resp. $\cos^2\theta$), and that we call
by complex points (resp. Lagrangian points) of infinite order,
 the above estimates in Prop.4.1, Lemma 4.3 
and Prop.4.2 leads to the conclusion (1) and (2) below:
\begin{Cr} Let $F:M^{2n}\!\ra\! N^{2n}$ be  an immersion with   
parallel mean  curvature and e.k.as. Then:\\[1mm]
$(1)$ If $H=0$,  and $F$ is not a complex
submanifold,  the set ${\cal C}$ of complex points is a 
set of $M$ of Hausdorff codimension at least 1. \\[1mm]
$(2)$ If $N$ is KE and $F$ is not a 
Lagrangian submanifold, the Lagrangian points is a
set of $M$ of Hausdorff codimension at least 1.\\[1mm] 
$(3)$ If $n=2$, $M$ is closed, and $F$ is any immersion with e.k.as, the set 
${\cal L}$ of Lagrangian points  is a countably $(n-2)$-rectifiable
set and so has Hausdorff codimension at least 2.
\end{Cr}
\noindent
\em Proof. \em  (3) If $n=2$ and $F$ has e.k.as then $\Fw$ 
 satisfies $D\Fw=0$ where $D=d+\delta$ is the
usual Dirac operator for forms on $M$, from \cite{[Ba]} we obtain the
result.\qed .
\\[4mm]
\em Remark. \em In case (3) $\Fw$ is an harmonic self-dual $2$-form.
It is known that the zero set of generic harmonic self-dual $2$-form,
is a disjoint union of curves diffeomorphic to $S^1$.
\section{A residue-type formula}
\setcounter{Th}{0}
\setcounter{Pp}{0}
\setcounter{Cr} {0}
\setcounter{Lm} {0}
\setcounter{equation} {0}
\subsection{Curvature tensors and characteristic classes}
Recall that if  $(E, g,\lnabe{})$ is a rank $k$ Riemannian vector bundle
over  a Riemannian manifold $M$ of dimension $4$, 
the   first Pontrjagin class  $p_1(E)$ can
be represented in the cohomology class  $H^{4}(M,\RR^{})$ 
 by the 4-form 
 defined by using the curvature tensor $R^E$ of $\lnab{}^E$ 
\begin{eqnarray*}
p_1(E)={p}_1(R^E)&=&\sm{\sum_{i<j}}
\frac{1}{4\pi^2} R^E_{ij}\wedge R^E_{ij}
\end{eqnarray*}
where the curvature components $R^E_{ij}\in C^{\infty}(\bigwedge^2T^*M)$
are defined w.r.t. a local o.n. frame $B=\{E_i\}_{1\leq i\leq k}$.
If $k=4$ and $E$ is oriented, 
the  Euler characteristic class ${\cal X}(E)$ is given by
\begin{eqnarray*}
{\cal X}(E)={\cal X}(R^E)
&=&\frac{1}{4\pi^2}(R^E_{12}\wedge R^E_{34}-R^E_{13}\wedge R^E_{24}
+R^E_{14}\wedge R^E_{23}).
\end{eqnarray*}
If we take $\{1,2,3\}=\Lambda^{\pm}_1, \Lambda^{\pm}_2, \Lambda^{\pm}_3$
the usual corresponding basis of $ \bigwedge^2_{\pm}E$ built from
$B$ (see (3.14) for the selfdual case) we easily verify that
$R^{\wedge^2_{\pm}}_{12}=2R^E_{\Lambda^{\pm}_3}$,
$R^{\wedge^2_{\pm}}_{13}=-2R^E_{\Lambda^{\pm}_2}$,
$R^{\wedge^2_{\pm}}_{23}=2R^E_{\Lambda^{\pm}_1}$.
With this basis one derives   the  well known relation:
\begin{equation}
p_1(\sm{\bigwedge^2_{\pm}}E)=p_1(E)\pm 2{\cal X}(E).
\end{equation}
Set, for   direct orthonormal bases $e_i$ of $T_pM$ and $E_i$ of $E_p$,
\[
z_1=\ha(e_1-ie_2),~~ z_2=\ha(e_3-ie_4),~~
w_1=\ha(E_1-iE_2),~~  w_2=\ha(E_3-iE_4)
\]
Then $\mbox{Vol}_M(z_1,z_{\bar{1}},z_2,z_{\bar{2}})
=-\frac{1}{4}$, and if $R_{AB}^E$ denotes the curvature components
w.r.t. this basis, i.e. with  $A,B\in \{w_1,w_{\bar{1}}, w_2,w_{\bar{2}}\}$, 
we have
\begin{eqnarray}
{\cal X}(E)
&=&\frac{1}{\pi^2}( R^E_{12}\wedge R^E_{\bar{1}\bar{2}}
-R^E_{1\bar{1}}\wedge R^E_{2\bar{2}}+R^E_{1\bar{2}}\wedge R^E_{2\bar{1}})
\\
p_1(E) &=&\frac{1}{\pi^2}(-R^E_{1\bar{1}}\wedge R^E_{1\bar{1}}-
R^E_{2\bar{2}}\wedge R^E_{2\bar{2}}+2R^E_{12}\wedge R^E_{\bar{1}\bar{2}}
-2R^E_{1\bar{2}}\wedge R^E_{2\bar{1}})\\
p_1(\sm{\bigwedge^2_{+}}E)&=&- 
\frac{1}{\pi^2}\LA{(} (R^E_{1\bar{1}}+R^E_{2\bar{2}})
\wedge (R^E_{1\bar{1}}+R^E_{2\bar{2}})-
4 R^E_{12}\wedge R^E_{\bar{1}\bar{2}}
\LA{)}\\
p_1(\sm{\bigwedge^2_{-}}E)&=&-
\frac{1}{\pi^2}\LA{(} (R^E_{1\bar{1}}-R^E_{2\bar{2}})
\wedge (R^E_{1\bar{1}}-R^E_{2\bar{2}})+
4 R^E_{1\bar 2}\wedge R^E_{2\bar{1}}\LA{)}
\end{eqnarray}
Herman  Weyl introduced some curvature invariants $\kappa_{2c}(M)$, $1\leq
c\leq[\scr{\frac{n}{2}}]$, of a manifold $M$ of dimension $n$ embedded
in a Euclidean space, that appear in its formula on the volume of a tube 
of radius $r$ about $M$.
These invariants are defined in the same way for any Riemannian manifold $M$
(see e.g. \cite{[G2]}). For $c=1$, and $c=2$ they are respectively
\begin{eqnarray*}
\kappa_2(M) = \ha\int_{M}s^M \mbox{Vol}_M, ~~~~~&&
\kappa_4(M)=\frac{1}{8}
\int_{M}({(s^M)}^2 - 4 \|Ricci^M\|^2 + \|R^M\|^2)\mbox{Vol}_M\non
\end{eqnarray*}
where $\|R^M\|$
is the Hilbert-Schmidt norm of $R^M$ as a  4-tensor.
Thus, for  $dim(M)=4$, $\kappa_{4}$ reads the Gauss-Bonnet formula
$\frac{1}{4\pi^2}\kappa_4(M)={\cal X}(M)$ (see e.g \cite{[Be]}).
\\[3mm]
If $E$ and $F$ are vector bundles over $M$
and $T:TM\ra E$ a 1-tensor,
$l:E\times F\ra \RR_M$ and  $R:TM\times TM \ra F$ 2-tensors, then
$l(T \wedge R)\in \bigotimes^3TM^*$ denotes the 3-tensor
\[\sm{l(T \wedge R)(X,Y,Z)=\cyclic_{X,Y,Z}l(T(X),R(Y,Z))=
l(T(X),R(Y,Z))+l(T(Z),R(X,Y))+l(T(Y),R(Z,X)).}\]
If $R$ is symmetric (resp. skew symmetric), then so it is $l(T \wedge R)$.
 We also recall  the Kulkarni-Nomizu
operator, a symmetric product for two 2-tensors $\phi,\xi\in \bigotimes^2
TM^*$
\[\phi\bullet \xi (X,Y,Z,W)=\phi(X,Z)\xi(Y,W)+\phi(Y,W)\xi(X,Z)
-\phi(X,W)\xi(Y,Z)-\phi(Y,Z)\xi(X,W)\]
Assume $(E, g_E)$ is a Riemannian vector bundle with a Riemannian
connection $\lnabe{}$.
The curvature tensor $\bar{R}$ of $TM^*\otimes E$
is given by
\[(\bar{R}(X,Y)\Phi)(Z)=-\lnab{X,Y}^2\Phi(Z)+\lnab{Y,X}^2\Phi(Z)=
R^E(X,Y)(\Phi(Z))-\Phi(R^M(X,Y)Z)\]
where $\lnab{X,Y}^2\Phi=\lnab{X}(\lnab{Y}\Phi)
-\lnab{\lnab{X}Y}\Phi$, for a smooth section $\Phi$ of $TM^*\otimes E$.
\begin{Lm} Let $(E,\lnabe{},g_E)$ be a rank-4 Riemannian vector bundle over
$M$ and $\Phi:TM\ra E$ a conformal morphism, with $g_E(\Phi(X),
\Phi(Y)) = h g(X,Y)$, and denote by
$\Xi=-h\Phi^{-1}$. Then:\\[2mm]
$(1)$~$g( \Xi(U),Y)=-g_E( U,\Phi(Y)),$ and $g(\lnab{X}\Xi(U),Y)
=-g( \lnab{X}\Phi(Y),U).$\\[1mm]
$(2)$~$g_E( \lnab{X}\Phi(Y),\Phi(Z)) + g_E( \Phi(Y),\lnab{X}
\Phi(Z)) =dh(X)g(Y,Z)$.\\[1mm]
$(3)$~$ d^2\Phi(X,Y,Z)=
-(\bar{R}(X,Y)\Phi)(Z) -(\bar{R}(Z,X)\Phi)(Y) -(\bar{R}(Y,Z)\Phi)(X)=
-R^E\wedge \Phi(X,Y,Z).$\\[1mm]
$(4)$~$g_E((\bar{R}(X,Y)\Phi)(Z),\Phi(W))=-
g_E((\bar{R}(X,Y)\Phi)(W),\Phi(Z)).$\\[1mm]
$(5)$~$R^E(X,Y,\Phi(Z),\Phi(W))
= hR^M(X,Y,Z,W)+g_E((\bar{R}(X,Y)\Phi)(Z),\Phi(W))$.
\end{Lm}
\noindent
\em Proof. \em Using $\Phi\circ \Xi=-h Id_E$, 
(1) and (5)  are obvious. (2) is obtained from differentiation of
$g_E(\Phi(Y),\Phi(Z))=hg(Y,Z)$.
Since $\lnabe{}$ is a $g_E$-Riemannian connection, from (5) we derive (4).
(3) follows from the definitions  $d\Phi(X,Y)=\lnab{X}\Phi(Y)
-\lnab{Y}\Phi(X)$ and
$d^2\Phi(X,Y,Z)=\sm{\cyclic_{\ti{X,Y,Z}}}
(\lnab{X}d\Phi) (Y,Z)$,
and that $R^M$ satisfies first Bianchi.
 \qed \\[6mm]
We consider the degenerated metric on $M$,
$\hat{g}(X,Y)=g_E(\Phi(X),\Phi(Y))$, and singular
connection $\lnabp{}=\Phi^{-1*}\lnab{}$ with torsion $T$
that makes $\Phi:(TM, \lnabp{}, \hat{g})\ra (E, \lnabe{}, g_E)$ parallel,
 namely $\lnabp{X}Y=\lnab{X}Y+ S(X,Y)$, where
\[
S(X,Y)=\Phi^{-1}\lnab{X}
\Phi (Y)\mbox{~~~and~~~}T(X,Y)=\Phi^{-1}d\Phi(X,Y).
\]
It is a Riemannian connection w.r.t
$\hat{g}$. Since $\Phi$ is conformal
then $\hat{g}=hg$. Let $\hlnab{}$ denote the Levi-Civita connection
of $(M,\hat{g})$, $\varphi=\log h$, and set
\[\hat{S}(X,Y)=\hlnab{X}Y-\lnab{X}Y=\ha \varphi_XY+\ha\varphi_YX-\ha
g(X,Y)\nabla \varphi
,~~~S'(X,Y)=\lnabp{X}Y-\hlnab{X}Y.\]
where $\varphi_X=d\varphi(X)$.
Then $S=\hat{S}+S'$, and $T(X,Y)=S(X,Y)-S(Y,X)=S'(X,Y)-S'(Y,X)$.
The curvature tensor $R':\bigwedge^2TM \ra \bigwedge^2TM$
 of $\lnabp{}$, that is given by $\Phi(R^E)$,
i.e.
\begin{eqnarray}
R'(X,Y,Z,W) &=&
\hat{g}(R'(X,Y)Z,W)=g_E(R^E(X,Y)\Phi(Z), \Phi(W))\non\\
&=& h g(\Phi^{-1}R^E(X,Y)\Phi(Z),W)=\hat{g}( \Phi(R^E)(X,Y)Z,W)
\end{eqnarray}
may not be a curvature-type tensor.
The   Bianchi map for $R: \bigwedge^2TM\ra \bigwedge^2TM$  is defined as
$g(b(R)(X,Y,Z),W)=\!\!\!\cyclic_{\ti{X,Y,Z}}R(X,Y,Z,W).$
 Note that $b(R)\in \bigwedge^3TM^*\otimes TM\,\bigcap\,
L(\bigwedge^2TM;\bigwedge^2TM)$. 
\begin{Pp} In the conditions of previous lemma, $b(R')=-\Phi^{-1}d^2\Phi$.
So $R'$
satisfies the first Bianchi identity iff  $(\bar{R}(X,Y)\Phi)(Z)$
does so, iff $d^2\Phi=0$. In that case $R'$
is also symmetric. Thus, $R'=\Phi(R^E)$ is a curvature operator at a point
$p\in M$ (i.e $R'\in {\cal B}$, see notation in \cite{[S-PV]})
 iff $d^2\Phi(p)=0$.
\end{Pp}
\noindent
\em Proof. \em  Using the fact that $R^M$ satisfies the first Bianchi 
identity  and lemma 5.1(3) gives
\begin{eqnarray*}
b(R')(X,Y,Z,W) &=& \cyclic_{\ti{X,Y,Z}}
R^{E}(X,Y,\Phi(Z),\Phi(W))=\cyclic_{\ti{X,Y,Z}}
g_E((\bar{R}(X,Y)\Phi)(Z), \Phi(W))\\[-1mm]
&=&
 -g_E(d^2\Phi(X,Y,Z),\Phi(W))=-\hat{g}(\Phi^{-1}d^2\Phi(X,Y,Z), W).
\end{eqnarray*}
Now we have from symmetry of $R^M$,
\begin{eqnarray*}
R'(Z,W,X,Y)-R'(X,Y,Z,W) &=&
g((\bar{R}(X,Y)\Phi)(Z),\Phi(W))
-g((\bar{R}(Z,W)\Phi)(X),\Phi(Y))\\
&=&-g_E( d^2\Phi(X,Y,Z),\Phi(W))
+g_E( d^2\Phi(X,Y,W),\Phi(Z))\\
&&-g_E( d^2\Phi(Y,Z,W),\Phi(X))
+g_E(d^2\Phi(X,Z,W),\Phi(Y)).\qed\\[-2mm]
\end{eqnarray*}
\begin{Pp} If $(Y,Z)\ra g(\lnab{X}\Phi(Y), \Phi(Z))$ is symmetric
then $\bar{R}(X,Y)\Phi=0$.
\end{Pp}
\noindent
\em Proof. \em  Set $\lh=\log h$.
From lemma 5.1(2) we have $\lnab{X}\Phi(Y)=\ha d\lh(X)\Phi(Y)$.
It follows that
$\lnab{X,Y}^2\Phi(Z)=\ha Hess \lh (X,Y)\Phi(Z)+ \frac{1}{4}\lh_X\lh_Y\Phi(Z)
.$
That implies $\lnab{X,Y}^2\Phi=\lnab{Y,X}^2\Phi$.\qed
\subsection{Proof of (1.5) of Theorem 1.1}
If $N$ is
KE with $Ricci^N=Rg$,  (3.4) says that
\begin{equation}
\sm{ \sum_{\al}}R^N(\al,\Phi(\bal))=\frac{\sin^2\theta}{4}R\,\omega
 = \sm{\sum_{\al}}R^N(\bal,\Phi(\al))
\end{equation}
Note that $JX=(JX)^{\top}+(JX)^{\bot}=
\cos\theta \Jw(X)+\Phi(X)$, and so
\begin{eqnarray*}
R^N(X,Y, \Phi(Z), \Phi(W)) &=&
R^N(X,Y,Z,W) +\cos^2\theta\, R^N(X,Y, \Jw Z, \Jw W)\\
&&+\cos\theta\, R^N(X,Y,Z,J\Jw W) + \cos\theta\, R^N(X,Y,J\Jw Z,W)
\non
\end{eqnarray*}
 Since $\Jw(\al)=i\al$,
we have on $M\sim {\cal L}$, ~
\[\begin{array}{l}
R^N(\Phi(\al),\Phi(\be))=\sin^2\theta\,R^N(\al,\be),~~~~~~
R^N(\Phi(\al),\bar{\ga})=-R^N(\al,\Phi(\bar{\ga})),\\
R^N(\Phi(\al),\ga)=-R^N(\al,\Phi(\ga)) -2i\cos\theta\, R^N(\al,\ga),\\
R^N(\Phi(\al),\Phi(\bar{\ga}))=(1+\cos^2\theta)R^N(\al,\bar{\ga})-
2i\cos\theta R^N(\al,J\bar{\ga})
=\sin^2\theta \,R^N(\al,\bar{\ga})-
2i\cos\theta\,R^N(\al,\Phi(\bar{\ga})).\end{array}\]
\begin{equation}
\sm{\sum_{\al}}R^N(\Phi(\al),\Phi(\bal))=
\sm{\sum_{\al}}\sin^2\theta R^N(\al,\bal)
-2i\cos\theta R^N(\al,\Phi(\bal))
= \sin^2\theta\La{(} \sm{\sum_{\al}}R^N(\al,\bal)-\frac{i}{2}\cos\theta\,R\,
\omega\La{)}.
\end{equation}
The Gauss and the Ricci equations (2.8)-(2.9) gives
\begin{Lm} On $M\sim {\cal L}$
\begin{eqnarray}
\!\!\!\!
R^{\bot}(\Phi(\al),\Phi(\be))&=&\sin^2\theta\, R^M(\al, \be)-\sin^2\theta
\,\langle \lnab{\al}dF\wedge \lnab{\be}dF\rangle
+\langle A^{\Phi(\al)}\wedge A^{\Phi(\be)}\rangle~~~~~~~~~~~\\[3mm]
R^{\bot}(\Phi(\al),\Phi(\bbe))&=&\sin^2\theta\, R^M(\al, \bbe)
-2i\cos\theta R^N(\al,\Phi(\bbe))\\
&&~-\sin^2\theta
\,\langle \lnab{\al}dF\wedge \lnab{\bbe}dF\rangle
+\langle A^{\Phi(\al)}\wedge A^{\Phi(\bbe)}\rangle. \non\\\non
\end{eqnarray}
\end{Lm}
\noindent
We have for $A,B\in T^c_pM$
\begin{eqnarray*}
\lefteqn{\langle \lnab{A}dF\wedge \lnab{B}dF\rangle (X,Y)-
\frac{1}{\sin^2\theta}\langle A^{\Phi(A)}\wedge A^{\Phi(B)}\rangle(X,Y)=}\\
&=&\frac{2}{\sin^2\theta}\sum_{\al}\La{(}~(g_XA\al\, g_YB\bal
-g_X\al A\, g_Y\bal B) +
(g_XA\bal\, g_YB\al -g_X\bal A g_Y \al B)\\[-3mm]
&&~~~~~~~~~~~~~~
-(g_YA\al\, g_XB\bal-g_Y\al A\,g_X\bal B)
-( g_YA\bal g_XB\al -g_Y\bal A g_X\al B)~\La{)}\\[-4mm]
\end{eqnarray*}
Using lemma 3.1 applied to the above equation we have
\begin{Lm} For a Cayley submanifold $F:M\ra N$, we have on
$M\sim ({\cal L}\cup {\cal C})$,
\begin{eqnarray*}
\sm{\ha\langle A^{\Phi(1)}\wedge A^{\Phi(\bar{1})}\rangle} &=&
\sm{\frac{\sin^2\theta}{2}\langle \lnab{1}dF\wedge \lnab{\bar{1}}dF\rangle
 +\frac{i}{2}d\cos\theta\wedge
 (g_{~\cdot~}\bar{1}1+g_{~\cdot~}1\bar{1})
}\non\\
&& \sm{+\cos\theta g_{~\cdot~}21\wedge g( \lnab{\cdot}
\Jw(\bar{1}),\bar{2})-\cos\theta g_{~\cdot~}\bar{1}\bar{2}\wedge
 g( \lnab{\cdot}\Jw({1}),{2})}\non 
\end{eqnarray*}
\begin{eqnarray*}
\sm{\ha\langle A^{\Phi(2)}\wedge A^{\Phi(\bar{2})}\rangle} &=&
\sm{\frac{\sin^2\theta}{2}\langle \lnab{2}dF\wedge \lnab{\bar{2}}dF\rangle
 +\frac{i}{2}d\cos\theta\wedge
 (g_{~\cdot~}\bar{2}2+g_{~\cdot~}2\bar{2})}\non\\
&& \sm{+\cos\theta g_{~\cdot~}12\wedge g( \lnab{\cdot}
\Jw(\bar{2}),\bar{1})-\cos\theta g_{\cdot}\bar{2}\bar{1}\wedge
 g( \lnab{\cdot}\Jw({2}),{1})}\non \\[3mm]
\sm{\ha\langle A^{\Phi(1)}\wedge A^{\Phi({2})}\rangle} &=&
\sm{\frac{\sin^2\theta}{2}\langle \lnab{1}dF\wedge \lnab{{2}}dF\rangle}\\[3mm]
\sm{\ha\langle A^{\Phi(1)}\wedge A^{\Phi(\bar{2})}\rangle} &=&
\sm{\frac{\sin^2\theta}{2}\langle \lnab{2}dF\wedge \lnab{\bar{2}}dF\rangle
+id\cos\theta \wedge g_{~\cdot ~}1\bar{2}}\\
&&\sm{-\cos\theta g_{~\cdot~}11\wedge g( \lnab{\cdot}\Jw(\bar{1}),
\bar{2})
-\cos\theta g_{~\cdot~}\bar{2}\bar{2}\wedge
 g( \lnab{\cdot}\Jw({1}),{2}) }
\end{eqnarray*}
\end{Lm}
\noindent
\em Proof. \em This is a long but straightforward proof using
lemma 3.1(5)(6)(8) and (3.3).
We only prove one of the equalities, for the other ones
are similar. 
\begin{eqnarray*}
\lefteqn{\La{(}\frac{\sin^2\theta}{2}\langle
\lnab{1}dF\wedge\lnab{2} dF\rangle-
\ha \langle A^{\Phi(1)}\wedge A^{\Phi(2)}\rangle\La{)}(X,Y)=}\\
&=&\sm{ \sum_{\al} (g_X1\al g_Y2\bal - g_X\al 1 g_Y\bal 2)
+(g_X1\bal g_Y2\al - g_X\bal 1 g_Y\al 2)}\\ 
&&\sm{+ \sum_{\al} (-g_Y1\al g_X2\bal + g_Y\al 1 g_X\bal 2)
+(-g_Y1\bal g_X2\al + g_Y\bal 1 g_X\al 2)}
\end{eqnarray*}
\begin{eqnarray*}
&=& \sm{\sum_{\al} g_X1\al g_Y2\bal -( g_X1\al +\cos\theta g(
\lnab{X}\Jw(1),\al) )(g_Y2\bal
+\frac{i}{2}\delta_{\al 2} d\cos\theta(Y))}\\
&&\sm{~~~~+g_X1\bal g_Y2\al -( g_X1\bal +\frac{i}{2}
\delta_{\al 1} d\cos\theta(X))
(g_Y2\al +\cos\theta g(\lnab{Y}\Jw(2),\al))}\\
&&\sm{ \sum_{\al} -g_Y1\al g_X2\bal +( g_Y1\al +\cos\theta g(
\lnab{Y}\Jw(1),\al) )(g_X2\bal
+\frac{i}{2}\delta_{\al 2} d\cos\theta(X))}\\
&&\sm{~~~~-g_Y1\bal g_X2\al +( g_Y1\bal +\frac{i}{2}
\delta_{\al 1} d\cos\theta(Y))
(g_X2\al +\cos\theta g(\lnab{X}\Jw(2),\al))}\\
&=&\sm{ -\frac{i}{2}g_X12d\cos\theta(Y)
-\cos\theta g( \lnab{X} \Jw (1),2) g_Y2\bar{2}
-\cos\theta g_X1\bar{1} g(\lnab{Y}\Jw(2),1)
-\frac{i}{2}d\cos\theta(X)g_Y21 }\non\\
&&\sm{ +\frac{i}{2}g_Y12d\cos\theta(X)
+\cos\theta g( \lnab{Y} \Jw (1),2) g_X2\bar{2}
+\cos\theta g_Y1\bar{1}g( \lnab{X}\Jw(2),1)
+\frac{i}{2}d\cos\theta(Y)g_X21 }\non\\
&&\sm{-\frac{i}{2}\cos\theta d\cos\theta (Y)g( \lnab{X}\Jw(1),2)
-\frac{i}{2}\cos\theta d\cos\theta (X)g(\lnab{Y}\Jw(2),1)}\non\\
&&\sm{+\frac{i}{2}\cos\theta d\cos\theta (X)g( \lnab{Y}\Jw(1),2)
+\frac{i}{2}\cos\theta d\cos\theta (Y)g(
\lnab{X}\Jw(2),1)}\non\\
&=&\sm{ \frac{i}{2} d\cos\theta(X)(-g_Y21 +g_Y12)
-\frac{i}{2} d\cos\theta(Y)(-g_X21 +g_X12)}\\
&&\sm{-\cos\theta g(\lnab{X}\Jw(1),2) g_Y2\bar{2}
-\cos\theta g(\lnab{Y}\Jw(2),1) g_X1\bar{1}}\\
&&\sm{ +\cos\theta g(\lnab{Y}\Jw(1),2) g_X2\bar{2}
+\cos\theta g(\lnab{X}\Jw(2),1) g_Y1\bar{1}}\\
&&\sm{+i\cos\theta d\cos\theta(X) g( \lnab{Y}\Jw(1),2)
-i\cos\theta d\cos\theta(Y) g(\lnab{X}\Jw(1),2) }
\end{eqnarray*}
\begin{eqnarray*}
&=&\sm{ \frac{i}{2}\cos\theta d\cos\theta(X) g(\lnab{Y}\Jw(2),1)
-\frac{i}{2}\cos\theta d\cos\theta(Y))g(\lnab{X}\Jw(2),1)}\\
&&\sm{-\cos\theta g(\lnab{X}\Jw(1),2)(g_Y2\bar{2}+g_Y1\bar{1})
+\cos\theta g(\lnab{Y}\Jw(1),2 )(g_X2\bar{2}+g_X1\bar{1})}\\
&&\sm{+i\cos\theta d\cos\theta(X) g( \lnab{Y}\Jw(1),2)
-i\cos\theta d\cos\theta(Y)g(\lnab{X}\Jw(1),2)}\\
&=&\sm{ -\frac{i}{2}\cos\theta d\cos\theta(Y)g(\lnab{X}\Jw(1),2)
+\frac{i}{2}\cos\theta d\cos\theta(X) g(\lnab{Y}\Jw(1),2)}\\
&&\sm{+\cos\theta g( \lnab{X}\Jw(1),2)(\frac{i}{2}d\cos\theta(Y))
+\cos\theta g(\lnab{Y}\Jw(1),2)(-\frac{i}{2}d\cos\theta(X))}\\
&=&0~~~~~~~~~~\qed.\\[-2mm]
\end{eqnarray*}
The two previous lemmas, (5.7) and (5.8) give us
\begin{Pp} If $F:M\ra N$ is a Cayley submanifold of
a Ricci-flat  $N$, then on $M\sim {\cal L}$
\begin{eqnarray*}
\sm{R^{\bot}(\Phi(1),\Phi(\bar{1}))} &=&
\sm{\sin^2\theta R^M(1,\bar{1})+2id\cos\theta\wedge
g_{~\cdot~}1\bar{1}-2i\cos\theta R^N(1,\Phi(\bar{1}))}\\
&& \sm{+2\cos\theta g_{~\cdot~}21\wedge g( \lnab{\cdot}
\Jw(\bar{1}),\bar{2}) -2\cos\theta g_{~\cdot~}\bar{1}\bar{2}\wedge
 g( \lnab{\cdot}\Jw({1}),{2})}\non \\
\sm{ R^{\bot}(\Phi(2),\Phi(\bar{2}))} &=&
\sm{\sin^2\theta R^{M}(2,\bar{2}) -2id\cos\theta\wedge g_{~\cdot~}1\bar{1}
-2i\cos\theta R^N(2,\Phi(\bar{2}))}\\&& \sm{+2\cos\theta g_{~\cdot~}12\wedge g(\lnab{~\cdot~}
\Jw(\bar{2}),\bar{1})-2\cos\theta g_{~\cdot~}\bar{2}\bar{1}\wedge
g(\lnab{\cdot}\Jw({2}),{1})}\non 
\end{eqnarray*}
\begin{eqnarray*}
\sm{R^{\bot}(\Phi(1),\Phi(2))} &=& \sin^2\theta R^M(1,2)\\[3mm]
\sm{R^{\bot}(\Phi(1),\Phi(\bar{2}))\rangle} &=&
\sm{\sin^2\theta R^M(1,\bar{2})
+2id\cos\theta \wedge g_{~\cdot ~}1\bar{2}-2i\cos\theta R^N(1,\Phi(\bar{2}))}\\
&&\sm{-2\cos\theta g_{~\cdot~}11\wedge g( \lnab{\cdot}\Jw(\bar{1}),
\bar{2})
-2\cos\theta g_{~\cdot~}\bar{2}\bar{2}\wedge
g( \lnab{\cdot}\Jw({1}),{2}) }
\end{eqnarray*}
Furthermore,
$\sum_{\al}R^{\bot}(\Phi(\al), \Phi(\bal))= \sum_{\al}\sin^2\theta
R^M(\al,\bal)$.\\[-3mm]
\end{Pp}
\noindent
\begin{Pp} If $F:M\ra N$ is a non-$J$-holomorphic Cayley submanifold
and $N$ is Ricci-flat, then $(1.5)$ holds, that is
$p_1(\bigwedge^2_+NM)=p_1(\bigwedge^2_+TM)$.
\end{Pp}
\noindent
\em Proof. \em 
If $M$ is a Lagrangian submanifold,
then $\Phi:TM\ra NM$ is an orientation preserving isometry,
and so characteristic classes of $M$ and $NM$ are the same.
Now we assume
$M$ is neither Lagrangian nor complex submanifold. 
We consider  the formulas (5.2)-(5.4) using the
curvature tensors $R^M$ of $M$ and $R^{\bot}$ of $NM$ w.r.t the 
connections $\lnab{}^M$ of $M$
and  $\lnabo{}$ of $NM$, respectively, and away from complex and Lagrangian 
points we may take
$e_1,e_2,e_3,e_4$ as $X_1,Y_1,X_2,Y_2$ and $E_1,E_2,E_3, E_4$ as
$\frac{\Phi(X_1)}{\sin\theta}, \frac{\Phi(Y_1)}{\sin\theta},
\frac{\Phi(X_2)}{\sin\theta}, \frac{\Phi(Y_2)}{\sin\theta}$.
By the previous Proposition 5.3  and (5.4) we easily see that the equality
(1.5) is valid on $M\sim {\cal L}\cup {\cal C}$, as forms
( and not only as chomology classes). Moreover the expressions in (5.2)-(5.5)
do not depend on the o.n.\ basis used, and are smoothly defined on all $M$.
Since the set of complex
and Lagrangian points have empty interior (corollary 4.1).
Then (5.4) and so (1.5) stays
valid on all $M$. \qed.\\[4mm]
From (5.1) and  the previous proposition we obtain:
\begin{Cr} In the conditions of the Prop.5.4, $~{\cal X}(M)-
{\cal X}(NM)=\ha(p_1(NM)-p_1(M))$.
\end{Cr}
\subsection{Proof. of (1.6) of Theorem 1.1}
Since $\lnabp{}$ is a $\hat{g}$-Riemannian connection,
 and $\lnabp{}=\hlnab{}+ S'$, by Theorem 1.1 of \cite{[S-PV]}, we have
\begin{equation}
p_1(R')-p_1(R^M)= -\frac{1}{2\pi^2}d\left(\langle {\cal S}'\wedge
(\hat{R}^M-\ha d{\cal S}'-\sm{\frac{1}{3}}({\cal S}')^2)
\rangle_{\hat{g}}\right)
\end{equation}
where
${\cal S}':TM\!\ra\! \bigwedge^2TM$,
$({\cal S}')^2$, $\hat{R}^M:\bigwedge^2TM \!\ra\! \bigwedge^2TM$ are
defined by
\[\begin{array}{rcl}
\langle{\cal S}'(X),Y\wedge Z\rangle_{\hat{g}}&=&\hat{g}(S'(X,Y),Z)
\\
\langle ({\cal S}')^2(X\wedge Y), Z\wedge W\rangle_{\hat{g}}&=&
\hat{g}(S'(X,Z), S'(Y,W))- \hat{g}(S'(X,W), S'(Y,Z))
\\
\langle \hat{R}^M(X\wedge Y),Z\wedge W\rangle_{\hat{g}}&=&
h \La{(}R^M(X, Y, Z,W) + \phi\bullet g(X,Y, Z, W)\La{)}
\end{array}\]
where $\hat{R}^M$ is the curvature tensor of $(M,\hat{g}=h g)$,
and
\begin{equation}
\phi=\frac{1}{2}\La{(}-\frac{\|\nabla \log h\|^{2}}{4}g
+\frac{1}{2}d \log h \otimes d \log h- Hess( \log h)\La{)}.
\end{equation}
The inner product $\langle, \rangle_{\hat{g}}$ is the
usual inner product on $\bigwedge^2TM$, defined
w.r.t. $\hat{g}$. So we will compute all the terms in (5.11).
\\[2mm]
Let $\lh=\log h$. The letters $A,B,...$ denote vector fields 
$X,Y,Z, \nabla \lh$ or $e_A,...$, and we denote by
\[
\Phi_AB=\Phi^{-1}\lnab{A}\Phi(B)~~~~~~
 \Phi_ABC=g(\Phi^{-1}\lnab{A}\Phi(B),C)~~~~~~~\lh_A=d\lh(A).
\]
 The gradient $\nabla \lh$ is w.r.t. $g$. From lemma 5.1, 
\begin{eqnarray}
\Phi_ABC+\Phi_ACB &=& \lh_Ag(B,C).
\end{eqnarray}
We easily derive 
\begin{eqnarray}
\lh_Z\Phi_YXA -\lh_Y\Phi_ZXA
&=&-\lh_Z\Phi_YA X+\lh_Y\Phi_ZA X\\
\Phi_YAZ-\Phi_ZAY &=& g(\Phi^{-1}d\Phi(Z,Y),A)+g(\lh_Y Z-\lh_ZY, A)\\
\Phi_YZ\nabla\lh -\Phi_ZY\nabla \lh
&=& -\Phi_Y\nabla\lh Z +\Phi_Z\nabla \lh Y
=g(\Phi^{-1}d\Phi(Y,Z),\nabla\lh).
\end{eqnarray}
Now we have
\begin{eqnarray*}
S'(X,Y)&=& \Phi^{-1}\lnab{X}\Phi(Y)-\ha \lh_XY-\ha \lh_YX
+\ha g(X,Y)\nabla \lh\\
\hlnab{X}Y&=& \lnab{X}Y+\ha \lh_XY +\ha \lh_YX
-\ha g(X,Y)\nabla \lh.
\end{eqnarray*}
Since $\lnabp{}~$ is a $\hat{g}$-Riemannian connection
and $S'(X,Y)=\lnabp{X}Y-\hlnab{X}Y$, then  $\hat{g}(S'(X,Y),Z)=
-\hat{g}(S'(X,Z),Y)$, and so, the same holds w.r.t. $g$. Let
$e_i$ be a $g$-o.n. basis of $T_pM$. We have
\begin{equation}
g(S'(X,e_i),e_j)=\Phi_Xij-\ha\lh_X \delta_{ij}
-\ha \lh_i g(X,e_j)+ \ha g(X,e_i)\lh_j.
\end{equation}
We consider from now on  $TM$ with the metric $g$, and 
from a tensor $\varrho\in C^{\infty}(\bigwedge^2TM
\otimes \bigwedge^2TM)$ we  define a 4-tensor on $M$
 as 
$\varrho(X,Y,Z,W)=\langle \varrho(X\wedge Y), Z\wedge W\rangle_g$, 
where $\langle X\wedge Y, Z\wedge W\rangle_g= g(X,Z)g(Y,W)-g(X,W)g(Y,Z)$
is the Riemannian structure in $\bigwedge^2TM$ defined w.r.t. $g$.
For each tangent vector $X$ of $T_pM$ we denote by 
$\hat{X}= h^{-\frac{1}{2}}X$.\\[-3mm] 
\begin{Lm} If $\varrho\in C^{\infty}(\bigwedge^2TM\otimes \bigwedge^2TM)$ 
and $e_i$ is a $g$-o.n. basis of $T_pM$,
$X,Y,Z\in T_pM$
\[\cyclic_{\ti{X,Y,Z}}\sm{\sum_{ij}}g(S'(X,e_i), e_j)\varrho(Y,Z,e_i,e_j)=
b(\varrho)(X,Y,Z,\nabla \lh)
-\cyclic_{\ti{X,Y,Z}}\sm{\sum_{i}}\varrho
(Y,Z,\Phi^{-1}\lnab{X}\Phi(e_i), e_i).\]
\end{Lm}
\noindent
\em Proof. \em  
Since $\varrho(Y,Z,e_i,e_j)$ is skew symmetric on $(e_i,e_j)$
\begin{eqnarray*}
\cyclic_{\ti{X,Y,Z}}\sm{\sum_{ij}}g(S'(X,e_i), e_j)\varrho(Y,Z,e_i,e_j)
&=&\cyclic_{\ti{X,Y,Z}}~~\sm{\sum_{i}}
\varrho(Y,Z,e_i,\Phi^{-1}\lnab{X}\Phi(e_i))
+\varrho(Y,Z,X,\nabla\lh).~\qed\\[-1mm]
\end{eqnarray*}
From previous lemma and (5.13) we obtain for any  2-tensor $\xi\in
C^{\infty}(\bigotimes^2TM^*)$
\begin{eqnarray*}
\lefteqn{\sm{\cyclic_{\ti{X,Y,Z}}\sum_{ij}g( S'(X,{e}_i),{e}_j)
\xi\bullet g(Y,Z,e_i,e_j)}=}\\
&=&
\sm{\cyclic_{\ti{X,Y,Z}}-\xi(Y,\Phi^{-1}\lnab{X}\Phi(Z))
+\xi(Z,\Phi^{-1}\lnab{X}\Phi(Y))+\sum_i
\xi(Y,e_i)\Phi_XiZ-\xi(Z,e_i)\Phi_XiY}\\
&=&
\sm{\cyclic_{\ti{X,Y,Z}}-\xi(Y,\Phi^{-1}\lnab{X}\Phi(Z))
+\xi(Z,\Phi^{-1}\lnab{X}\Phi(Y))}\\
&&+\sm{\cyclic_{\ti{X,Y,Z}}\sum_i
\xi(Y,e_i)(-\Phi_XZi+\lh_Xg(Z,e_i))
-\xi(Z,e_i)(-\Phi_XYi+\lh_X g(Y,e_i))}\\
&=& \sm{\cyclic_{\ti{X,Y,Z}}2\xi(Z,\Phi^{-1}d\Phi(X,Y)).}
\end{eqnarray*}
Thus
\begin{Lm} 
\begin{eqnarray}
\sm{\cyclic_{\ti{X,Y,Z}}\sum_{ij}g( S'(X,{e}_i),{e}_j)R^M(Y,Z,e_i,e_j)}&=&
\sm{\cyclic_{\ti{X,Y,Z}}\sum_i-R^M(Y,Z,\Phi^{-1}\lnab{X}\Phi(e_i), e_i)}\\
\sm{\cyclic_{\ti{X,Y,Z}}\sum_{ij}g( S'(X,{e}_i),{e}_j)
\xi\bullet g(Y,Z,e_i,e_j)}
&=& \sm{  \cyclic_{\ti{X,Y,Z}}2\xi(X,\Phi^{-1}d\Phi(Y,Z)).}\\\non
\end{eqnarray}
\end{Lm}
\begin{Lm} Let $~\varrho(Y,Z,A,B)= {g}(S'(Y, A), S'(Z, B))
-{g}(S'(Y, B), S'(Z,A))$.  We have
\begin{eqnarray}
\lefteqn{\cyclic_{\ti{X,Y,Z}} \sm{\sum_{ij}}{g(S'(X,e_i),e_j)
\varrho(Y,Z,e_i,e_j)}=}\\
&=& \!\!\!\sm{\cyclic_{\ti{X,Y,Z}}\left(\sum_{ij}~(\Phi_Xij-\Phi_Xji)\,
g(\Phi^{-1}\lnab{Y}\Phi(e_i), \Phi^{-1}\lnab{Z}\Phi(e_j))
-3g(\Phi^{-1}\lnab{X}\Phi(\nabla \lh),
\Phi^{-1}d\Phi(Y,Z))\right)}\non\\
&&\sm{+\cyclic_{\ti{X,Y,Z}}\left(~\frac{3}{2}d\lh\otimes 
d\lh(X, \Phi^{-1}d\Phi(Y,Z)) 
+\frac{3}{4}\|\nabla \lh\|^2g(X, \Phi^{-1}d\Phi(Y,Z))~\right).}\non
\end{eqnarray}
\end{Lm}
\noindent
\em Proof. \em  
\begin{eqnarray*}
\sm{g(S'(Y,e_i), S'(Z,e_j))}&=&\sm{g\La{(}
\Phi^{-1}\lnab{Y}\Phi(e_i)-\ha \lh_Ye_i-\ha \lh_iY+
\ha g(Y,e_i)\nabla \lh~,~}\\[-2mm] 
&&~~~~~~~~~\sm{\Phi^{-1}\lnab{Z}\Phi(e_j)-\ha \lh_Z e_j-\ha \lh_jZ+
\ha g(Z,e_j)\nabla \lh\La{)}}.
\end{eqnarray*}
Using (5.13) and  the fact that $g(S'(X,e_i),e_j)$ is skew symmetric
on $(i,j)$ we have, after
interchanging $i$ with $j$ in some terms,
\begin{eqnarray}
\lefteqn{\cyclic_{\ti{X,Y,Z}} \sum_{ij}g(S'(X,e_i),e_j)
\varrho(Y,Z,e_i,e_j)=}\non\\
&=&\!\!\!\cyclic_{\ti{X,Y,Z}} \sum_{ij}\sm{g(S'(X,e_i),e_j)\La{(}
2g(\Phi^{-1}\lnab{Y}\Phi(e_i), \Phi^{-1}\lnab{Z}\Phi(e_j))
-\lh_Z\Phi_Yij +\lh_Y \Phi_Zij}\\[-5mm]
&&~~~~~~~~~~~~~~~~~~~~~~~~
\sm{-\lh_j\Phi_YiZ +\lh_j \Phi_{Z}iY
+g(Z,e_j)\Phi_{Y}i\nabla \lh -g(Y,e_j)\Phi_{Z}i\nabla \lh}~~\La{)}\\
&&\!\!\!+\cyclic_{\ti{X,Y,Z}} \sum_{ij}\sm{g(S'(X,e_i),e_j)\LA{(}
\sm{(-\frac{3}{4} d\lh \otimes d \lh+\frac{1}{8}\|\nabla \lh\|^2g)
\bullet g (Y,Z,e_i,e_j)}\LA{)}}
\end{eqnarray}
\begin{eqnarray}
\lefteqn{(5.21)+(5.22)=}\non\\
&=&\!\!\!\cyclic_{\ti{X,Y,Z}} \sum_{ij}\sm{\LA{(} 
\Phi_{X}ij-\ha \lh_X\delta_{ij} 
-\ha \lh_ig(X,e_j)+ \ha g(X,e_i)\lh_j
\LA{)}}
\LA{(}\sm{
2g(\Phi^{-1}\lnab{Y}\Phi(e_i), \Phi^{-1}\lnab{Z}\Phi(e_j))}~~~~~~\non\\[-4mm]
&&~~~~~~~~~~~\sm{-\lh_Z \Phi_{Y}ij
+\lh_Y \Phi_{Z}ij
-\lh_j\Phi_{Y}iZ
+\lh_j \Phi_{Z}iY}
\sm{+g(Z,e_j)\Phi_{Y}i\nabla \lh}
\sm{-g(Y,e_j)\Phi_{Z}i\nabla \lh}~\LA{)}\non\\[2mm]
&=& \!\!\!\sm{\cyclic_{\ti{X,Y,Z}}\sum_{ij}~+2\Phi_Xij\,
g(\Phi^{-1}\lnab{Y}\Phi(e_i), \Phi^{-1}\lnab{Z}\Phi(e_j))
-\lh_Z\Phi_Xij\Phi_Yij
+\lh_Y\Phi_Xij\Phi_Zij}\non\\
&&\!\!\!+\sm{\cyclic_{\ti{X,Y,Z}}\sum_i-\Phi_Xi\nabla \lh\, \Phi_YiZ
+\Phi_Xi\nabla \lh\, \Phi_ZiY 
+\Phi_XiZ\Phi_Yi\nabla \lh
-\Phi_XiY\Phi_Zi\nabla \lh}\non\\
&&\!\!\!+\sm{\cyclic_{\ti{X,Y,Z}}\sum_i-\lh_X
g(\Phi^{-1}\lnab{Y}\Phi(e_i), \Phi^{-1}\lnab{Z}\Phi(e_i))
+\lh_X\lh_Z\lh_Y-\lh_X\lh_Y\lh_Z }\non\\
&&\!\!\!+\sm{\cyclic_{\ti{X,Y,Z}}+\ha\lh_X\Phi_Y\nabla\lh Z
-\ha\lh_X\Phi_Z\nabla\lh Y-\ha\lh_X\Phi_YZ\nabla\lh 
+\ha\lh_X\Phi_ZY\nabla\lh}\\
&&\!\!\!+\sm{\cyclic_{\ti{X,Y,Z}}-g(\Phi^{-1}\lnab{Y}\Phi(\nabla \lh),
\Phi^{-1}\lnab{Z}\Phi(X))+\ha\lh_Z\Phi_Y\nabla\lh X
-\ha\lh_Y\Phi_Z\nabla\lh X}\\
&&\!\!\!+\sm{\cyclic_{\ti{X,Y,Z}}
+\ha\lh_X\Phi_Y\nabla \lh Z-\ha\lh_X\Phi_Z\nabla \lh Y
-\frac{1}{4}g(X,Z)\lh_Y\|\nabla\lh\|^2
+\frac{1}{4}g(X,Y)\lh_Z\|\nabla\lh\|^2}\\
&&\!\!\!+\sm{\cyclic_{\ti{X,Y,Z}}+g(\Phi^{-1}\lnab{Y}\Phi(X),
\Phi^{-1}\lnab{Z}\Phi(\nabla \lh))-\frac{1}{2}\lh_Z \Phi_YX\nabla \lh
+\frac{1}{2}\lh_Y \Phi_ZX\nabla \lh}\\
&&\!\!\!+\sm{\cyclic_{\ti{X,Y,Z}}-\frac{1}{2}\|\nabla\lh\|^2\Phi_YXZ
+\frac{1}{2}\|\nabla\lh\|^2\Phi_ZXY+\frac{1}{2}\lh_Z\Phi_YX\nabla \lh
-\frac{1}{2}\lh_Y\Phi_ZX\nabla \lh}
\end{eqnarray}
The last two terms of (5.27) cancel with the last two of (5.28).
From (5.14), (5.15) and (5.16)
$\sm{\cyclic_{\ti{X,Y,Z}}\ha\lh_Z\Phi_Y\nabla \lh X-\ha\lh_Y\Phi_Z\nabla \lh X=
\cyclic_{\ti{X,Y,Z}}\ha\lh_X\Phi_Z\nabla \lh Y-\ha\lh_X\Phi_Y\nabla \lh Z
=\cyclic_{\ti{X,Y,Z}}\ha\lh_Xg(\Phi^{-1}d\Phi(Y,Z), \nabla \lh)}$
that we replace in (5.25), and 
$\ha\lh_X\Phi_Y\nabla\lh Z
-\ha\lh_X\Phi_Z\nabla\lh Y =
-\ha\lh_X\Phi_YZ\nabla\lh
+\ha\lh_X\Phi_ZY\nabla\lh=
\ha\lh_X g(\Phi^{-1}d\Phi(Z,Y),
\nabla \lh)$
that we replace in (5.24) and (5.26).
We also have $\forall ij$
\begin{equation}
\cyclic_{\ti{X,Y,Z}}\lh_Z\Phi_Xij\Phi_Yij
-\lh_Y\Phi_Xij\Phi_Zij=0.
\end{equation}
Thus,
\begin{eqnarray}
\lefteqn{(5.21)+(5.22)=}\non\\
&=& \!\!\!\sm{\cyclic_{\ti{X,Y,Z}}\sum_{ij}~+2\Phi_Xij\,
g(\Phi^{-1}\lnab{Y}\Phi(e_i), \Phi^{-1}\lnab{Z}\Phi(e_j))}~~~~\\
&&\!\!\!+\sm{\cyclic_{\ti{X,Y,Z}}\sum_i~-\Phi_Xi\nabla \lh\, \Phi_YiZ
+\Phi_Xi\nabla \lh\, \Phi_ZiY
+\Phi_XiZ\Phi_Yi\nabla \lh
-\Phi_XiY\Phi_Zi\nabla \lh}\\
&&\!\!\!+\sm{\cyclic_{\ti{X,Y,Z}}\sum_i~-\lh_X
g(\Phi^{-1}\lnab{Y}\Phi(e_i), \Phi^{-1}\lnab{Z}\Phi(e_i))}\\
&&\!\!\!+\sm{\cyclic_{\ti{X,Y,Z}}-\lh_Xg(\Phi^{-1}d\Phi(Y,Z),
\nabla \lh)}\\
&&\!\!\!+\sm{\cyclic_{\ti{X,Y,Z}}
-g(\Phi^{-1}\lnab{Y}\Phi(\nabla \lh),
\Phi^{-1}\lnab{Z}\Phi(X)) +g(\Phi^{-1}\lnab{Y}\Phi(X),
\Phi^{-1}\lnab{Z}\Phi(\nabla \lh))}\\
&&\!\!\!+\sm{\cyclic_{\ti{X,Y,Z}}\frac{1}{2}\|\nabla\lh\|^2(\Phi_XYZ-\Phi_YXZ)
-\frac{1}{4}\|\nabla \lh\|^2
(\lh _Y g(Z,X)-\lh_Z g(Y,X)).}
\end{eqnarray}
Moreover
\[(5.34)=\cyclic_{\ti{X,Y,Z}}-g(\Phi^{-1}\lnab{X}\Phi(\nabla \lh), 
\Phi^{-1}d\Phi(Y, Z)).\]
 Since $\cyclic_{\ti{X,Y,Z}}\Phi_Zi\nabla \lh \Phi_XiY=
\cyclic_{\ti{X,Y,Z}}\Phi_Yi\nabla \lh \Phi_ZiX$, then
\begin{eqnarray*}
(5.31)&=&\sm{\cyclic_{\ti{X,Y,Z}}\sum_i-\Phi_Xi\nabla \lh (\Phi_YiZ-\Phi_ZiY)
+\Phi_Yi\nabla \lh(\Phi_XiZ-\Phi_ZiX)}\\
&=&\sm{\cyclic_{\ti{X,Y,Z}}\sum_i-2\Phi_Xi\nabla \lh(\Phi_YiZ-\Phi_ZiY)}\\
&=&\sm{\cyclic_{\ti{X,Y,Z}}\sum_i-2\Phi_Xi\nabla \lh g(\Phi^{-1}d\Phi(Z,Y),e_i)
-2\lh_Y\Phi_XZ\nabla \lh + 2\lh_Z\Phi_XY\nabla \lh}\\
&=& \sm{\cyclic_{\ti{X,Y,Z}}\sum_i-2\Phi_Xi\nabla \lh g(\Phi^{-1}d\Phi(Z,Y),e_i)
-2\lh_Yg(\Phi^{-1}d\Phi(X,Z), \nabla \lh)}\\
&=&\sm{\cyclic_{\ti{X,Y,Z}}\sum_i-2(-\Phi_X\nabla \lh i
+\lh_X\lh_i) g(\Phi^{-1}d\Phi(Z,Y),e_i)
-2\lh_Yg(\Phi^{-1}d\Phi(X,Z), \nabla \lh)}\\
&=&\sm{\cyclic_{\ti{X,Y,Z}}
 2g(\Phi^{-1}d\Phi(Z,Y),\Phi^{-1}\lnab{X}\Phi(\nabla\lh))
-2\lh_Xg(\Phi^{-1}d\Phi(Z,Y),\nabla \lh)
-2\lh_Yg(\Phi^{-1}d\Phi(X,Z), \nabla \lh)}\\
&=&\sm{\cyclic_{\ti{X,Y,Z}}
 -2g(\Phi^{-1}\lnab{X}\Phi(\nabla\lh)), \Phi^{-1}d\Phi(Y,Z))
+4\lh_Xg(\Phi^{-1}d\Phi(Y,Z),\nabla \lh).}
\end{eqnarray*}
We have $\sm{\Phi_XYZ-\Phi_YXZ=g(\Phi^{-1}d\Phi(X,Y),Z)}$ and 
$\sm{\cyclic_{\ti{X,Y,Z}}\frac{1}{4}\|\nabla \lh\|^2
(\lh _X g(Y,Z)-\lh_Z g(Y,X))=0}$, that we will replace in (5.35).
Moreover
\begin{equation}
(5.30)+(5.32)=\cyclic_{\ti{X,Y,Z}}\sm{\sum_{ij}}~(\Phi_Xij-\Phi_Xji)\,
g(\Phi^{-1}\lnab{Y}\Phi(e_i), \Phi^{-1}\lnab{Z}\Phi(e_j)).
\end{equation}
Therefore,
\begin{eqnarray*}
(5.21)+(5.22) &=&
\sm{\cyclic_{\ti{X,Y,Z}}\sum_{ij}~(\Phi_Xij-\Phi_Xji)\,
g(\Phi^{-1}\lnab{Y}\Phi(e_i), \Phi^{-1}\lnab{Z}\Phi(e_j))}~~~\non\\
&&\!\!\!\sm{+\cyclic_{\ti{X,Y,Z}} -3g(\Phi^{-1}\lnab{X}\Phi(\nabla \lh),
\Phi^{-1}d\Phi(Y,Z))+ 3g(\lh_X\nabla\lh, \Phi^{-1}d\Phi(Y,Z))} \non\\
&&\sm{ +\cyclic_{\ti{X,Y,Z}}
\ha\|\nabla\lh\|^2g(X,\Phi^{-1}d\Phi(Y,Z))}.~~~~~~~~~~~~~\qed
\end{eqnarray*}
\begin{Pp}
\begin{eqnarray}
\lefteqn{
\langle {\cal S}'\wedge (\hat{R}^M-\ha d{\cal S}'-\sm{\frac{1}{3}}
({\cal S}')^2)\rangle_{\hat{g}} (X,Y,Z)=~~~~~~~~~~~~}\\[2mm]
&=&\sm{\cyclic_{\ti{X,Y,Z}}\frac{1}{4}\langle\Phi^{-1}R^{\bot}(Y,Z)
+ R^M(Y,Z),\Phi^{-1}\lnab{X}\Phi\rangle -
\frac{1}{4}d\La{(}g(\Phi^{-1}d\Phi(\cdot, \cdot), \nabla \lh)\La{)}(X,Y,Z)}\\
&&\sm{+\cyclic_{\ti{X,Y,Z}}
\frac{1}{12}\langle \Phi^{-1}\lnab{X}\Phi,
[\Phi^{-1}\lnab{Y}\Phi,\Phi^{-1}\lnab{Z}\Phi]\rangle.}
\end{eqnarray}
\end{Pp}
\noindent
\em Proof. \em 
Let $e_i$  a $g$-orthonormal frame of $TM$. Then
$\hat{e}_i$ is a  $\hat{g}$-orthonormal frame. We have
\begin{eqnarray}
\lefteqn{\sm{
\langle {\cal S}'\wedge (\hat{R}-\ha d{\cal S}'-\sm{\frac{1}{3}}
({\cal S}')^2)\rangle_{\hat{g}} (X,Y,Z)=\cyclic_{\ti{X,Y,Z}}
\langle {\cal S}'(X), (\hat{R}-\ha d{\cal S}'-\sm{\frac{1}{3}}
({\cal S}')^2)(Y,Z)\rangle_{\hat{g}}=}}\\
&=&\!\! \sm{\cyclic_{\ti{X,Y,Z}}
\sum_{ij}\ha\langle {\cal S}'(X), \hat{e}_i\wedge
\hat{e}_j\rangle_{\hat{g}}\left(
\hat{R}(Y,Z,\hat{e}_i, \hat{e}_j)-\ha \langle d{\cal S}'(Y\wedge Z),
\hat{e}_i\wedge \hat{e}_j\rangle_{\hat{g}}
-\sm{\frac{1}{3}} \langle ({\cal S}')^2(Y\wedge Z),
\hat{e}_i\wedge \hat{e}_j\rangle_{\hat{g}}\right)}\non \\
&=&\sm{ \cyclic_{\ti{X,Y,Z}}\sum_{ij}\ha \hat{g}(S'(X, \hat{e}_i),
\hat{e}_j)\left( hR^M(Y,Z, \hat{e}_i, \hat{e}_j)+
h \phi\bullet g(Y,Z, \hat{e}_i, \hat{e}_j)
-\ha \langle d{\cal S}'(Y\wedge Z),
\hat{e}_i\wedge \hat{e}_j\rangle_{\hat{g}}\right)}\non \\[-2mm]
&&~~~~~~~~~~~~-\sm{\frac{1}{6} \hat{g}(S'(X, \hat{e}_i),
\hat{e}_j)\left( \hat{g}(S'(Y, \hat{e}_i), S'(Z, \hat{e}_j))
-\hat{g}(S'(Y, \hat{e}_j), S'(Z, \hat{e}_i))\right)}\non\\[1mm]
&=&\sm{ \cyclic_{\ti{X,Y,Z}}\sum_{ij}\ha {g}(S'(X, {e}_i),
{e}_j)\left( R^M(Y,Z, {e}_i, {e}_j)+
 \phi\bullet g(Y,Z, {e}_i, {e}_j)\right)}\\[-3mm]
&&~~~~~~~~~~~~-\sm{\frac{1}{4}{g}(S'(X, {e}_i),{e}_j)
 \La{(}\langle \hlnab{Y}{\cal S}(Z)-\hlnab{Z}{\cal S}(Y),
\hat{e}_i\wedge \hat{e}_j\rangle_{\hat{g}}\La{)}}\\
&&~~~~~~~~~~~~-\sm{\frac{1}{6} {g}(S'(X, {e}_i),
{e}_j)\left( {g}(S'(Y, {e}_i), S'(Z, {e}_j))
-{g}(S'(Y, {e}_j), S'(Z,{e}_i))\right).}
\end{eqnarray}
We assume that at a given point $p_0$, $\sm{\hlnab{}X=\hlnab{}Y=
\hlnab{}Z=\hlnab{}\hat{e}_i=0}$.
Thus, at $p_0$,
$\sm{\lnab{X}Y=-\hat{S}(X,Y)}$ $\sm{=-\ha \lh_XY-\ha \lh_YX
+\ha g(X,Y)\nabla\lh}$, 
and similarly for the other vector fields.
The following computations
are computed at $p_0$.
\begin{eqnarray}
\sm{d(g(X,Y))(Z)} &=& \sm{-\lh_Z\, g(X,Y)}\\
\sm{d(\lh_X)(Y)}&=&\sm{ Hess \lh (X,Y)-\lh_X \lh_Y
+\ha \|\nabla \lh \|^2 g(X,Y)}~~~~
\end{eqnarray}
and since $\sm{\bar{R}(X,Y)\Phi=-\lnab{X,Y}^2\Phi+\lnab{Y,X}^2\Phi}$,
we have
\begin{eqnarray}
\lefteqn{\sm{d(\Phi_XZW)(Y)
-d(\Phi_{Y}ZW)(X)=}}\\
&=& \sm{-2\lh_Y \Phi_{X}ZW+
2\lh_X \Phi_{Y}ZW} 
\sm{+g(\Phi^{-1}(\bar{R}(X,Y)\Phi)(Z), W)}\non\\
&&
\sm{-\ha \lh_Zg(\Phi^{-1}d\Phi(X,Y), W)+\ha g(Y,Z)\Phi_{X}\nabla \lh W
-\ha g(X,Z)\Phi_{Y}\nabla \lh W}\non\\
&&\sm{+g(\Phi^{-1}\lnab{X}\Phi(Z), \Phi^{-1}\lnab{Y}\Phi(W))
-g(\Phi^{-1}\lnab{Y}\Phi(Z), \Phi^{-1}\lnab{X}\Phi(W))}\non\\
&&\sm{-\ha \lh_W \Phi_{X}ZY
+\ha g(Y,W)\Phi_{X}Z\nabla \lh +\ha \lh_W \Phi_{Y}ZX
-\ha g(X,W)\Phi_{Y}Z\nabla \lh}.\non
\end{eqnarray}
Applying eqs (5.44),(5.45), (5.46) and (5.17) we get
\begin{eqnarray}
(5.42)&=&\sm{ \cyclic_{\ti{X,Y,Z}}\sum_{ij}-\sm{\frac{1}{4}}g(S'(X,e_i),e_j)
\LA{(}\lnab{Y}(\langle {\cal S}'(Z),\hat{e}_i\wedge
\hat{e}_j\rangle_{\hat{g}})-\lnab{Z}(\langle {\cal S}'(Y),\hat{e}_i\wedge
\hat{e}_j\rangle_{\hat{g}})\LA{)}}\non\\
&=&\sm{\cyclic_{\ti{X,Y,Z}}\sum_{ij}-\sm{\frac{1}{4}}g(S'(X,e_i),e_j)
\LA{(}\lnab{Y}(g( S'(Z,{e}_i),{e}_j)
-\lnab{Z}(g( S'(Y,{e}_i),{e}_j)\LA{)}}\non\\
&=&\!\!\!\sm{\cyclic_{\ti{X,Y,Z}}\sum_{ij}-\frac{1}{4}\La{(}
\Phi_Xij -\ha \delta_{ij}\lh_X -\ha \lh_i g(X,e_j)+\ha g(X,e_i)\lh_j
\La{)}\cdot}\\[-2mm]
&&~~~~\sm{\cdot \La{(} -2\lh_Y\Phi_{Z}ij +2\lh_Z \Phi_{Y}ij
+ g(\Phi^{-1}(\bar{R}(Z,Y)\Phi)(e_i), e_j)}\\
&&~~~~~~\sm{-\ha \lh_i g(\Phi^{-1}d\Phi(Z,Y), e_j)
+\ha g(Y, e_i)\Phi_{Z}\nabla \lh j
-\ha g(Z, e_i)\Phi_{Y}\nabla \lh j}\\
&&~~~~~~\sm{+g(\Phi^{-1}\lnab{Z}\Phi(e_i), 
\Phi^{-1}\lnab{Y}\Phi(e_j)) -g(\Phi^{-1}\lnab{Y}\Phi(e_i), 
\Phi^{-1}\lnab{Z}\Phi(e_j))}\\
&&~~~~~~\sm{-\ha \lh_j \Phi_{Z}iY+\ha 
g(Y,e_j)\Phi_{Z}i\nabla \lh
+\ha \lh_j\Phi_{Y}iZ
-\ha g(Z,e_j)\Phi_{Y}i\nabla \lh}~\La{)}\\
&&~~~-\sm{\frac{1}{8}g(S'(X,e_i),e_j)
\cdot \LA{(}-Hess \lh + 2d\lh \otimes d \lh 
-\ha \|\nabla \lh\|^2 g\LA{)} \bullet g (Y,Z, e_i, e_j)}
\end{eqnarray}
\begin{eqnarray}
\lefteqn{(5.47)+\ldots +(5.51)=}\non\\
&=&\sm{\cyclic_{\ti{X,Y,Z}}\sum_{ij}(\ha \lh_Y\Phi_Xij\Phi_Zij
-\ha \lh_Z\Phi_Xij\Phi_Yij)}\\
&&+\sm{\cyclic_{\ti{X,Y,Z}} 
- \frac{1}{4}\langle\Phi^{-1}\bar{R}(Z,Y)\Phi,
\Phi^{-1}\lnab{X}\Phi\rangle 
+\frac{1}{8}\langle\Phi^{-1}\lnab{X}\Phi(\nabla \lh),\Phi^{-1}d\Phi(Z,Y)
\rangle}\\
&&+\sm{\cyclic_{\ti{X,Y,Z}}-\frac{1}{8}\langle \Phi^{-1}\lnab{X}\Phi(Y),
\Phi^{-1}\lnab{Z}\Phi(\nabla \lh)\rangle 
+\frac{1}{8}\langle \Phi^{-1}\lnab{X}\Phi(Z),
\Phi^{-1}\lnab{Y}\Phi(\nabla \lh)\rangle}\\
&&+\sm{\cyclic_{\ti{X,Y,Z}}\sum_{ij}-\frac{1}{4}\Phi_Xij \,
g(\Phi^{-1}\lnab{Z}\Phi(e_i), \Phi^{-1}\lnab{Y}\Phi(e_j))
+\frac{1}{4}\Phi_Xij \,
g(\Phi^{-1}\lnab{Y}\Phi(e_i), \Phi^{-1}\lnab{Z}\Phi(e_j))}~~~~~~~~\non\\
&&+\sm{\cyclic_{\ti{X,Y,Z}}\sum_{i} \frac{1}{8}\Phi_Xi\nabla \lh
\Phi_ZiY -\frac{1}{8}\Phi_XiY\Phi_Zi\nabla \lh 
-\frac{1}{8}\Phi_Xi\nabla \lh
\Phi_YiZ +\frac{1}{8}\Phi_XiZ\Phi_Yi\nabla \lh}\\
&& +\sm{\cyclic_{\ti{X,Y,Z}}-\frac{1}{4}\lh_Y\Phi_Z\nabla \lh X 
+\frac{1}{4}\lh_Z\Phi_Y\nabla \lh X
+\frac{1}{4}\lh_Y\Phi_ZX\nabla \lh 
-\frac{1}{4}\lh_Z\Phi_YX\nabla \lh }\\
&&+\sm{\cyclic_{\ti{X,Y,Z}}+\frac{1}{8}
g(\Phi^{-1}(\bar{R}(Z,Y)\Phi)(\nabla \lh) ,X)-
\frac{1}{8}g(\Phi^{-1}(\bar{R}(Z,Y)\Phi)(X), \nabla \lh)}\non\\
&&+\sm{\cyclic_{\ti{X,Y,Z}}
-\frac{1}{16}\|\nabla \lh\|^2g(\Phi^{-1}d\Phi(Z,Y),X)
+\frac{1}{16}\lh_Xg(\Phi^{-1}d\Phi(Z,Y),\nabla \lh)}\non
\end{eqnarray}
\begin{eqnarray}
&&+\sm{\cyclic_{\ti{X,Y,Z}}+\frac{1}{16}\lh_Y\Phi_Z\nabla \lh X
-\frac{1}{16}\lh_Z\Phi_Y\nabla \lh X
-\frac{1}{16}g(Y,X)\Phi_Z\nabla \lh \nabla \lh
+\frac{1}{16}g(Z,X)\Phi_Y\nabla \lh \nabla \lh}\\
&&+\sm{\cyclic_{\ti{X,Y,Z}}\frac{1}{8}g(\Phi^{-1}\lnab{Z}\Phi(\nabla \lh),
\Phi^{-1}\lnab{Y}\Phi(X))-\frac{1}{8}g(\Phi^{-1}\lnab{Y}\Phi(\nabla \lh),
\Phi^{-1}\lnab{Z}\Phi(X))}\\
&&+\sm{\cyclic_{\ti{X,Y,Z}}
-\frac{1}{8}g(\Phi^{-1}\lnab{Z}\Phi(X),
\Phi^{-1}\lnab{Y}\Phi(\nabla\lh))
+\frac{1}{8}g(\Phi^{-1}\lnab{Y}\Phi(X),
\Phi^{-1}\lnab{Z}\Phi(\nabla \lh))}\\
&&+\sm{\cyclic_{\ti{X,Y,Z}}-\frac{1}{16}\lh_X\Phi_Z\nabla \lh Y
+\frac{1}{16}g(Y,X)\Phi_Z\nabla \lh \nabla \lh
+\frac{1}{16}\lh_X\Phi_Y\nabla \lh Z
-\frac{1}{16}g(Z,X)\Phi_Y\nabla \lh \nabla \lh}\\
&&+\sm{\cyclic_{\ti{X,Y,Z}}\frac{1}{16}\|\nabla \lh\|^2\Phi_ZXY
-\frac{1}{16}\lh_Y\Phi_ZX\nabla \lh
-\frac{1}{16}\|\nabla \lh\|^2\Phi_YXZ
+\frac{1}{16}\lh_Z\Phi_YX\nabla \lh}.
\end{eqnarray}
Note that by $(5.29)$, $(5.53)=0$, and the second term of 
$(5.54)$ is equal to 
$(5.55)=(5.59)=(5.60)$.
We also have using $(5.14) (5.16)$
\begin{eqnarray*}
\lefteqn{(5.56)=
\sm{\cyclic_{\ti{X,Y,Z}}\sum_{i} \frac{1}{8}(-\Phi_X\nabla \lh i
+\lh_X\lh_i)(\Phi_ZiY-\Phi_YiZ)}}\\[-1mm]
&&~~~~~~~~~~\sm{-
\frac{1}{8}(-\Phi_XYi+\lh_Xg(i,Y))(-\Phi_Z\nabla \lh i+\lh_Z\lh_i)
+\frac{1}{8}(-\Phi_XZi + \lh_X g(i,Z))
(-\Phi_Y\nabla \lh i +\lh_Y\lh_i)}\\[1mm]
&=&\sm{\cyclic_{\ti{X,Y,Z}}-\frac{1}{8}g(\Phi^{-1}\lnab{X}\Phi(\nabla \lh),
\Phi^{-1}d\Phi(Y,Z))
+ \frac{1}{8}\lh_Xg(\Phi^{-1}d\Phi(Y,Z), \nabla \lh)}\\[-2mm]
&&~~~~~~~\sm{-\frac{1}{8}g(\Phi^{-1}\lnab{X}\Phi(Y), \Phi^{-1}\lnab{Z}\Phi
(\nabla \lh))+\frac{1}{8}g(\Phi^{-1}\lnab{X}\Phi(Z), 
\Phi^{-1}\lnab{Y}\Phi(\nabla \lh))}\\
&&~~~~~~~\sm{+\frac{1}{8}\lh_Z(-\Phi_X\nabla \lh  Y+\Phi_XY\nabla \lh)
+\frac{1}{8}\lh_Y(\Phi_X\nabla \lh Z-\Phi_X Z\nabla \lh)
+\frac{1}{8}\lh_X(\Phi_Z\nabla \lh Y-\Phi_Y\nabla \lh Z)}\\[1mm]
&=&\sm{\cyclic_{\ti{X,Y,Z}}-\frac{1}{4}g(\Phi^{-1}\lnab{X}\Phi(\nabla \lh),
\Phi^{-1}d\Phi(Y,Z))
+ \frac{1}{8}\lh_Xg(\Phi^{-1}d\Phi(Y,Z), \nabla \lh)}\\[-1mm]
&&+\sm{\cyclic_{\ti{X,Y,Z}}
+\frac{1}{8}\lh_Z(-\Phi_X\nabla \lh  Y+\Phi_XY\nabla \lh)
+\frac{1}{8}\lh_Y(\Phi_X\nabla \lh Z-\Phi_X Z\nabla \lh)
+\frac{1}{8}\lh_X(\Phi^{-1}d\Phi(Y,Z), \nabla \lh).}
\end{eqnarray*}
And using again $(5.16)$
\begin{eqnarray*}
\lefteqn{(5.57)+(5.58)+(5.61)+(5.62)=}\\
&=& \sm{\cyclic_{\ti{X,Y,Z}}
\frac{3}{16}\lh_Y(-\Phi_Z\nabla \lh X + \Phi_Z X\nabla \lh)
+\frac{3}{16}\lh_Z(\Phi_Y\nabla \lh X - \Phi_Y X\nabla \lh)}\\
&&+\sm{\cyclic_{\ti{X,Y,Z}}\frac{1}{16}\lh_X(-\Phi_Z\nabla \lh Y +\Phi_Y\nabla
\lh Z) +\frac{1}{16}\|\nabla \lh\|^2 g( \Phi^{-1}d\Phi(Y,Z),X)}\\
&=&\sm{\!\!\!\!\!\!\!\!\!\cyclic_{\ti{X,Y,Z}}
\frac{3}{16}\lh_Y\la{(}-\Phi_X\nabla \lh Z-g( \Phi^{-1}d\Phi(X,Z),\nabla \lh)
-g(\lh_ZX\!-\!\lh_XZ, \nabla \lh)+ g( \Phi^{-1}d\Phi(Z,X),\nabla \lh)
+\Phi_XZ\nabla \lh\la{)}}\\
&+&\sm{\!\!\!\!\!\!\!\!\!\!\cyclic_{\ti{X,Y,Z}}
\frac{3}{16}\lh_Z\la{(}\Phi_X\nabla \lh Y+g( \Phi^{-1}d\Phi(X,Y),\nabla \lh)
+g(\lh_YX\!-\!\lh_XY, \nabla \lh)- g( \Phi^{-1}d\Phi(Y,X),\nabla \lh)
-\Phi_XY\nabla \lh\la{)}}\\
&+&\sm{\!\!\!\!\!\!\!\!\!\cyclic_{\ti{X,Y,Z}}
\frac{1}{16}\lh_X(-\Phi_Z\nabla \lh Y +\Phi_Y\nabla \lh Z)
+\frac{1}{16}\|\nabla \lh\|^2 g(X, \Phi^{-1}d\Phi(Y,Z))}\\
&=&\sm{\cyclic_{\ti{X,Y,Z}}
\frac{3}{16}\lh_Y(-\!\!\Phi_X\nabla \lh Z+\Phi_XZ\nabla \lh) 
+\frac{3}{8}\lh_Yg(\Phi^{-1}d\Phi(Z,X),\nabla \lh)
+\frac{3}{16}\lh_Z(\Phi_X\nabla \lh Y-\Phi_XY\nabla \lh)}\\
&&+\sm{\cyclic_{\ti{X,Y,Z}}\frac{3}{8}\lh_Zg(\Phi^{-1}d\Phi(X,Y),\nabla \lh)
+\frac{1}{16}\lh_X(\Phi^{-1}d\Phi(Z,Y),\nabla \lh)
+\frac{1}{16}\|\nabla \lh\|^2 g(X, \Phi^{-1}d\Phi(Y,Z))}\\
&=&\sm{\cyclic_{\ti{X,Y,Z}}
\frac{3}{16}\lh_Y(-\Phi_X\nabla \lh Z+\Phi_XZ\nabla \lh)
+\frac{3}{16}\lh_Z(\Phi_X\nabla \lh Y-\Phi_XY\nabla \lh)}\\
&&+\sm{\cyclic_{\ti{X,Y,Z}}+\frac{11}{16}\lh_X(\Phi^{-1}d\Phi(Y,Z),\nabla \lh)
+\frac{1}{16}\|\nabla \lh\|^2 g(X, \Phi^{-1}d\Phi(Y,Z)).}
\end{eqnarray*}
Thus, 
\begin{eqnarray}
\lefteqn{(5.47)+\ldots +(5.51)=}\non\\
&=&\sm{\cyclic_{\ti{X,Y,Z}} 
-\frac{1}{4}\langle\Phi^{-1}\bar{R}(Z,Y)\Phi,
\Phi^{-1}\lnab{X}\Phi\rangle
+\frac{1}{2}\langle\Phi^{-1}\lnab{X}\Phi(\nabla \lh),\Phi^{-1}d\Phi(Z,Y)
\rangle}\non\\
&&+\sm{\cyclic_{\ti{X,Y,Z}}\sum_{ij}\frac{1}{4}(\Phi_Xij-\Phi_Xji) \,
g(\Phi^{-1}\lnab{Y}\Phi(e_i), \Phi^{-1}\lnab{Z}\Phi(e_j))}\non\\
&&+\sm{\cyclic_{\ti{X,Y,Z}}-\frac{1}{4}g(\Phi^{-1}\lnab{X}\Phi(\nabla \lh),
\Phi^{-1}d\Phi(Y,Z))
+ \frac{1}{8}\lh_Xg(\Phi^{-1}d\Phi(Y,Z), \nabla \lh)}\non\\
&&+\sm{\cyclic_{\ti{X,Y,Z}}
+\frac{1}{8}\lh_Z(-\Phi_X\nabla \lh  Y+\Phi_XY\nabla \lh)
+\frac{1}{8}\lh_Y(\Phi_X\nabla \lh Z-\Phi_X Z\nabla \lh)
+\frac{1}{8}\lh_X g(\Phi^{-1}d\Phi(Y,Z), \nabla \lh)}~~~~~~\non\\
&&+\sm{\cyclic_{\ti{X,Y,Z}}+\frac{1}{8}g(\Phi^{-1}(\bar{R}(Z,Y)\Phi)
(\nabla \lh),X)-
\frac{1}{8}g(\Phi^{-1}(\bar{R}(Z,Y)\Phi)(X), \nabla \lh)}\non\\
&&+\sm{\cyclic_{\ti{X,Y,Z}}
-\frac{1}{16}\|\nabla \lh\|^2g(\Phi^{-1}d\Phi(Z,Y),X)
+\frac{1}{16}\lh_Xg(\Phi^{-1}d\Phi(Z,Y),\nabla \lh)}\non\\
&&+\sm{\cyclic_{\ti{X,Y,Z}}
\frac{3}{16}\lh_Y(-\Phi_X\nabla \lh Z+\Phi_XZ\nabla \lh)
+\frac{3}{16}\lh_Z(\Phi_X\nabla \lh Y-\Phi_XY\nabla \lh)}\non\\
&&+ \sm{\cyclic_{\ti{X,Y,Z}}
\frac{11}{16}\lh_Xg(\Phi^{-1}d\Phi(Y,Z),\nabla \lh)
+\frac{1}{16}\|\nabla \lh\|^2g(X, \Phi^{-1}d\Phi(Y,Z))}\non\\[1mm]
&=&\sm{\cyclic_{\ti{X,Y,Z}}  -\frac{1}{4}\langle\Phi^{-1}(\bar{R}(Z,Y)\Phi),
\Phi^{-1}\lnab{X}\Phi\rangle
-\frac{3}{4}\langle\Phi^{-1}\lnab{X}\Phi(\nabla \lh),\Phi^{-1}d\Phi(Y,Z)
\rangle}\non\\
&&+\sm{\cyclic_{\ti{X,Y,Z}}\sum_{ij}\frac{1}{4}(\Phi_Xij-\Phi_Xji) \,
g(\Phi^{-1}\lnab{Y}\Phi(e_i), \Phi^{-1}\lnab{Z}\Phi(e_j))
+ \frac{7}{8}\lh_Xg(\Phi^{-1}d\Phi(Y,Z), \nabla \lh)}\non\\
&&+\sm{\cyclic_{\ti{X,Y,Z}}
+\frac{1}{4}g(\Phi^{-1}(\bar{R}(Y,Z)\Phi)X, \nabla \lh)
-\frac{1}{8}\|\nabla \lh\|^2g(\Phi^{-1}d\Phi(Z,Y),X)
}\non\\
&&+\sm{\cyclic_{\ti{X,Y,Z}}
\frac{1}{16}\lh_Y(-\Phi_X\nabla \lh Z+\Phi_XZ\nabla \lh)
+\frac{1}{16}\lh_Z(\Phi_X\nabla \lh Y-\Phi_XY\nabla \lh).}
\end{eqnarray}
Now,
\begin{eqnarray*}
(5.63)&=&\sm{\cyclic_{\ti{X,Y,Z}}
\frac{1}{16}\lh_Y(\Phi_XZ\nabla \lh -\lh_X\lh_Z+\Phi_XZ\nabla \lh)
+\frac{1}{16}\lh_Z(-\Phi_XY\nabla \lh +\lh_X\lh_Y-\Phi_XY\nabla \lh)}\\
&=& \sm{\cyclic_{\ti{X,Y,Z}} \frac{1}{8}\lh_Y\Phi_XZ\nabla \lh -
\frac{1}{8}\lh_Z\Phi_XY\nabla \lh
=\cyclic_{\ti{X,Y,Z}}-\frac{1}{8}\lh_Xg(\Phi^{-1}d\Phi(Y,Z), \nabla \lh).}
\end{eqnarray*}
Finally using Lemma 5.5  with
$\phi=\ha(-\frac{\|\nabla \lh\|^2}{4}g + \ha d\lh \otimes d \lh -Hess\, \lh)$
\begin{equation}
\sm{(5.52)+\!\cyclic_{\ti{X,Y,Z}}\sum_{ij}\ha g(S'(X,e_i),e_j)\phi\bullet 
g(Y,Z,e_i,e_j)}
=\sm{\!\cyclic_{\ti{X,Y,Z}}\!-\frac{1}{4}(Hess\, \lh +d\lh\otimes d\lh)
(X,\Phi^{-1}d\Phi(Y,Z)).}
\end{equation}
Therefore from (5.41), (5.42), (5.43), and Lemma 5.6
\begin{eqnarray*}
\lefteqn{\sm{
\langle {\cal S}'\wedge (\hat{R}-\ha d{\cal S}'-\sm{\frac{1}{3}}
({\cal S}')^2)\rangle_{\hat{g}} (X,Y,Z)=}}\non\\
&=&\!\!\!
\sm{\cyclic_{\ti{X,Y,Z}}\sum_{i}-\ha R^M(Y,Z,\Phi^{-1}\lnab{X}\Phi(e_i),e_i)
+(5.47)+(5.48)+(5.49)+(5.50)+(5.51)+(5.43)+ (5.64)
}\\
&=&\sm{\cyclic_{\ti{X,Y,Z}}\sum_{i}-\ha R^M(Y,Z,\Phi^{-1}\lnab{X}\Phi(e_i),e_i)
  -\frac{1}{4}\langle\Phi^{-1}\bar{R}(Z,Y)\Phi,
\Phi^{-1}\lnab{X}\Phi\rangle}\\
&&+\sm{\cyclic_{\ti{X,Y,Z}}-\frac{3}{4}\langle\Phi^{-1}
\lnab{X}\Phi(\nabla \lh),\Phi^{-1}d\Phi(Y,Z)
\rangle
+\sum_{ij}\frac{1}{4}(\Phi_Xij-\Phi_Xji) \,
g(\Phi^{-1}\lnab{Y}\Phi(e_i), \Phi^{-1}\lnab{Z}\Phi(e_j))}\\
&&+\sm{\cyclic_{\ti{X,Y,Z}}
+ \frac{6}{8}\lh_Xg(\Phi^{-1}d\Phi(Y,Z), \nabla \lh)
+\frac{1}{4}g(\Phi^{-1}(\bar{R}(Y,Z)\Phi)(X), \nabla \lh)
-\frac{1}{8}\|\nabla \lh\|^2g(\Phi^{-1}d\Phi(Z,Y),X)
}\\
&&+\sm{\cyclic_{\ti{X,Y,Z}}-\frac{1}{4}Hess\, \lh(X, \Phi^{-1}d\Phi(Y,Z))
-\frac{1}{4}\lh_X g(\Phi^{-1}d\Phi(Y,Z), \nabla \lh)}
\end{eqnarray*}
\begin{eqnarray*}
&&+\sm{\cyclic_{\ti{X,Y,Z}}\sum_{ij}-\frac{1}{6}(\Phi_Xij-\Phi_Xji)\,
g(\Phi^{-1}\lnab{Y}\Phi(e_i), \Phi^{-1}\lnab{Z}\Phi(e_j))
+\frac{1}{2}g(\Phi^{-1}\lnab{X}\Phi(\nabla \lh),\Phi^{-1}d\Phi(Y,Z))}\\
&&+\sm{\cyclic_{\ti{X,Y,Z}}
-\frac{1}{4}\lh_Xg( \Phi^{-1}d\Phi(Y,Z), \nabla \lh)
-\frac{1}{8}\|\nabla \lh\|^2g(X, \Phi^{-1}d\Phi(Y,Z))}\\
&=&\sm{ \cyclic_{\ti{X,Y,Z}}
\sum_{i}-\ha R^M(Y,Z,\Phi^{-1}\lnab{X}\Phi(e_i), e_i)
-\frac{1}{4}\langle\Phi^{-1}\bar{R}(Z,Y)\Phi,
\Phi^{-1}\lnab{X}\Phi\rangle
+\frac{1}{4}g(\Phi^{-1}(\bar{R}(Y,Z)\Phi)(X), \nabla \lh)}\\
&&+\sm{\cyclic_{\ti{X,Y,Z}}
-\frac{1}{4}\langle\Phi^{-1}\lnab{X}\Phi(\nabla \lh),\Phi^{-1}d\Phi(Y,Z)
\rangle
+\sum_{ij}\frac{1}{12}(\Phi_Xij-\Phi_Xji) \,
g(\Phi^{-1}\lnab{Y}\Phi(e_i), \Phi^{-1}\lnab{Z}\Phi(e_j))}\\
&&+\sm{\cyclic_{\ti{X,Y,Z}}
+ \frac{1}{4}\lh_Xg(\Phi^{-1}d\Phi(Y,Z), \nabla \lh)
-\frac{1}{4}Hess\, \lh(X, \Phi^{-1}d\Phi(Y,Z)).}
\end{eqnarray*}
Note that
\[\sm{  -\langle\Phi^{-1}\bar{R}(Z,Y)\Phi,
\Phi^{-1}\lnab{X}\Phi\rangle = \sum_i
g(\Phi^{-1}R^{\bot}(Y,Z)\Phi(e_i)
, \Phi^{-1}\lnab{X}\Phi(e_i))-g(R^M(Y,Z)e_i, \Phi^{-1}\lnab{X}
\Phi(e_i)).}\]
We have
$\sm{\cyclic_{\ti{X,Y,Z}} g(\Phi^{-1}(\bar{R}(Y,Z)\Phi)(X),
\nabla \lh)=\cyclic_{\ti{X,Y,Z}} g(\Phi^{-1}({R}^{\bot}(Y,Z)\Phi(X)),
\nabla \lh),}$  for $R^M$ satisfies Bianchi equality, and
recall that $\sm{\cyclic_{\ti{X,Y,Z}} \Phi^{-1}({R}^{\bot}(Y,Z)\Phi(X))
=-\Phi^{-1}d^2\Phi(X,Y,Z)}$. Now,
\begin{eqnarray*}
\lefteqn{ \sm{\cyclic_{\ti{X,Y,Z}}\sum_{ij}(\Phi_Xij-\Phi_Xji) \,
g(\Phi^{-1}\lnab{Y}\Phi(e_i), \Phi^{-1}\lnab{Z}\Phi(e_j))}=}\\
&=& \sm{\cyclic_{\ti{X,Y,Z}}\langle \Phi^{-1}\lnab{Y}\Phi, (\Phi^{-1}\lnab{Z}\Phi)
\circ (\Phi^{-1}\lnab{X}\Phi)\rangle 
-\langle (\Phi^{-1}\lnab{Y}\Phi)\circ (\Phi^{-1}\lnab{X}\Phi), 
\Phi^{-1}\lnab{Z}\Phi\rangle}\\
&=&\sm{\cyclic_{\ti{X,Y,Z}} \langle \Phi^{-1}\lnab{X}\Phi,
(\Phi^{-1}\lnab{Y}\Phi)\circ (\Phi^{-1}\lnab{Z}\Phi)
-(\Phi^{-1}\lnab{Z}\Phi)\circ (\Phi^{-1}\lnab{Y}\Phi)\rangle.}
\end{eqnarray*}
Thus
\begin{eqnarray}
\lefteqn{\sm{
\langle {\cal S}'\wedge (\hat{R}-\ha d{\cal S}'-\sm{\frac{1}{3}}
({\cal S}')^2)\rangle_{\hat{g}} (X,Y,Z)=}}\non\\
&=&\sm{ \cyclic_{\ti{X,Y,Z}}
\frac{1}{4}\langle\Phi^{-1}R^{\bot}(Y,Z)\Phi + R^M(Y,Z),
\Phi^{-1}\lnab{X}\Phi\rangle
+\frac{1}{4}g(\Phi^{-1}(R^{\bot}(Y,Z)\Phi(X)), \nabla \lh)}\non\\
&&+\sm{\cyclic_{\ti{X,Y,Z}}
-\frac{1}{4}g(\Phi^{-1}\lnab{X}\Phi(\nabla \lh),\Phi^{-1}d\Phi(Y,Z))
 +\frac{1}{12}\langle \Phi^{-1}\lnab{X}\Phi,
[\Phi^{-1}\lnab{Y}\Phi,\Phi^{-1}\lnab{Z}\Phi]\rangle} \non\\
&&+\sm{\cyclic_{\ti{X,Y,Z}}+\frac{1}{4}(d\lh\otimes d\lh-
Hess\, \lh)(X,\Phi^{-1}d\Phi(Y,Z))}.\non
\end{eqnarray}
Using (5.13)
\[
\sm{\cyclic_{\ti{X,Y,Z}}g(\Phi^{-1}\lnab{X}\Phi(\nabla \lh),
\Phi^{-1}d\Phi(Y,Z))= \cyclic_{\ti{X,Y,Z}}-g(\Phi^{-1}\lnab{X}
\Phi(\Phi^{-1}d\Phi(Y,Z)), \nabla \lh) 
+\lh_Xg(\Phi^{-1}d\Phi(Y,Z),\nabla \lh)}
\]
with $\sm{\lh_Xg(\Phi^{-1}d\Phi(Y,Z),\nabla \lh)=(d\lh\otimes d\lh)(X,
\Phi^{-1}d\Phi(Y,Z))}$.
Since $\lnab{X}\Phi^{-1}=-\Phi^{-1}(\lnab{X}\Phi)\Phi^{-1}$, we have
\[\sm{-\Phi^{-1}d^2\Phi(X,Y,Z)
=-d(\Phi^{-1}d\Phi)(X,Y,Z)+\cyclic_{\ti{X,Y,Z}}
-\Phi^{-1}\lnab{X}\Phi(\Phi^{-1}d\Phi(Y,Z)).}\]
Since $~b(\Phi(R^{\bot}))=-\Phi^{-1}d^2\Phi$ we obtain 
\begin{eqnarray*}
\lefteqn{\sm{\cyclic_{\ti{X,Y,Z}}
g(\Phi^{-1}(R^{\bot}(Y,Z)\Phi(X)), \nabla \lh)=
-g(d(\Phi^{-1}d\Phi)(X,Y,Z),\nabla \lh)-\cyclic_{\ti{X,Y,Z}}g(
\Phi^{-1}\lnab{X}\Phi(\Phi^{-1}d\Phi(Y,Z)), \nabla \lh)=}}\\
&=&\sm{-g(d(\Phi^{-1}d\Phi)(X,Y,Z),\nabla \lh)
+\cyclic_{\ti{X,Y,Z}}g(\Phi^{-1}\lnab{X}
\Phi(\nabla \lh), \Phi^{-1}d\Phi(Y,Z))-(d\lh\otimes d\lh)(X,
\Phi^{-1}d\Phi(Y,Z))}
\end{eqnarray*}
Thus,
\begin{eqnarray}
\lefteqn{
\langle {\cal S}'\wedge (\hat{R}^M-\ha d{\cal S}'-\sm{\frac{1}{3}}
({\cal S}')^2)\rangle_{\hat{g}} (X,Y,Z)=~~~~~~~~~~~~}\\[2mm]
&=&\sm{\cyclic_{\ti{X,Y,Z}}\frac{1}{4}\langle\Phi^{-1}R^{\bot}(Y,Z)
+ R^M(Y,Z),\Phi^{-1}\lnab{X}\Phi\rangle -
\frac{1}{4}g(d(\Phi^{-1}d\Phi)(X,Y,Z), \nabla \lh)}\\
&&\sm{+\cyclic_{\ti{X,Y,Z}}
\frac{1}{12}\langle \Phi^{-1}\lnab{X}\Phi,
[\Phi^{-1}\lnab{Y}\Phi,\Phi^{-1}\lnab{Z}\Phi]\rangle
-\frac{1}{4}Hess\, \lh (X, \Phi^{-1}d\Phi(Y,Z)).}
\end{eqnarray}
Finally
$~~~\sm{d\La{(}g(\Phi^{-1}d\Phi(\cdot, \cdot), \nabla \lh)\La{)}(X,Y,Z)=
g(d(\Phi^{-1}d\Phi)(X,Y,Z), \nabla \lh)+\cyclic_{\ti{X,Y,Z}}
Hess\, \lh (X,\Phi^{-1}d\Phi(Y,Z))}$
 what
proves the Proposition.\qed\\[-3mm]
\begin{Pp} If $F:M\ra N$ is a non-$J$-holomorphic Cayley submanifold 
and $N$ is Ricci-flat then $(1.6)-(1.7)$ holds.
\end{Pp}
\noindent
\em Proof. \em 
 To prove (1.6) we note that from (1.5) and (5.1) and Corollary 5.1
\[p_1(\sm{\bigwedge^2_-}NM)=p_1(\sm{\bigwedge^2_-}TM)+ 4( {\cal X}(M)
-{\cal X}(NM))=p_1(\sm{\bigwedge^2_-}TM)+2(p_1(NM)-p_1(M)).\]
Since $ R'(X,Y,{Z}, {W})=R^{\bot}(X,Y,\Phi(Z), \Phi(W))$ 
and $\Phi:(TM,\hat{g}, \lnabp{}\, )\ra (NM,g,\lnabo{})$ is a parallel
isometry along $M\sim {\cal C}$, then on this open set
$p_1(NM)=p_1(R^{\bot})=p_1(R')$, as forms defined by
the formulas (5.3). From (5.11), Proposition 5.5, and that
$d\La{(}d(g(\Phi^{-1}d\Phi(\cdot,\cdot),\nabla \lh))\La{)}=0$
we obtain (1.6)-(1.7). \qed 
\begin{Pp} If $(Y,Z)\ra g(\lnab{X}\Phi(Y), \Phi(Z))$ is symmetric
then $p_1(\sm{\bigwedge^2_-}NM)=p_1(\sm{\bigwedge^2_-}TM)$.
\end{Pp}
\noindent
\em Proof. \em From Proposition 5.2,
 $R^{\bot}(X,Y,\Phi(\hat{Z}), \Phi(\hat{W}))= R^M(X,Y,Z,W)$
and so the characteristic classes induced by $(R^{\bot}, g)$ 
are the same has the ones
induced by $(R^M,g)$. We could also check directly from 
$\Phi_ABC=\ha \lh_A g(B,C)$
that all terms of $\eta$ in (1.6) i.e (5.38)-(5.39) vanish. \qed.
\subsection{Homogeneous complex points}
Let us assume that $F:M\ra N$ is a compact 
Cayley submanifold, and let $\eta$ be the 3-form
on $M\sim {\cal C}$ defined as in (1.7).
Since   ${\cal C}$ is the zero set of
$\sin^2\theta$, that has only zeros of finite order, this
set as some regularity. Indeed, 
by the 
 Malgrange's preparation theorem for smooth functions
with zeros of finite order, locally we can find 
a coordinate chart $x=(x',x_4)$ onto an open set $U$ of $\RR^4$
 such that $\sin^2\theta$ can be written as
$h(x)\La{(} \sum_{0\leq a\leq k-1}
w_{a}(x')x_4^a + x_4^k\La{)},~$
where $h$ never vanish on  $U$ and $w_{a}$ vanish of
order $k-a$ at $0$.  Thus, the zero set of $\sin\theta$ can be 
locally parametrised as
$\Sigma=\{x=(x',x_4): \sum_{0\leq a \leq k-1}
w_{a}(x')x_4^a +x_4^k=0\}$. This set represents the zeros of a polinomial
function on the variable $x_4$, with coeficients on
the variable $x'$,  so it is, in general, still quite complicate
to handle. 
A simpler case is when we have a polinomial function of the type
$(x_{d+1}^2+\ldots  +x_4^2)^{k'} =0$, as is it is the case  with $k'=1$, of 
$x$ a 
Farmi coordinate chart of a submanifold $\Sigma$ of dimension $d$.
\\[-2mm]

Assume now $f$ is a nonnegative  continuous function 
defined on a open set  $V$ of $M$ containing $\Sigma=f^{-1}(0)$ and
smooth on $V\sim \Sigma$.
 For each $\epsilon>0$ sufficiently small define
\[V_f(\Sigma, \epsilon)=\{q\in M: f(q)<\epsilon\},~~~~~~~~~~
C_f(\Sigma, \epsilon)= \{q\in M: f(q) =\epsilon\}\]
For a dense set of regular values $\epsilon$, 
 $C_f(\Sigma, \epsilon)$ is a smooth hypersurface and  is the boundary of
$cl(V_f(\Sigma, \epsilon))$ and
for each $q\in C_f(\Sigma, \epsilon)$, $T_qC_f(\Sigma, \epsilon)=
[\nabla f(q)]^{\bot}$.
If the decreasing sequence 
$d_M(V_f(\Sigma,\epsilon))$ converges to $0$ 
when $\epsilon\ra 0$, where $d_M$ is the Lebesgue measure of $M$,  then
\begin{eqnarray}
\int_{M}d\la{(}\eta(\Phi)\la{)}
=\lim_{\epsilon\ra 0}\int_{M\sim V_f(\Sigma,\epsilon)}
d(\eta(\Phi))=-\lim_{\epsilon\ra 0}
\int_{C_f(\Sigma, \epsilon)}{\eta}(\Phi).
\end{eqnarray}

In case of $\Sigma$ is a smooth hypersurface of $M$,
$f$ is smooth on $V$, and $\nabla f$ does not vanish on
$\Sigma$, then for $\epsilon$ sufficiently
small,  $C_f(\Sigma,\epsilon)$ is connected,
diffeomorphic to $\Sigma$, and converges (in the 
Lebesgue sense) to $\Sigma$ and $M\sim V_f(\Sigma, \epsilon)$
to $M$ when $\epsilon\ra 0$. To see this let
$\rho:M\ra [0,1]$ be a smooth map s.t.
$\rho$ values 1 on $V_f(\Sigma,2r)$ and zero away from $V_f(\Sigma,3r)$
for $r$ sufficiently small, and let $\xi_t:M\ra M$ be the
one parameter family of diffeomorphisms
generated by the vector field globally defined on $M$,
 $X_f=\rho\frac{\nabla f}{\|\nabla f\|^2}$.
Then, there exist $\epsilon_0, \delta>0$ such that $\forall |t|<\epsilon_0$
and $q\!\in V_f(\Sigma,\delta)$,  $\xi_t(q)\!\in V_f(\Sigma,2r)$, 
and so
$\frac{\partial}{\partial t}f(\xi_t(q)) =
df(\xi_t(q))(\frac{\partial}{\partial t}\xi_t(q))
=1$~(\cite{[M]}). That is $f(\xi_t(q))=t +f(\xi_0(q))=t +f(q)$. In particular
$\forall 0<\!\epsilon\! <\epsilon_0$, $p\in \Sigma$, 
and $q\in C_f(\Sigma, \epsilon)$,
$f(\xi_{\epsilon}(p))=\epsilon$ and $f(\xi_{-\epsilon}(q))=0$. 
This means that
${\xi_{\epsilon}}_{|\Sigma}:\!\Sigma\!\ra\! C_f(\Sigma,\epsilon)$ is a
diffeomorphism with inverse $\xi_{-\epsilon}$. 
Let $\vartheta(\epsilon)$  and $ \tilde{\vartheta}(\epsilon)$ be the 
coefficients of dilatation of $\xi_{\epsilon}$ and 
$\xi_{\epsilon}|_{\Sigma}$. From $\xi_{0}(q)\!=\!q$, $\forall q$, we
easily  see that both  $\vartheta(\epsilon)(q)\!\ra\! 1$,
$\tilde{\vartheta}(\epsilon)(p)\!\ra\! 1$,
when $\epsilon\!\ra\! 0$. Moreover
if $\eta(\Phi)$ can be defined as an $L^1$-form along $\Sigma$, then
$(5.68) =-\int_{\Sigma}
\eta(\Phi).$ Unfortunately the case $\Sigma$ 
a hypersurface
is the least interesting, for, non $J$-complex Cayley submanifolds
of $\RR^8$ cannot have ${\cal C}$ as an analytic hypersurface \cite{[H-L1]}.
\\[-1mm]

A key example is  of $f= \sigma$ the
 intrinsic distance function to a smooth 
submanifold  $\Sigma$ of dimension $d$,
$\sigma(q)=d(q, \Sigma)=\inf_{p\in \Sigma}d(q,p)$. In this case,
$\nabla f$ is not well defined at each point $p\in\Sigma$, but
$\|\nabla f\|=1$, everywhere.
In fact $\nabla f$ it is multivalued, with
sublimits all unit normal vectors to $\Sigma$ in $M$. 
Nevertheless 
the flow can be smoothly extended to $\Sigma$, in all directions
of $T_p\Sigma^{\bot}$.  We explain as follows.
Let $N\Sigma$ denote the total space of the normal bundle
of $\Sigma$ in $TM$ and $N^1\Sigma$ the 
spherical subbundle of the unit orthogonal vectors. 
For each $\epsilon>0$ let
\[G_{\epsilon}=\{(p,w)\in N\Sigma: p\in \Sigma, w\in T_p\Sigma^{\bot}, 
\|w\|< {\epsilon}\},~~
C_{\epsilon}=\{(p,w)\in N\Sigma: p\in \Sigma, w\in T_p\Sigma^{\bot}, 
\|w\|={\epsilon}\}.\]
For $0< \epsilon\leq \epsilon_0$, with $\epsilon_0$
sufficiently small, 
the restriction of the exponential map of
$M$,  $exp:G_{\epsilon}\ra M$, $exp(p,w)=exp_p(w)$  
defines a diffeomorphism onto
$V_{\sigma}(\Sigma,{\epsilon})$ and
$exp(C_{\epsilon})=C_{\sigma}(\Sigma,{\epsilon})$  is its boundary.
For each  $w\in T_p\Sigma^{\bot}$,
$\gamma_{(p,w)}({\epsilon})=exp_p({\epsilon}w)$ is
the  geodesic normal to
$\Sigma$, starting at $p$ with initial velocity $w\in T_p\Sigma^{\bot}$.
Thus, $s(p,w):=\sigma(exp(p,w))=\|w\|$, is just the Euclidean norm
in $T_p\Sigma^{\bot}$. Since $N\Sigma$ is the total space of a Riemannian
vector bundle, then it has a natural Riemannian structure
such that $\pi:N\Sigma\ra \Sigma$ is a Riemannian submersion.
The volume element $Vol_{N\Sigma}$ for such metric satisfies
$Vol_{N\Sigma}(p,w)=Vol_{\Sigma}(p)\wedge ds(p,w)$ and
$Vol_{C_{\epsilon}}(p,w)= Vol_{\Sigma}(p)
\wedge Vol_{S(p,{\epsilon})}(w)$, where ${\epsilon}=\|w\|$,
and $S(p,{\epsilon})$ is the sphere of $T_p\Sigma^{\bot}$ of radius 
${\epsilon}$. For each $u\in N^1\Sigma_p$,  
$\vartheta_{(p,u)}({\epsilon})= \langle Vol_{N\Sigma}(p,{\epsilon}u),
exp^{*}Vol_M(p,{\epsilon}u)\rangle$ is the coefficient of dilatation
that  measures the volume
distortion by $exp$ in the direction $u$. It satisfies 
$\vartheta_{(p,u)}(0)=1$.
We recall the following (see \cite{[G2]}):
(1) $\nu(q)=\nabla \sigma(q)$
is the unit outward of $C(\Sigma,{\epsilon})$,
(2) $\nu(\gamma_{(p,u)}({\epsilon}))=  \gamma'_{(p,u)}({\epsilon})$,
(3) $ds\wedge *ds$ and $d\sigma\wedge *d\sigma$ are
the volume elements of $N\Sigma$ and $M$ respectively.
(4) $*ds$ and $*d\sigma$ are
the volume elements of each hypersurface
$C_{\epsilon}$ of $N\Sigma$ and $C(\Sigma,{\epsilon})$ of $M$ 
respectively, and $exp^*(*d\sigma)(p,w)
=\vartheta_{(p,\frac{w}{{\epsilon}})}({\epsilon})(*ds)(p,w)$,
where ${\epsilon}=\|w\|$. 
The set of sublimits of $\nabla \sigma$ at a point 
$p\in \Sigma$ is the
entire sphere $S(p,1)$ of $T_p\Sigma^{\bot}$ and by (1) and (2) 
for each $u\in N^1\Sigma$,  
$\gamma_{(p,u)}(t)$ is an integral curve of $\nabla \sigma$
smoothly extended at $t=0$ by $p$ and initial velocity $u$.
The map $\xi:N\Sigma\ra V$, $\xi(p,tu)=\gamma_{(p,u)}(t)$, 
can be seen as the flow of $\nabla \sigma$, a vector field
multivalued at $\Sigma$. \\[2mm]
This example motivates the following.
We consider functions $f:V\ra \RR^+_0$
satisfying the following conditions (E-1) and (E-2), that generalizes
the case of $f=\sigma$.
Let $X_f=\frac{\nabla f}
{\|\nabla f\|^2}$ and $\Sigma$ a smooth closed submanifold
of dimension $d$.\\[3mm]
{\bf (E-1)} \em  $f$ is a nonnegative  continuous function 
with zero set $\Sigma$, 
smooth on $V\sim \Sigma$ and with  $\|\nabla f\|$ defined
 $\forall p\in \Sigma$, 
giving a positive function of class $C^{\mu}$ on $M$, and such that
$\{ u\in T_pM: ~u\mbox{~is~a~sublimit~of~}\frac{\nabla f}
{\|\nabla f\|}\mbox{~at~}p ~\}= N^1\Sigma_p.$ \em \\[4mm]
Set for $p\in \Sigma$, $c(p)= \lim_{q\ra p} \|\nabla f(q)\|>0$.
We are considering the sublimits  defined through orthogonal 
curves to $\Sigma$, 
$\rho:[0, 1]\ra M$, such that $\rho(]0,1])\subset V\sim \Sigma$ and
$\rho(0)\in \Sigma$, and exist $u=\lim_{t\ra 0}\frac{\nabla f}
{\|\nabla f\|}(\rho(t))\in N^1\Sigma_p$. The set $E\Sigma_p=
\{\frac{u}{c(p)}: u\in N^1\Sigma_p\}$
is just the set  of sublimits of $X_f$  at $p$.
 An integral curve $\gamma:]0,b[\ra V\sim \Sigma$ of $X_f$ has an end 
point at $0$ converging  to $\Sigma$ with initial
velocity $\frac{u}{c(p)}$ where $u\in N\Sigma_p'$, 
 if $\exists \lim_{t\ra 0^+} \gamma(t)=p$ and 
$\exists \lim_{t\ra 0^+} \gamma'(t)=\frac{u}{c(p)}$.\\[4mm]
{\bf (E-2)} \em  $X_f$ 
has an extensible  flow to 
$\Sigma$, i.e  $\forall (p,u)\in N^1\Sigma$ there exist a 
$C^{\mu +1}$ curve $\gamma_{(p,u)}(t)$, defined $\forall t\in[0, t_0]$, 
smooth for $t>0$, that satisfies:\\[1mm]
$(a)$ for $t>0$, $\gamma_{(p,u)}(t)$ is an integral curve of 
$X_f$ on $V\!\sim\! \Sigma$, and
 $\gamma_{(p,u)}(0)=p$,  $\gamma_{(p,u)}'(0)=\frac{u}{c(p)}$.\\[1mm]
$(b)$ The flow at $\Sigma$, $\xi: N\Sigma\ra V $, 
defined for  $(p,w)$ with $\|w\|< t_0$, by $\xi(p,0)=p$,
$\xi(p,w)=\gamma_{(p,u)}(\epsilon)$, where $u=\frac{w}{\|w\|}$ 
and $\epsilon=\|w\|$, 
is a diffeomorphism of class $C^{\mu+1}$. \em \\[4mm]
So we have for $0\leq \epsilon<t_0$,
$\xi_{\epsilon}:C_{\epsilon}\ra C_f(\Sigma,\epsilon)$. The coefficient of
dilatation of  $\xi $ at $(p,\epsilon u)$,
$\vartheta_{(p,u)}(\epsilon)=\langle
Vol_{N\Sigma}(p,\epsilon u), (\xi^*Vol_M)(p,\epsilon u)\rangle$,
satisfies $\vartheta_{(p,u)}(0)=\frac{1}{c(p)^{4-d}}$, as we will see below.
The volume element of $C_f(\Sigma,\epsilon)$
is $\frac{*df}{\|\nabla f\|}$, for  $\frac{\nabla f}{\|\nabla f\|}$ is
the outward unit.
Since $\frac{df}{\|\nabla f\|}\wedge \frac{*df}{\|\nabla f\|}$ is the volume
element of $M$ and $\xi^*(\frac{df}{\|\nabla f\|})(p,\epsilon u)(0,u)=
\frac{df(\gamma'_{(p,u)}(\epsilon))}{\|\nabla f((p,\epsilon u)\|}
=\frac{1}{\|\nabla f(\xi(p,\epsilon u))\|}$, then on
 $T_{(p,\epsilon u)}C_{\epsilon}$ we have
$\frac{1}{\|\nabla f(p, \epsilon u)\|}\xi^{*}( \frac{*df}{\|\nabla f\|})
(p, \epsilon u)=\vartheta_{(p,u)}(\epsilon)Vol_{C_{\epsilon}}=
\vartheta_{(p,u)}(\epsilon)*ds$. Therefore, 
$\tilde{\vartheta}_{(p,u)}(\epsilon):=\|\nabla f(\xi(p,\epsilon u))\|
\vartheta_{(p,u)}(\epsilon)$ is the coefficient of dilatation 
of $\xi$ restricted to $C_{\epsilon}$, and
$\lim_{\epsilon\ra 0}\tilde{\vartheta}_{(p,u)}(\epsilon)=
{c(p)^{d-3}}$. \\[2mm]
From (E-2) we have a coordinate system of class $C^{\mu +1}$
of Farmi-type. Let
$O$ be an open set of $\Sigma$ where a coordinate system exist and
a d.o.n. frame $E_{d+1}, \ldots, E_4$ of $N\Sigma$. 
Then for each $w\in N\Sigma_p$, $w=\sum_{ d+1\leq i\leq 4}t_iE_i(p)$.
Define on $V'$ the image by $\xi$ of the restriction to $O$,
\[ \begin{array}{rccccc}
   x:&V' & \ra &N\Sigma & \ra & \RR^4\\
     &\xi(p,w)&\ra &(p,\sum_{ d+1\leq i\leq 4}t_iE_i(p)) & \ra &(y(p), t_{d+1},
\ldots, t_4)\end{array}\]
Thus, for $q=\xi(p, \epsilon u)\in C_f(\Sigma, \epsilon)$, $\|u(p)\|=1$,
$u(p)=\sum_{1+d\leq i\leq 4}u_i(p)E_i(p)$
\begin{eqnarray}
\partial_i f(q) &=&
\left\{\begin{array}{ll}
0  &\forall i\leq d\\
u_i(p)=\frac{x_i(q)}{f(q)}&\forall i\geq d+1
\end{array}\right.\\[3mm]
Hess\, f(q) (\partial_i, \partial_j)&=&
\left\{\begin{array}{l}
-\sum_{s\geq d+1}\Gamma^s_{ij}(q) u_s(p)~~~~~~ \mbox{if}~i\leq d
\mbox{~or~}j\leq d\\
\frac{1}{\epsilon}(\delta_{ij}-u_i(p)u_j(p))
-\sum_{s\geq d+1}\Gamma^s_{ij}(q)u_s(p)~~~ \mbox{if~}i,j\geq 1+d
\end{array}\right.
\end{eqnarray}
In fact for $s\leq \mu $,  
it is defined a tensor $T^{s}\in C^{\mu-s}(\pi^{-1}\bigotimes^{s}
TM^*)$, where $\pi:N^1\Sigma\ra O\subset\Sigma$, 
and s.t. $\exists \lim_{\epsilon\ra 0} \epsilon^{s-1}
\lnab{\partial_{i_{s}}, \ldots,
\partial_{i_1}}^{s}f\,\sm{(\xi(p, \epsilon u))}= 
:T^{s}(p,u)(\partial_{i_1}(p), \ldots,\partial_{i_{s}}(p))$
where $\lnab{X}f=df(X)$,
\begin{equation}
\lnab{X_k, \ldots,X_1}^kf=\lnab{X_k}(\lnab{X_{k-1},
\ldots, X_1}^{k-1}f)-\sm{\sum_{k-1\geq i\geq 1}}\lnab{X_{k-1}, \ldots, 
\lnab{X_k}X_i, \ldots, X_1 }^{k-1}f.
\end{equation}
Recall that $T_{(p,0)}N\Sigma= T_p\Sigma\times T_p\Sigma^{\bot}$, and
from $\xi(p,0)=\gamma_{(p,u)}(0)=p$
we get $\forall X\in T_p\Sigma$, $d\xi_{(p,0)}(X,0)=X$. Now,
if $0\neq h\in T_p \Sigma^{\bot}$, the curve $\tau(s)=\xi(p, s h)=\gamma_{(
p,\frac{h}{\|h\|})}(\|h\| s)$
satisfies  $\tau'(0)=d\xi_{(p,0)}(0,h)=\|h\|\gamma'_{(p,\frac{h}{\|h\|})}(0)=
\frac{h}{c(p)}$. Thus, $d\xi_{(p,0)}(X,h)=X+\frac{h}{c(p)}$,
and so $\vartheta_{(p,u)}(0)= (c(p))^{d-4}$.
So we conclude:
\begin{Pp} If a continuous function $f$ satisfies (E-1) and (E-2) then for each
$p\in \Sigma$ there exist a $C^{\mu +1}$ coordinate chart $x$ of $M$, 
adapted  to $\Sigma$,  and such that $f^2= x_{d+1}^2
+ \ldots + x_{4}^2$ and for $i\geq d+1$, $E_i(p)
=c(p)\frac{\partial}{\partial x_i}(p)$ is an o.n.\ basis of 
$T_p\Sigma^{\bot}$.\\[1mm]
\end{Pp}
\noindent

\em From now on we assume $\sin\theta=f^r$, with
$\Sigma=f^{-1}={\cal C}$ and $f$
 satisfying conditions (E-1) and (E-2), with $\mu\geq r+1$.~ \em
Set $\tilde{\Phi}=\frac{\Phi}{\|\Phi\|}=\frac{\Phi}{\sin\theta}$.
Then $\tilde{\Phi}$ is an isometry and $\tilde{\Phi}^{-1}=\tilde{\Xi}
:=\frac{\Xi}{\|\Xi\|}$.\\ 

Using  a $C^{\mu +1}$ coordinate chart with $\mu \geq r$,
 $\Phi$ has a zero of order $r$ at $0$ iff $\frac{\Phi(x)}{\|x\|^{r-1}}
\ra 0$  and $ \frac{\Phi(x)}{\|x\|^{r}}$ does not converge to $0$
when $x\ra 0$, in other words, $\Phi$ ( or equivalently $\|\Phi \|$) is
an $O(\|x\|^r)$. This is  
equivalent to $D^s\Phi(0)=0$ $\forall s\leq r-1$
and $D^r\Phi(0)\neq 0$. Thus, at all points  $p\in\Sigma$, $r$ is the 
order of the zero of $\Phi$ at $p$.
Fix $p\in \Sigma$ and $y$ a coordinate
system of $\Sigma$ with $y(p)=0$, and consider $x$ the corresponding
Farmi coordinate system.
Let $V''=x(V')$ open set of $\RR^4$
 and $dx^{-1}:\RR^4_{V''}\ra T{V''}$, and an isomorphism
$\tau:N{V''}\ra\RR^4_{V''}$. Then $P=\tau\circ \Phi_{x^{-1}}\circ
d(x^{-1})
:V''\ra L(\RR^{4}; \RR^{4})$
has at $0$ a zero of order $r$. Thus, for $v$ sufficiently close to $0$
\begin{equation}
P(v)=\frac{1}{r!}D^rP(0)(v)^r + \int_0^1\frac{(1-t)^r}{r!}
D^{r+1}P(tv)(v)^{r+1}dt
\end{equation}
where $D^rP(0)(v)^r= \tau_p \circ
\lnab{(d(x^{-1})(0)(v))^r}^r\Phi(p)\circ d(x^{-1})(0)$. This
term does not vanish for some $v$.
If we take $u=\sum_{i\geq d+1}t_iE_i\in N^1\Sigma_p$, $\epsilon$
small enough, and $v=(0,\epsilon t_{d+1}, \ldots ,\epsilon t_4)$,
then $x^{-1}(sv)=\gamma_{(p,u)}(s\epsilon)$ and so
$d(x^{-1})(tv)(v)=\epsilon \gamma_{(p,u)}'(t\epsilon )$,
 $d(x^{-1})(0)(v)=\epsilon c(p)^{-1}u$.
Therefore
\begin{equation}\begin{array}{l}
\Phi(\xi(p,\epsilon u))=\epsilon^r\LA{(}A(\xi(p, \epsilon u))
+\epsilon Q(\xi(p,\epsilon u))\LA{)}\\[1mm]
A(\xi(p, \epsilon u))=\frac{1}{r!c(p)^r}\tau^{-1}_{\xi(p, \epsilon u)}
\circ\tau_p\lnab{u^r}^r\Phi(p)\circ dx^{-1}(0)\circ
dx(\xi(p,\epsilon u))\mbox{~~~of~class~}C^{\mu}\end{array}
\end{equation}
for some $Q$ of class $C^{\mu-r}$. Let $\lnab{X_k,\ldots, X_1}^k\Phi(p)$
be defined as (5.71).
\begin{Lm} 
$\forall X_i\in T_p\Sigma$, $\lnab{X_r,\ldots, X_1}^r\Phi(p)=0$.
\end{Lm}
\noindent
\em Proof. \em  If $X$ is
a vector of $T_p\Sigma$, we  can assume $X(p)=\gamma'(0)$ for some
curve $\gamma(t)$ on $\Sigma$. Then $\lnab{X}\Phi(p)=\lnab{\frac{d}{dt}}
\gamma^{-1}\Phi(0)$. But $\gamma^{-1}\Phi$ is constantly equal to $0$.
So if $X$ is a vector field of $\Sigma$, $\lnab{X}\Phi$ vanish along $\Sigma$.
Moreover,
 $\lnab{X_s,\ldots,X_1}^s\Phi(p)=0$ for any $s\leq r-1$ and $X_i\in T_pM$.
If $r\geq 2$ and
$X,Y\in T_p\Sigma$, we extend $X$ to a vector field along $\Sigma$,
and so, $\lnab{Y,X}^2\Phi(p)=\lnab{Y}(\lnab{X}\Phi)(p)=0$,
that is  $\lnab{Y,X}^2\Phi$ vanish along
$\Sigma$. This implies that if $r\geq 3$, for $X,Y,Z$ vector fields
of $\Sigma$, 
$\lnab{Z,Y,X}^3\Phi(p)=\lnab{Z}(\lnab{Y,X}^2\Phi)(p)=0$.
 The same for any $r$.
\qed\\[3mm]
A tensor $\varsigma$ in $C^{\infty}
( N\Sigma^* \otimes (T^*\Sigma\otimes N\Sigma))$
is defined by: if  $u\in T_p\Sigma^{\bot}$,
$X\in T_p\Sigma$ then
$\varsigma(u)(X)=(\lnab{u}\tilde{X})^{\bot}$, where $\tilde{X}$ is any vector
field of $M$ with $\tilde{X}_p=X$.
\begin{Pp}   $\forall (p,u)\in N^1\Sigma$, 
 $\exists \lim_{\epsilon\ra 0}\tilde{\Phi}(\xi(p, \epsilon u))=
\frac{1}{r!c(p)^r}\lnab{u^r}^r\Phi(p)=:\tUpsilon(p,u)$,
is an isometry. Moreover  
 $\forall X\in T_p \Sigma$, $Y_i\in T_pM$,
$\lnab{Y_{r-1}, \ldots, Y_1,X}^r\Phi(p)
=\lnab{Y_{r-1}, \ldots, Y_s,X, Y_{s+1}, \ldots, Y_1}^r\Phi(p)=0$,
and $\lnab{u^r}^r(\lnab{X}\Phi)=\lnab{u^r, X}^{r+1}\Phi(p)
+r\lnab{u^{r-1},\varsigma(u)(X)}^r\Phi$.
\end{Pp}
\noindent
\em Proof. \em $\frac{\Phi((\xi(p,\epsilon u))}{\epsilon^r}=
\tilde{\Phi}(\xi(p,\epsilon u)$ is an isometry.
 By (5.73), making $\epsilon \ra 0$ we conclude that
$\forall u\in N^1\Sigma_p,$ 
 $\tUpsilon(p,u)=\frac{1}{r!c(p)^r}\lnab{u^r}^r\Phi(p)$ is an isometry. 
Extend $X$ to a local section of $T\Sigma$, and then
extend $X$  and $Y_i$ to
local sections of $TM$.
Recall that $\lnab{}^s\Phi$ vanish along $\Sigma,$
$\forall s\leq r-1$.
Then
$\lnab{X}\La{(}\lnab{Y_{r-1}, \ldots, Y_1}^{(r-1)}\Phi\La{)}(p)=0$, 
and
\begin{eqnarray*}
\sm{0}
&=&\sm{\lnab{X}(\lnab{Y_{r-1}, \ldots, Y_1}^{(r-1)}\Phi)(p)
=\lnab{X}(\lnab{Y_{r-1}}(
\lnab{Y_{r-2}, \ldots, Y_1}^{(r-2)}\Phi))(p)}\\
&=&\sm{\lnab{Y_{r-1}}{(}\lnab{X}(\lnab{Y_{r-2}, \ldots, Y_1}
^{(r-2)}\Phi){)}(p)
+\bar{R}(Y_{r-1},X)(\lnab{Y_{r-2},\ldots, Y_1 }^{(r-2)}\Phi(p))}\\
&=&\sm{\lnab{Y_{r-1}}(\lnab{X}(\lnab{Y_{r-2}, \ldots, Y_1}
^{(r-2)}\Phi))(p)=\lnab{Y_{r-1}}(\lnab{X}(\lnab{Y_{r-2}}
(\lnab{Y_{r-3}, \ldots, Y_1}^{(r-3)}\Phi)))(p)}\\
&=&\sm{\lnab{Y_{r-1}}\La{(}\lnab{Y_{r-2}}
(\lnab{X}(\lnab{Y_{r-3}, \ldots , Y_1}^{(r-3)}\Phi))
+ \bar{R}(Y_{r-2},X)(\lnab{Y_{r-3}, \ldots, Y_1}^{(r-3)}\Phi)\La{)}(p)}\\
&=& \sm{\lnab{Y_{r-1}}\La{(}\lnab{Y_{r-2}}
(\lnab{X}(\lnab{Y_{r-3}, \ldots, Y_1}^{(r-3)}\Phi))\La{)}(p)
+\bar{R}(Y_{r-2},X)(\lnab{Y_{r-1}}(\lnab{Y_{(r-3)}, \ldots, Y_1}
^{(r-3)}\Phi)(p))}\\
&=&\sm{\lnab{Y_{r-1}}\La{(}\lnab{Y_{r-2}}
(\lnab{X}(\lnab{Y_{r-3},\ldots, Y_1}^{(r-3)}\Phi))(p)}
\end{eqnarray*}
and successively, 
$0=\lnab{X}(\lnab{Y_{r-1}, \ldots, Y_1}^{(r-1)}\Phi)(p)
=\lnab{X,Y_{r-1},\ldots, Y_1}^{r}\Phi(p)
=\lnab{Y_{r-1}, \ldots, Y_s, X, Y_{s+1}, \ldots, Y_1}^{(r-1)}
\Phi(p)= \lnab{Y_{r-1}, \ldots, Y_1}^r(\lnab{X}\Phi)(p)
=\lnab{Y_{r-1}, \ldots, Y_1,X}^{r}\Phi(p).$
Thus, if $r=2$ then $\lnab{u^2}^2(\lnab{X}\Phi)(p)=\lnab{u^2, X}^3\Phi(p)
+ 2\lnab{u, \lnab{u}X}^2\Phi(p)$, and 
$\lnab{u, \lnab{u}X}^2\Phi(p)=\lnab{u, (\lnab{u}X)^{\bot}}^2\Phi(p)
=\lnab{u, \varsigma(u)(X)}^2\Phi(p)$. The proof for $r\geq 3$ is
similar, slightly more complicate.
\qed \\
\begin{Pp} $\forall X$ vector field on $M$ and 
 $\forall (p,u)\in N^1\Sigma$, set $X^{\bot_u}=X_p-g(X_p,u)u$.
Then
$\exists \lim_{\epsilon \ra 0}\epsilon
\lnab{X}\tilde{\Phi}(\xi(p, \epsilon u))=\frac{1}{(r-1)!c(p)^{r-1}}
\lnab{u^{r-1},X^{\bot_u}}^r\Phi(p)
=:\tilde{\Psi}(p,u)(X^{\bot_u}).$ 
If $d\geq 1$ and $X_p\in T_p\Sigma$ and
for $q$ near $p$,
$X(q)\bot \nabla f(q)$ 
then 
$\exists\lim_{\epsilon \ra 0}
\lnab{X}\tilde{\Phi}(\xi(p, \epsilon u))=\frac{1}{r!c(p)^{r}}
\La{(}\lnab{u^{r},X}^{r+1}\Phi(p) 
+r \lnab{u^{r-1},\varsigma(u)(X)}^{r}\Phi(p)\La{)}
=:\tilde{G}(p,u)(X_p).$
\end{Pp}
\noindent
\em Proof. \em 
There exist some $\tilde{Q}$ of class $C^{\mu-r-1}$, s.t.\\[-2mm]
\begin{equation}
\begin{array}{l}
\lnab{X}\Phi(\xi(p,\epsilon u))=\epsilon^{r-1}\LA{(}B(X,\xi(p,\epsilon u))
+\epsilon \tilde{Q}(\xi(p,\epsilon u))\LA{)}\\[1mm]
B(X,\xi(p,\epsilon u))=\frac{1}{(r-1)!c(p)^{(r-1)}}
\tau^{-1}_{\xi(p, \epsilon u)}\circ\tau_p
\lnab{u^{(r-1)},X}^r\Phi(p)
\circ dx^{-1}(0)\circ dx(\xi(p,\epsilon u)).
\end{array}
\end{equation}
Thus\\[-7mm]
\begin{eqnarray*}
\lim_{\epsilon\ra 0}\epsilon\lnab{X}\tilde{\Phi}(\xi(p,\epsilon u))&=&
\lim_{\epsilon\ra 0}\epsilon^{-r+1}\lnab{X}\Phi(\xi(p,\epsilon u))-
rg(\nabla f(\xi(p,\epsilon u)),X)\tilde{\Phi}
(\xi(p,\epsilon u))\non\\
&=& \lim_{\epsilon\ra 0}B(X,\xi(p,\epsilon u))
+ \epsilon \tilde{Q}(\xi(p,\epsilon u))-
rg(\nabla f(\xi(p,\epsilon u)),X)\tilde{\Phi}
(\xi(p,\epsilon u))\non\\
&=&\frac{1}{(r-1)!c(p)^{(r-1)}}\lnab{u^{(r-1)},X}^r\Phi(p)
-\frac{1}{(r-1)!c(p)^{r-1}}g(u,X)\lnab{u^{r}}^r\Phi(p)\non\\
&=&\frac{1}{(r-1)!c(p)^{(r-1)}}\lnab{u^{(r-1)},X^{\bot_u}}^r\Phi(p). 
\end{eqnarray*}
If $X_p\in T_p\Sigma$ then by Prop. 5.9, $B(X,\xi(p,\epsilon u))=0$,
and $\lnab{X}\Phi(\xi(p,\epsilon u))=\epsilon^r\tilde{Q}(\xi(p,\epsilon
u))$, and  $\lim_{\epsilon\ra 0}\tilde{Q}(\xi(p,\epsilon u))
=\frac{1}{r!c(p)^r}\lnab{u^r}^r(\lnab{X}\Phi)(p)$. Set $X=\partial_i$ where
$i\leq d$. By (5.69) those vector fields span exactly the ones that are
orthogonal to $\nabla f$ at $\xi(p,\epsilon u)$.
  Note that $X_p= \frac{\partial}{\partial y_i}(p)$.
Hence, by Prop.5.9,~ 
$\lim_{\epsilon \ra 0}\lnab{X}\tilde{\Phi}(\xi(p, \epsilon u))
=\lim_{\epsilon \ra 0}\epsilon^{-r}\lnab{X}\Phi(\xi(p, \epsilon u))
=\frac{1}{r!c(p)^r}
\lnab{u^r}^{r}(\lnab{X}\Phi)(p)=\tilde{G}(p,u)(X_p).$\qed\\[4mm]
 Therefore if  $X_p\bot u$,\\[-4mm]
\begin{equation}
\lim_{\epsilon \ra 0}
\epsilon\tilde{\Phi}^{-1}\lnab{X}\tilde{\Phi}
(\xi(p,\epsilon u))=\tUpsilon(p,u)^{-1}\tPsi(p,u)(X_p)
=rc(p)\LA{(}\lnab{u^{r}}^r\Phi(p)\LA{)}^{-1}\!\!\!\circ
\left(\lnab{u^{(r-1)},X}^r\Phi(p)\right),
\end{equation}
and if $X\in T_p\Sigma$ and $X(q)\bot \nabla f(q)$ for $q$ near $p$,
\\[-4mm]
\begin{equation}
\lim_{\epsilon \ra 0}
\tilde{\Phi}^{-1}\lnab{X}\tilde{\Phi}
(\xi(p,\epsilon u))=\tUpsilon(p,u)^{-1}\tilde{G}(p,u)(X_p)
= \LA{(}\lnab{u^{r}}^r\Phi(p)\LA{)}^{-1}\!\!\!\circ
\left(\lnab{u^{r},X}^{r+1}\Phi(p)
+ r\lnab{u^{r-1},\varsigma(u)(X)}^{r}\Phi(p)\right).
\\[5mm]
\end{equation}
We can write 
 $\Phi=\sin\theta\, \tilde{\Phi}$ where $\tilde{\Phi}:TM\ra NM$
is an isometry, away from $\Sigma$. More generally, 
\begin{Lm} If  $V$ is an open set  containing $\Sigma$, and on 
$V\sim \Sigma$ ~$\Phi= \zeta^r\tilde{\Phi}$,  where $\tilde{\Phi}$
is a section of $TM^*\otimes NM$ defined on $V\sim \Sigma$ and
$\zeta:V\sim \Sigma \ra \RR^+_{0}$ some function , then where 
 $\zeta$ and $\tilde{\Phi}$ are  differentiable and 
do not vanish, 
$\eta(\Phi)=\eta(\tilde{\Phi})$
\end{Lm}
\noindent
\em Proof. \em  $\tilde{\Phi}:TM\ra NM$ is a conformal morphism. 
So (5.13) holds
for $\tilde{\Phi}$ with $\tilde{\lh}=\log(\frac{\|\tilde{\Phi}\|^2}{4})$.
and
\begin{eqnarray*}
&&\Phi^{-1}\lnab{X}\Phi = rd\log \zeta(X) Id_{TM}+ \tilde{\Phi}^{-1}
\lnab{X}\tilde{\Phi}\\
&&\Phi(R^{\bot})= \tilde{\Phi}(R^{\bot}),~~~~~~~
b(\Phi(R^{\bot}))= b(\tilde{\Phi}(R^{\bot}))
=\tilde{\Phi}^{-1}d^2\tilde{\Phi}\\
&&d^2\Phi=\zeta^rd^2\tilde{\Phi},~~~~~~~~
d(\Phi^{-1}d\Phi)=d(\tilde{\Phi}^{-1}d\tilde{\Phi})
\end{eqnarray*}
The last two equalities are proved using the symmetry of $Hess\, \log \zeta$.
Now
\[\begin{array}{c}
\langle Id_{TM}, \tilde{\Phi}(R^{\bot})(Y,Z)\rangle
 =\sm{\sum_i} \zeta^{-r} h^{-1} g( \tilde{\Phi}(e_i), R^{\bot}(Y,Z)
( \tilde{\Phi}(e_i)))=0,\\[1mm]
\langle Id_{TM}, R^M(Y,Z)\rangle =
\sm{\sum_i} g(e_i, R^M(Y,Z)(e_i))=0,\\[1mm]
 [ (r d\log \zeta(\cdot)Id_{TM}+ \tilde{\Phi}^{-1}\lnab{}\tilde{\Phi}),
(rd\log \zeta(\cdot)Id+ \tilde{\Phi}^{-1}\lnab{}\tilde{\Phi}) ]
= [ \tilde{\Phi}^{-1}\lnab{}\tilde{\Phi},
\tilde{\Phi}^{-1}\lnab{}\tilde{\Phi} ].
\end{array}\]
Thus,
$\eta(\Phi)=
 \eta(\tilde{\Phi})-\frac{r}{12\epsilon}\langle (d\zeta(\cdot)Id_{TM})
\wedge [\tilde{\Phi}^{-1}\lnab{}\tilde{\Phi},
\tilde{\Phi}^{-1}\lnab{}\tilde{\Phi}]\rangle. 
~$ Moreover
\begin{eqnarray*}
\lefteqn{\langle (d\zeta(\cdot)Id_{TM})
\wedge [\tilde{\Phi}^{-1}\lnab{}\tilde{\Phi},
\tilde{\Phi}^{-1}\lnab{}\tilde{\Phi}]\rangle(X,Y,Z)=~~~~~~~}~~~~~~~\\
&=&
\cyclic_{\ti{X,Y,Z}}\zeta_X\langle Id, \tilde{\Phi}^{-1}\lnab{Y}\tilde{\Phi}
\circ\tilde{\Phi}^{-1}\lnab{Z}\tilde{\Phi}
-\tilde{\Phi}^{-1}\lnab{Z}\tilde{\Phi}\circ
\tilde{\Phi}^{-1}\lnab{Y}\tilde{\Phi}\rangle\\[-1mm]
&=&\cyclic_{\ti{X,Y,Z}}\sm{\sum_i}~
\zeta_X g(e_i, \tilde{\Phi}^{-1}\lnab{Y}\tilde{\Phi}
(\tilde{\Phi}^{-1}\lnab{Z}\tilde{\Phi}(e_i)))
-\zeta_Xg(e_i,\tilde{\Phi}^{-1}\lnab{Z}\tilde{\Phi}(
\tilde{\Phi}^{-1}\lnab{Y}\tilde{\Phi}))\\[-1mm]
&=&\cyclic_{\ti{X,Y,Z}}\sm{\sum_i}~
-\zeta_X g(\tilde{\Phi}^{-1}\lnab{Y}\tilde{\Phi}(e_i),
\tilde{\Phi}^{-1}\lnab{Z}\tilde{\Phi}(e_i)))+\zeta_X\tilde{\lh}_Yg(e_i,
\tilde{\Phi}^{-1}\lnab{Z}\tilde{\Phi}(e_i))\\[-1mm]
&&
\cyclic_{\ti{X,Y,Z}}\sm{\sum_i}~
+\zeta_Xg\tilde{\Phi}^{-1}\lnab{Z}\tilde{\Phi}(e_i),
\tilde{\Phi}^{-1}\lnab{Y}\tilde{\Phi}(e_i)))
-\zeta_X\tilde{\lh}_Zg(e_i,
\tilde{\Phi}^{-1}\lnab{Y}\tilde{\Phi}(e_i))\\[-1mm]
&=& 2\zeta_X\tilde{\lh}_Y\tilde{\lh}_Z-2\zeta_X\tilde{\lh}_Z\tilde{\lh}_Y=0
\end{eqnarray*}
Thus,
$\eta(\Phi)=\eta(\tilde{\Phi})$.\qed.\\[5mm]
\em Proof of Corollary 1.1. \em 
For simplicity of notation we assume $\Sigma_i= \Sigma$.
By (1.6) of Theorem 1.1 and by Lemma 5.8 we have
\[p_1(\sm{\bigwedge}^2_{-}NM)[M]-p_1(\sm{\bigwedge}^2_{-}TM)[M]
=\int_M d\eta(\Phi)
=-\lim_{\epsilon\ra 0}\int_{C_f(\Sigma, \epsilon)}\eta(\Phi)
=-\lim_{\epsilon\ra 0}\int_{C_f(\Sigma, \epsilon)}\eta(\tilde{\Phi})
\]
Now\\[-8mm] 
\begin{eqnarray*}
\lefteqn{\int_{C_f(\Sigma, \epsilon)}{\eta}(\tilde{\Phi}) = 
\int_{C(\Sigma, \epsilon)}\frac{1}{\|\nabla f\|^2}\langle {\eta}
(\tilde{\Phi})(q), * df(q)\rangle *df(q)
= \int_{C_{\epsilon}}\xi^*\La{(}\frac{1}{\|\nabla f\|^2}\langle {\eta}
(\tilde{\Phi}), * df \rangle *df\La{)}(p,w)}\non\\
&=& \int_{C_{\epsilon}}\frac{1}{\|\nabla f(\xi(p,w))\|}
\langle \eta(\tilde{\Phi})(\xi(p,w)), 
* df(\xi(p,w)))\rangle 
\vartheta_{(p,\frac{w}{\epsilon})}(\epsilon)(* ds)(p,w)\\
&=& \int_{\Sigma}\LA{(}\int_{S(p,\epsilon)} \frac{1}{\|\nabla f(\xi(p,w))\|}
\langle \eta(\tilde{\Phi})(\xi(p,w)), 
* df(\xi(p,w))\rangle \vartheta_{(p,\frac{w}{\epsilon})}(\epsilon)
d_{S(p, \epsilon)}(w)\LA{)} d_{\Sigma}(p)\\
&=&\int_{\Sigma}\LA{(}\int_{S(p,1)}\frac{1}{\|\nabla f(\xi(p,\epsilon u))\|}
\langle  \eta(\tilde{\Phi})(\xi(p,\epsilon u)),
* df(\xi(p,\epsilon u))\rangle 
\vartheta_{(p, u)}(\epsilon)\epsilon^{3-d}d_{S(p,1)}(u)\LA{)} d_{\Sigma}(p)
\\[-1mm]
\end{eqnarray*}
and $~  \frac{1}{\|\nabla f(\xi(p,\epsilon u))\|}
\langle  \eta(\tilde{\Phi})(\xi(p,\epsilon u)),
* df(\xi(p,\epsilon u))\rangle
=\eta(\tilde{\Phi})(\xi(p,\epsilon u))(e_2,e_3,e_4),~~$
where $e_i\in T_{(\xi(p,\epsilon u))}C_f(\Sigma, \epsilon)$
is  a d.o.n. frame.
We take $e_1:=\frac{\nabla f(\xi(p,\epsilon u))}
{\|\nabla f(\xi(p,\epsilon u)\|}$  
the outward unit of $C_f(\Sigma, \epsilon)$ at $\xi(p,\epsilon u)$, and so
 $e_2\wedge e_3\wedge e_4=\frac{*\nabla f(\xi(p,\epsilon u))}
{\|\nabla f(\xi(p,\epsilon u))\|}$, 
giving $e_i$ a d.o.n. basis of $T_{\xi(p,\epsilon u)}M$. 
Then
\begin{eqnarray}
\int_{C_f(\Sigma, \epsilon)}{\eta}(\tilde{\Phi}) =
\int_{\Sigma}\left(\int_{S(p,1)}{\eta}(\tilde{\Phi})
\left(\frac{*\nabla f(\xi(p, \epsilon u))}
{\|\nabla f(\xi(p,\epsilon u))\|}\right)
{\vartheta_{(p,u)}(\epsilon)\epsilon^{3-d}}
d_{S(p,1)}(u)\right) d_{\Sigma}(p)
\end{eqnarray}
But $\frac{\nabla f}{\|\nabla f\|^2}(\xi(p,\epsilon
u))=\gamma'_{(p,u)}(\epsilon)$ where $\gamma_{(p,u)} (\epsilon)=\xi(p, \epsilon u)$.
Then $\frac{*\nabla f(\xi(p,\epsilon u))}{\|\nabla f(\xi(p,\epsilon u))\|}$
 converges to $*u$ when $\epsilon\ra 0$. Recall that $\vartheta(p,u)(0)
=\frac{1}{c(p)^{4-d}}$. Now, from Propositions 5.9 and 5.10, 
\begin{eqnarray*}
\mbox{if~}d=2~~~~~&&\lefteqn{\lim_{\epsilon \ra 0}\epsilon\langle 
\tilde{\Phi}^{-1}\lnab{}\tilde{\Phi}
(\xi(p,\epsilon u))\wedge (\tilde{\Phi}(R^{\bot})+R^M)(\xi(p,\epsilon u))
\rangle=~~~~~~~~~~~~~~}~~~~~~~~~~~~~\non\\[-2mm]
&&~~~~~~~~=\langle \tUpsilon(p,u)^{-1}\tPsi(p,u) \wedge {(}
\tUpsilon(p,u)(R^{\bot}(p))+R^M(p))\rangle
~~~~~~~~
\\[1mm]
\mbox{if~}d=0~~~~~&&\lefteqn{\lim_{\epsilon \ra 0}\epsilon^{3}\langle 
\tilde{\Phi}^{-1}\lnab{}\tilde{\Phi}
(\xi(p,\epsilon u))\wedge [\tilde{\Phi}^{-1}\lnab{}\tilde{\Phi}
(\xi(p,\epsilon u)), \tilde{\Phi}^{-1}\lnab{}\tilde{\Phi}
(\xi(p,\epsilon u))]\rangle=~~~}\non\\[-2mm]
&&~~~~~~~~~= \langle \tUpsilon(p,u)^{-1}\tPsi(p,u) \wedge
[\tUpsilon(p,u)^{-1}\tPsi(p,u), \tUpsilon(p,u)^{-1}\tPsi(p,u)]\rangle
~~~~~~~~
\end{eqnarray*}
Now if $d\geq 1$,  
$\frac{*\nabla f(\xi(p,\epsilon u))}{\|\nabla f(\xi(p,\epsilon u))\|}
=X_1\wedge X_2\wedge X_3$ with $X_i\bot \nabla f(\xi(p,\epsilon u)$.
By (5.69), 
 we may take  $X_1, \ldots, X_d$ an
o.n. basis of $ span \{\partial_i,
i\leq d\}$. Note that  for $i\leq d$, $\lim_{\epsilon \ra 0}
\partial_i(\xi(p,\epsilon u))=\frac{\partial}{\partial y_i}(p)\in T_p\Sigma$.
This implies by Prop. 5.9, if $d=3$
\begin{eqnarray*}
\lefteqn{\lim_{\epsilon \ra 0}\epsilon^{3-d}\langle
\tilde{\Phi}^{-1}\lnab{}\tilde{\Phi}
(\xi(p,\epsilon u))\wedge [\tilde{\Phi}^{-1}\lnab{}\tilde{\Phi}
(\xi(p,\epsilon u)), \tilde{\Phi}^{-1}\lnab{}\tilde{\Phi}
(\xi(p,\epsilon u))]\rangle (\frac{*\nabla f(\xi(p,\epsilon u))}
{\|\nabla f(\xi(p,\epsilon u))})=~~~~~~~~~~~~~~~~~~}~~~~~~~~~~~~~~~~~\non\\
&=& \langle \tUpsilon(p,u)^{-1}\tilde{G}(p,u) \wedge
[\tUpsilon(p,u)^{-1}\tilde{G}(p,u), \tUpsilon(p,u)^{-1}\tilde{G}(p,u)]\rangle
~~~~~~~~~~~~ ~~~~~~~~
\end{eqnarray*}
If $d=2$ or $d=1$ we get similar expressions, noting for example
that if $d=1$, $\tPsi(p,u)(X_1)=0$ (see Prop.5.10), and so
\[
\langle \tilde{G}(p,u)(X_1),[\tPsi(p,u)(X_2), \tPsi(p,u)(X_3)]\rangle=
\langle \tilde{G}(p,u)\wedge [\tPsi(p,u), \tPsi(p,u)]\rangle (X_1,X_2,X_3).
\]
Using (5.73) and (5.74) 
we see that $\eta(\Phi)=\eta(\tilde{\Phi})$ is bounded by an $L^1$ form
on $N^1\Sigma$, and so
we can apply the dominate convergence theorem to
interchange $\int$ with $\lim_{\epsilon\ra 0}$ in (5.80),
and the expression of Corollary 1.1 is proved \qed.\\[2mm]
\em Remark. \em
For any symmetric
tensor $S\in C^{\infty}(TM^*\otimes TM$, 
$\langle S\wedge [\Phi^{-1}\lnab{}\Phi,
\Phi^{-1}\lnab{}\Phi]\rangle =0$. Thus the condition
of $\langle \Phi^{-1}\lnab{}\Phi\wedge [\Phi^{-1}\lnab{}\Phi,
\Phi^{-1}\lnab{}\Phi]\rangle =0$ is a quite weaker condition then
$\Phi^{-1}\lnab{X}\Phi$ to be symmetric, for each vector field $X$.
In this case, if $d_i\neq 2$ $\forall i$ then 
$p_1(\sm{\bigwedge^2_-}NM)=p_1(\sm{\bigwedge^2_-}TM)$.
\section{ $J$-K\"{a}hler submanifolds}
\setcounter{Th}{0}
\setcounter{Pp}{0}
\setcounter{Cr} {0}
\setcounter{Lm} {0}
\setcounter{equation} {0}
Assume $M$ is a K\"{a}hler submanifold of $N$. 
If $E$ is a rank-4 Hermitian vector bundle over
$M$ with a complex structure $J^E$ and a unitary connection $\lnabe{}$,
the curvature is $J^E$-invariant. 
Let $B=(E_1, E_2=J^EE_1, E_3,E_4= J^EE_3)$ be a local o.n. frame of $E$,
and $\Xi_s^+$ defined as $J^B_s$ in (3.14), and $\Xi_s^-$ defined in
the same way, but replacing $E_4$ by $-E_4$.
Let $c_1(E)$ and $c_2(E)$ be the first and the second Chern classes
of $E$. Then
\begin{eqnarray}
2\pi c_1(E) &=& R^E(\Xi_1^+)~~~~~~~~~
R^E(\Xi_{s}^+)=0~~\mbox{for}~s=2,3.\\
c_2(E) &=&{\cal X}(E)~~~~~~~~~
~~~~p_1(E)=-c_2(E^c)=-2c_2(E)+ c_1(E)^2
\end{eqnarray}
If $E=TM$ we denote $J^E$ by $J$, $E_i$ by $e_i$, and
$\Xi_s^{\pm}$ by $\Lambda^{\pm}_s$.
Then ${2\pi}c_1(M)(X,Y)=Ricci^M(JX,Y)=R^M(\Lambda^+_1,X\wedge Y)
=R^M(X\wedge Y,\Lambda^+_1)$.
The first equation (6.1) implies that $c_1(E)=0$ iff $\sm{\bigwedge^2_+}E$ is
flat. 
We also recall that (see e.g. \cite{[Be]})
\begin{eqnarray}
{\cal X}(E)&=& \frac{1}{8\pi^2}\La{(} \|(R^E)^+_+\|^2-
\|(R^E)^+_-\|^2 -\|(R^E)^-_+\|^2+
\|(R^E)^-_-\|^2\La{)}Vol_M\\[-1mm]
p_1(E)&=& \frac{1}{4\pi^2}\La{(} \|(R^E)^+_+\|^2+
\|(R^E)^+_-\|^2 -\|(R^E)^-_+\|^2-
\|(R^E)^-_-\|^2\La{)}Vol_M
\end{eqnarray}
Since $M$ is K\"{a}hler , $NM$ is a Hermitian vector bundle and $\lnabo{}$
is  a unitary connection, and
\begin{eqnarray}
\sm{F^*c_1(N)=c_1(M)+c_1(NM)}&~~~ & 
\sm{F^*c_2(N)=c_2(M)+c_1(M)\wedge c_1(NM)+c_2(NM)}\\
\sm{p_1(M)+2{\cal X}(M)=c_1(M)^2}&~~~ & 
\sm{p_1(NM)+2{\cal X}(NM)=c_1(NM)^2}\\
\sm{p_1(M)-2{\cal X}(M)=c_1(M)^2-4c_2(M)}&~~~ & 
\sm{p_1(NM)-2{\cal X}(NM)=c_1(NM)^2-4c_2(NM)}
\end{eqnarray}
We define for  $U,V\in NM_p$,  $X,Y\in T_pM$
\begin{equation}
 Ricci^{\bot}(U\wedge V)=R^{\bot}(\Lambda^+_1, U\wedge V),~~~~~~
Ricci_{\bot}(X\wedge Y)= R^{\bot}(X\wedge Y, \Xi^+_1)=2\pi c_1(NM)
\end{equation}
\begin{Lm}
If $M$ is a K\"{a}hler submanifold
of $N$, then\\[3mm]
$(1)~~
(\|(R^{\bot})^+_+\|^2-\|(R^{\bot})^-_+\|^2)=
\ha Ricci^{\bot}\wedge Ricci^{\bot}(E_1,E_2,E_3,E_4)$\\
$~~~~~~~(\|(R^{\bot})^+_+\|^2-\|(R^{\bot})^+_-\|^2)=
\ha Ricci_{\bot}\wedge Ricci_{\bot}(e_1,e_2,e_3,e_4)$.\\[1mm]
$(2)~~p_1(\bigwedge^2_+NM)
=\frac{1}{4\pi^2} Ricci^{\bot}\wedge Ricci^{\bot}
(E_1,E_2,E_3,E_4)\mbox{Vol}_M
=\frac{1}{4\pi^2} Ricci_{\bot}\wedge
Ricci_{\bot}$.\\[1mm]
$(3)$~~ 
$\|(R^{\bot})^+_-\|=\|(R^{\bot})^-_+\|$
and  $\|(R^{\wedge^2_-NM})^+\|=\|(R^{\wedge^2_+NM})^-\|$.\\[2mm]
If we replace $R^{\bot}$ by $R^M$  the same
equalities holds.
\end{Lm}
\noindent
\em Proof. \em (1) Using (6.1)\\[-6mm] 
\begin{eqnarray*}
\lefteqn{4\La{(}\|{R^{\bot}}^+_+\|^2 - \|{R^{\bot}}^-_+\|^2\La{)}
=\sm{\sum_{ts}}\La{(}R^{\bot}_{\Xi^+_{s}}(\Lambda^+_{t})\La{)}^2
-\La{(}R^{\bot}_{\Xi^-_{s}}(\Lambda^+_{t})\La{)}^2}\\
&=&\La{(}{Ricci}^{\bot}(\Xi_1^+)\La{)}^2-
2\La{(}{Ricci}^{\bot}(\Xi_1^-)\La{)}^2
-2\La{(}{Ricci}^{\bot}(\Xi_2^-)\La{)}^2
-2\La{(}{Ricci}^{\bot}(\Xi_3^-)\La{)}^2\\
&=&\sm{\La{(} Ricci^{\bot} (E_1\wedge JE_1)
+Ricci^{\bot}(E_3,JE_3)\La{)}^2
-\La{(} Ricci^{\bot}(E_1\wedge JE_1)-
 Ricci^{\bot}(E_3\wedge JE_3)\La{)}^2}\\
&&\sm{-\La{(} Ricci^{\bot} (E_1\wedge E_3)
+Ricci^{\bot}(JE_1\wedge JE_3)\La{)}^2
-\La{(} Ricci^{\bot}(E_1\wedge E_4)
 +Ricci^{\bot}(JE_1\wedge JE_4)\La{)}^2}\\
&=&4\La{(}\sm{Ricci^{\bot}(E_1\wedge JE_1)Ricci^{\bot}
(E_3\wedge JE_3)-(Ricci^{\bot}(E_1\wedge E_3))^2
-(Ricci(E_1\wedge E_4))^2}\La{)}\\
&=&2\La{(}\sm{Ricci^{\bot}(E_1\wedge JE_1)Ricci^{\bot}
(E_3\wedge JE_3)
-Ricci^{\bot}(E_1\wedge E_3)Ricci^{\bot}(E_2\wedge E_4)}\\[-2mm]
&&\sm{
+Ricci^{\bot}(E_1\wedge E_4)Ricci^{\bot}(E_2\wedge E_3)}\La{)}\\
&=& 4 {Ricci}^{\bot}\wedge {Ricci}^{\bot}(E_1,E_2,E_3,E_4)
\end{eqnarray*}
Similar for the second equality. From (6.3),(6.4), (1) and (5.1),
$p_1(\sm{\bigwedge^2_+}NM)=
\frac{1}{2\pi^2}(\|(R^{\bot})^+_+\|^2-\|(R^{\bot})^-_+\|^2)
\mbox{Vol}
=\frac{1}{4\pi^2}Ricci^{\bot}\wedge Ricci^{\bot}
(E_1,E_2,E_3,E_4)\mbox{Vol}.$
But on the other hand by  (6.7) and (6.8)
$p_1(\sm{\bigwedge^2_+NM})=c_1^2(NM)=\frac{1}{4\pi^2} Ricci_{\bot}\wedge
Ricci_{\bot}$.
Thus $Ricci^{\bot}\wedge Ricci^{\bot}(E_1,E_2,E_3,E_4)
=Ricci_{\bot}\wedge Ricci_{\bot}(e_1,e_2,e_3,e_4)$.
So we have obtained  (2).
(3) Follows immediately from (1), (2) and that 
denoting by $\{1,2,3\}=\{\Xi_1^{\pm},\Xi_2^{\pm},
\Xi_3^{\pm}\}$, one has
$(R^{\wedge^2_{+}}_{ab})^+=\epsilon\sqrt{2}(R^E)^+_{\Xi^+_c}$, 
$(R^{\wedge^2_{+}}_{ab})^-=\epsilon\sqrt{2}(R^E)^-_{\Xi^+_c}$,
$(R^{\wedge^2_{-}}_{ab})^+=\epsilon\sqrt{2}(R^E)^+_{\Xi^-_c}$, 
$(R^{\wedge^2_{-}}_{ab})^-=\epsilon\sqrt{2}(R^E)^-_{\Xi^-_c}$, 
where $\{a,b,c\}$ is a permutation of $\{1,2,3\}$ of signature $\epsilon$.
\qed\\[-2mm]

\begin{Pp} If $M$ is a complex submanifold of $N$, and $c_1(N)=0$, then:
\\[2mm]
$(1)$~ $p_1(\sm{\bigwedge^2_+}NM)=p_1(\sm{\bigwedge^2_+}TM)$.\\[1mm]
$(2)$~$p_1(\sm{\bigwedge^2_-}NM)-p_1(\sm{\bigwedge^2_-}TM)=
4(-F^*c_2(N)+2c_2(M)-c_1(M)^2)$.\\[1mm]
$(3)$~If $c_1(M)=0$, then both $\bigwedge^2_+TM$
and $\bigwedge^2_+NM$ are flat and both
$\bigwedge^2_-TM$
and $\bigwedge^2_-NM$ are anti-self-dual. Moreover,
${\cal X}(M)\geq 0$ (resp. ${\cal X}(NM)\geq 0$)
 with equality to zero iff $M$ (resp. $NM$) is flat.
Furthermore, $F^*c_2(N)[M]\geq 0$ with equality to zero iff
$M$ and $NM$ are flat. 
\end{Pp}
\noindent
\em Proof. \em 
From  (6.5), if $c_1(N)=0$ then $c_1(M)=-c_1(NM)$ and so
by (6.6) $p_1(\bigwedge^2_+NM)=p_1(\bigwedge^2_+TM)$. Now
$p_1(\bigwedge^2_-NM)-p_1(\bigwedge^2_-TM)$ $=-4{\cal X}(NM)+4{\cal X}(M)$
$=-4(c_2(NM)+4c_2(M)$. (6.5) gives the first equality in (2).
Now we assume $c_1(M)=0$. Since $0=F^*c_1(N)=c_1(M)+c_1(NM)$ along $M$,
then $c_1(NM)=0$, and  both $\bigwedge^2_+TM$
and $\bigwedge^2_+NM$ are  flat.
By Lemma 6.1(3) $(R^{\wedge^2_-NM})^+=(R^{\wedge^2_-TM})^+=0$, and so
$\bigwedge^2_-NM$ and $\bigwedge^2_-TM$ are anti-selfdual.
From (6.8) $Ricci_{\bot}=0$, and since
$(R^{\bot})(\Xi^+_{s})=0~~\forall s=1,2,3$,
then $(R^{\bot})^-_+=(R^{\bot})^+_+=0$.
From Lemma 6.1(3) we get $(R^{\bot})^+_-=0$ as well.
 The same holds for $R^M$. 
The final statement is a consequence of the
previous ones and that by (6.3)(6.4), with $E=TM$ or $NM$,
$-p_1(E)=2{\cal X}(E)=\frac{1}{4\pi^2}\int_M\|(R^E)^-_-\|^2.\qed\\[4mm]$
\em Remark. \em Part of Prop.6.1 (4) is a particular case of some results in
\cite{[Hi]} and in \cite{[At-Hi-Si]}.
\section{$I$-K\"{a}hler submanifolds}
\setcounter{Th}{0}
\setcounter{Pp}{0}
\setcounter{Cr} {0}
\setcounter{Lm} {0}
\setcounter{equation} {0}
In subsection 3.2 we saw that
if $N$ is an HK manifold of complex dimension 4
and $M$ is an $I$-K\"{a}hler submanifold,
then  the zero set $\Sigma$ of $F^*\omega_J$ is the zero set of
 a globally defined $I$-holomorphic (2,0)-form $\varphi$ on $M$.
Thus, $\Sigma$ is a locally finite union of
irreducible $I$-complex hypersurfaces $\Sigma_i$, and $\varphi$ 
vanish to order $a_i$ along $\Sigma_i$. Since
$2\cos^2\theta= \|F^*\omega_J\|^2 = 2\|\varphi\|^2$,
 $\cos\theta$ vanish to homogeneous order $a_i$ along $\Sigma_i$. 
$D=\sum_i a_i \Sigma_i$ is a divisor of $\varphi$, and  for any 
closed 2-form  $\phi$ of $M$
\begin{equation}
\int_M-\frac{i}{\pi}\partial\bar{\partial}\log \|\varphi\|\wedge
\phi =\int_D\phi.\\[2mm]
\end{equation}
\noindent
\em Proof of Proposition 1.1 \em
By Theorem 3.1 we have
$-{i}\partial\bar{\partial}\log \|\varphi\|\wedge \omega_I=
Ricci(I(\cdot),\cdot)\wedge \omega_I= \ha s^MVol_M. $
If we take in (7.1) $\phi= \omega_I$, we get
$\frac{1}{\pi}\kappa_2(M)=\int_D\omega_I=\sum_i a_i\int_{\Sigma_i}\omega_I$.
 \qed \\[5mm]
If $I$ does not exist globally on $M$, we still can obtain a residue formula
under some conditions.
In \cite{[S-PV]} we introduced the notion of controlled zero set for a
function on $M$ with zero set a submanifold $\Sigma$.
For each $(p,u)\in N^1\Sigma$ define
$1\leq \kappa(p,u)\leq +\infty$
the order of the zero of $\varphi_{(p,u)}(r)= \cos^2\theta(exp_{p}(ru))$
at $r=0$.  We will say that $\cos^2\theta$
has a \em controlled zero set \em if
there exist a nonnegative integrable function $f:N^1\Sigma
\ra [0,+\infty]$ and $r_0>0$ s.t.
$~sup_{0<r<r_0}|r\frac{d}{d r}\log(\varphi_{(p,u)}(r))|\leq f(p,u)~$
a.e.\ $(p,u)\in N^1\Sigma$.
 For each $p\in \Sigma$, $S(p,1)$ denotes the
unit sphere of $T_p\Sigma^{\bot}$ and $\sigma_{d'}$ its volume. The function
 $\tilde{\kappa}(p)=\frac{1}{\sigma_{d'}}
\int_{S(p,1)}\kappa(p,u)d_{S(p,1)} u$  is the \em average order \em
of the zero $p$ of $\cos^2\theta$, in the normal direction. Next
proposition has a very similar proof to the one of Theorem 1.2 of 
\cite{[S-PV]}, so we omit it.
\begin{Pp} Assume  $(N,J,g)$ is   Ricci-flat KE
 and  $M$ is $I$-K\"{a}hler, closed, and 
$\Sigma$ is a finite disjoint union of closed
submanifolds $\Sigma_i$ with dimension $d_i\leq 2$  and let
$\bigcup_{\ga}k_{i\ga}$ be the range set of $\kappa$ on $N^1\Sigma^i$
 and define $N^1\Sigma_i^{\ga}=\kappa^{-1}(k_{i\ga})$.
If $\kappa$ is bounded a.e.\ and
$\cos\theta$ has controlled zero set, then
\[k_2(M)=-\sum_{i:d_i=2}\pi
\int_{\Sigma_i}\tilde{\kappa}(p)Vol_{\Sigma_i}=-\ha
\sum_{i: d_i=2}\sum_{\ga}k_{i\ga} \, Vol_{N^1\Sigma_i}(N^1\Sigma^{\ga}_i).
\]
\end{Pp}
\noindent
As a consequence we have got a removable high rank singularity theorem:
\begin{Cr} In the conditions of Prop.7.1, $k_2(M)\leq 0$,
with equality to zero iff
$\Sigma_i=0$ $\forall i: d_i=2$.
\end{Cr}
\noindent
Now we prove
\begin{Pp}  Let $M$ be  closed Cayley submanifold of a Ricci-flat
K\"{a}hler-Einstein 4-fold
 $(N,J,g)$, that is not $J$-complex
neither $J$-Lagrangian but it is
 $I$-K\"{a}hler on a open dense set $U$ of  $M\sim {\cal L}$. Then\\[1mm]
$(1)$~$\forall p\geq 1,~\int_M \cos^{2p}\theta s^M \mbox{Vol}_M \leq 0~$.
Consequently, $s^M\geq 0$ iff $s^M=0$. If that is the case then 
$\cos\theta$ is constant. \\[2mm]
$(2)$~If $M$ is immersed without $J$-Lagrangian  points, then
$2\kappa_2(M)=\int_M s^M \mbox{Vol}_M = 0.$
\end{Pp}
\noindent
\em Proof. \em  
From (3.20)  $\Delta\cos^{2p}\theta=p \cos^{2p}\theta s^M+
4 p^2\cos^{2p-2}\|\nabla \cos\theta\|^2$.
Integration  and  Stokes gives the inequalities in (1).
If $s^M\geq 0$ and since the set of $J$-Lagrangian 
points has empty interior, (1) implies $s^M=0$, and so
$\Delta\cos^{2p}\theta\geq 0$.
Thus,  $\cos\theta$
is constant. Integration of (3.20)  under the 
assumption of ${\cal L}=\emptyset$  proves (2). 
\qed
\begin{Cr} If $M$ is a closed $I$-complex 4-submanifold
of an HK manifold $(N,I,J,K,g)$ of real dimension 8,
 and if $s^M> 0$, then $M$ is a totally complex submanifold.
\end{Cr}
\noindent
\em Proof. \em
Quaternionic submanifolds are HK, and so
Ricci-flat, what is not possible. If we assume
$M$ is not totally complex, 
by  Proposition 7.1  $s^M$ should vanish.
\qed 
\begin{Pp} If $M$ is a closed $I$-complex 4-submanifold
of an HK manifold $(N,I,J,K,g)$ of real dimension 8,
 at quaternionic points $s^M\leq 0$.
\end{Pp}
\noindent
\em Proof. \em  Quaternionic points are maximum points of
$\cos\theta$. Thus, by (3.20) $s^M\leq 0$.\qed
\begin{Pp} Let $M$ be a Cayley submanifold of a Ricci flat KE
8-manifold $(N,J,g)$, that is neither $J$-complex
nor $J$-Lagrangian and it is a $I$-K\"{a}hler on
a open set $O$ of  $M$. If $\cos\theta$ is constant on $O$
then $(M,I,\Jw, I\Jw)$ is HK on $O$.
\end{Pp}
\noindent
\em Proof. \em Since $\cos\theta$ is constant on $O$, 
by Prop.3.6, $\Jw$ is K\"{a}hler on $O$, and so, 
 $M$ is HK on $O$.\qed\\[4mm]
The following proposition was already announced in \cite{[A-PV-S]} and
can be also seen as  a corollary of the above propositions:
\begin{Th}  (\cite{[A-PV-S]}) If $M$ is a closed $I$-complex 4-submanifold of
an HK manifold $(N,I,J,K,g)$ of real dimension 8, and $M$
is neither a quaternionic nor a totally complex submanifold,  
then the following assertions are equivalent to each other:\\[1mm]
$(a)$~~$s^M=0$\\
$(b)$~~ $\cos\theta_J$ is constant\\
$(c)$~~$(M,I,J_{\omega_J}, J_{\omega_K}, g)$ is  HK\\
$(d)$~~the quaternionic angle of $M$ is constant.
\end{Th}
\noindent
\em Remark. \em (b) implies (a) and (c) , and (b)$\Longleftrightarrow $
(c)$\Longleftrightarrow $(d)
do not need compactness of $M$. The proof (d)
$\Longleftrightarrow $ (b) is shown in \cite{[A-PV-S]}. 
In \cite {[A-M]} we can find
related results.
\section{Cayley submanifolds of $\RR^{8}$}
\setcounter{Th}{0}
\setcounter{Pp}{0}
\setcounter{Cr} {0}
\setcounter{Lm} {0}
\setcounter{equation} {0}
\subsection{Complex Cayley graphs}
We consider $\RR^{8}=\RR^{4}\times \RR^{4}$, with the euclidean metric $g_0$
( in $\RR^{4}$, and so in $\RR^{8}$).
The set of $g$-orthogonal complex structures of $\RR^{8}$ has two
connected components. Let us fix $J_0$ given by
$J_0(X,Y)=(-Y,X)$, and denote by $\omega_0$ the
K\"{a}hler form . If $f:\RR^{4}\ra \RR^{4}$ is a
smooth map, the graph of $f$ is the map
$\Gamma_f:\RR^{4}\ra \RR^{8}$, ~$\Gamma_f(x)=(x,f(x))$. 
In \cite{[D-S]} we compute the K\"{a}hler angles of 
$\Gamma_f$ with respect to $J_0$.
Let $g_M$ be the  the graph metric on $\RR^{4}$, $g_M=(\Gamma_f)^*g_0$. 
Note that $\Gamma_f^*\omega_0(X,Y)= g_0(-df(X)+df^t(X),Y)$,
where $g_0$ is w.r.t. $\RR^{4}$. Using the musical isomorphism w.r.t. 
the Euclidean metric  $g_0$ on $\RR^{4}$ we have
$g_M =Id +df^t\circ df$ and $\Gamma_f^*\omega_0 =-df + df^t$. 
The metric $g_M$ is complete if $f$ is defined on all $\RR^{4}$.
 The solutions of
\begin{equation}
det(\Gamma_f^*\omega_0- \lambda g_M)=0
\end{equation}
 are pure imaginary,
and $\lambda^2=-\cos^2\theta_{\al}$ give the K\"{a}hler angles.
We can compute explicitly (8.1).
 Set $f(x,y,z,w)=(u,v,s,t)$. Then
\[df=\left[\begin{array}{cccc}
\frac{\partial u}{\partial x} &\frac{\partial u}{\partial y}&
\frac{\partial u}{\partial z}&\frac{\partial u}{\partial w}\\
\frac{\partial v}{\partial x} &\frac{\partial v}{\partial y}&
\frac{\partial v}{\partial z}&\frac{\partial v}{\partial w}\\
\frac{\partial s}{\partial x} &\frac{\partial s}{\partial y}&
\frac{\partial s}{\partial z}&\frac{\partial s}{\partial w}\\
\frac{\partial t}{\partial x} &\frac{\partial t}{\partial y}&
\frac{\partial t}{\partial z}&\frac{\partial t}{\partial w}\\
\end{array}\right]\]
Now define
\begin{eqnarray*}
\begin{array}{llllll}
A=-\frac{\partial u}{\partial y}+\frac{\partial v}{\partial x}&
B=\frac{\partial s}{\partial x}-\frac{\partial u}{\partial z}&
C=\frac{\partial t}{\partial x}-\frac{\partial u}{\partial w}&
D=\frac{\partial s}{\partial y}-\frac{\partial v}{\partial z}&
E=\frac{\partial t}{\partial y}-\frac{\partial v}{\partial w}&
F=\frac{\partial t}{\partial z}-\frac{\partial s}{\partial w}~~~~\\[2mm]
l=\langle \frac{\partial f}{\partial y},
\frac{\partial f}{\partial w}\rangle &
m=\langle \frac{\partial f}{\partial z},
\frac{\partial f}{\partial w}\rangle &
p=\langle \frac{\partial f}{\partial x},
\frac{\partial f}{\partial y}\rangle &
q=\langle \frac{\partial f}{\partial x},
\frac{\partial f}{\partial z}\rangle &
r=\langle \frac{\partial f}{\partial x},
\frac{\partial f}{\partial w}\rangle &
k=\langle \frac{\partial f}{\partial y},
\frac{\partial f}{\partial z}\rangle \\[2mm]
h=(1\!+\!\|\frac{\partial f}{\partial x}\|^2)
& o=(1\!+\!\|\frac{\partial f}{\partial y}\|^2) &
 d=(1\!+\!\|\frac{\partial f}{\partial z}\|^2) &
 n=(1\!+\!\|\frac{\partial f}{\partial w}\|^2) & &
\end{array}\end{eqnarray*}
\begin{eqnarray*}
{\cal A}&=& 2hlkm +hodn-h(dl^2+om^2+nk^2)+
p^2(-dn+m^2)+q^2(l^2-no)+r^2(-od+k^2)\\
&&+2qm(-lp+or)+2pr(-mk+dl)+2qk(pn-rl)\\[2mm]
{\cal B}&=& 2DE(qr-hm)+2BE(-rk+pm)+2BD(lr-np)+2CE(-dp+qr)\\
&&+2AE(dr-qm)+2CF(-oq+pk)+2CB(-om+kl)+2DF(-rp+hl)\\
&&+2AF(-rk+ql)+2AD(-rm+nq)+2AC(-dl+mk)+2AB(ml-nk)\\
&&+2CD(-ql+mp)+2FE(qp-kh)+2FB(or-pl)+E^2(dh-q^2)\\
&&+B^2(no-l^2)+o(hF^2 +dC^2)+n(hD^2+dA^2)
-c^2k^2-r^2D^2-m^2A^2-p^2F^2\\[2mm]
{\cal D}&=&(AF-BE+CD)^2.
\end{eqnarray*}
The K\"{a}hler angles  of $\Gamma_f$ are the solutions
of (8.1) for
$-\lambda^2=\mu=\cos^2\theta$
what explicitly reads
  $\mu^2{\cal A}-\mu{\cal B}+{\cal D}=0$.
Thus $\Gamma_f$ has e.k.a. iff
${\cal A}\neq 0$ and ${\cal B}^2=4{\cal A}{\cal D}$,
and in this case $\cos^2\theta=\frac{{\cal B}}{2{\cal A}}=\sqrt{
\frac{{\cal D}}{{\cal A}}}$, or
$~{\cal A}=0$ and in this case 
$\cos^2\theta=\frac{{\cal D}}{{\cal B}}.$ 
We can find a very large
family of Cayley submanifolds $M$ in a hyper-K\"{a}hler ambient space $N$
by taking two different complex structures $J_x,J_y$, and considering
$M$ $J_x$-complex and $N$ $J_y$-complex, where $x,y$ are any elements
of $S^2$. Those submanifolds are
automatically minimal, and the expression of the k.a is simplified
(see Prop.3.5). We will restrict ourselves
to this case. \\[4mm]
If we consider on $\RR^{4}$ a $g_0$-orthogonal
complex structure $\Jw$, the complex structure
$(\Jw, -\Jw)$ of $\RR^{8}$,
$(\Jw, -\Jw)(X,Y)=(\Jw X, -\Jw Y)$
is $g_0$-orthogonal and anti-commutes with $J_0$.
Then $(J_0, (\Jw, -\Jw), J_0\times (\Jw, -\Jw))$ defines an Hyper-
K\"{a}hler structure on $\RR^{8}$. If $f:\RR^{4}\ra \RR^{4}$ is a
$\Jw$-anti-holomorphic map, then the graph of $f$
is a $(\Jw, -\Jw)$-complex submanifold of $\RR^{8}$.
 We consider $\Jw$ the   complex structure
$i=\frac{1}{\sqrt{2}}(e_1\wedge e_2 + e_3\wedge e_4)$,
where $e_1,e_2, e_3, e_4$ is
the canonic basis of $\RR^{4}\equiv \RR^{4}\times {0}\subset \RR^{8}$.
Recall that, considering $\CC$ with the usual complex structure, 
also denoted by $J_0$, $J_0(x,y)=(-y,x)$, and
if $f(x,y)=(u,v):\RR^{2}\equiv \CC \ra \RR^{2}\equiv \CC$
then $f(x,y)=(u,v)$ is anti-holomorphic iff $df\circ J_0 =-J_0 \circ df$,
iff $\frac{\partial u}{\partial y}=\frac{\partial v}{\partial x}$
and $\frac{\partial v}{\partial y}=-\frac{\partial u}{\partial x}$, that is,
iff $h(x,y)=(v,u)$ is holomorphic.
Let $f:\RR^{4}\ra\RR^{4}$, $f(x,y,z,w)=(s,t,u,v)$. Then $f$ is 
anti-$i$-holomorphic iff $(x,y)\ra(u,v)$, $(x,y)\ra(s,t)$,
$(z,w)\ra(u,v)$ and $(z,w)\ra(s,t)$ are 
anti-holomorphic, iff
\begin{eqnarray}
&& \left\{\begin{array}{cccc}
\frac{\partial u}{\partial x}=-\frac{\partial v}{\partial y}~~~&
\frac{\partial u}{\partial y}=\frac{\partial v}{\partial x}~~~&
\frac{\partial u}{\partial z}=-\frac{\partial v}{\partial w}~~~&
\frac{\partial u}{\partial w}=\frac{\partial v}{\partial z}\\[2mm]
\frac{\partial s}{\partial x}=-\frac{\partial t}{\partial y}~~~&
\frac{\partial s}{\partial y}=\frac{\partial t}{\partial x}~~~&
\frac{\partial s}{\partial z}=-\frac{\partial t}{\partial w}~~~&
\frac{\partial s}{\partial w}=\frac{\partial t}{\partial z}
\end{array}\right.
\end{eqnarray}
This implies
\[\begin{array}{l}
 A=F=0;~~~B=-E=\frac{\partial s}{\partial x}-\frac{\partial u}{\partial z};
~~~C=D=\frac{\partial s}{\partial y}-\frac{\partial u}{\partial w};~~~
 p=\langle \frac{\partial f}{\partial x},\frac{\partial f}{\partial y}
\rangle=0;~~~m=\langle \frac{\partial f}{\partial z},
\frac{\partial f}{\partial w}\rangle=0;\\[2mm]
 q=\langle \frac{\partial f}{\partial x},\frac{\partial f}{\partial z}
\rangle=\frac{\partial u}{\partial x}\frac{\partial u}{\partial z}
+\frac{\partial u}{\partial y}\frac{\partial u}{\partial w}
+\frac{\partial s}{\partial x}\frac{\partial s}{\partial z}
+\frac{\partial s}{\partial y}\frac{\partial s}{\partial w}
=l=\langle \frac{\partial f}{\partial y},\frac{\partial f}{\partial w}
\rangle\\[2mm]
 r=\langle \frac{\partial f}{\partial x},\frac{\partial f}{\partial w}
\rangle=-\frac{\partial u}{\partial y}\frac{\partial u}{\partial z}
-\frac{\partial v}{\partial y}\frac{\partial v}{\partial z}
-\frac{\partial t}{\partial y}\frac{\partial t}{\partial z}
-\frac{\partial s}{\partial y}\frac{\partial s}{\partial z}
=-k=-\langle \frac{\partial f}{\partial y},
\frac{\partial f}{\partial z}\rangle\\[2mm]
h=(1+\|\frac{\partial f}{\partial x}\|^2)=
(1+\|\frac{\partial f}{\partial y}\|^2)=o;~~~
d=(1+\|\frac{\partial f}{\partial z}\|^2)=
(1+\|\frac{\partial f}{\partial w}\|^2)=n
\end{array}\]
\[\begin{array}{l}
{\cal A} = (hd-q^2-k^2)^2\geq 1;~~~~~~~
{\cal B} = 4BCkq +2(B^2+C^2)(dh-k^2-q^2)\\[2mm]
{\cal D} = (B^2+ C^2)^2; ~~~~~~~~
\cos^2\theta = \frac{B^2+C^2}{hd-k^2-q^2}.
\end{array}\]
Note that,  the linear map
$(u_0,v_0):\RR^{4}\ra \CC\equiv \RR^{2}$
\begin{equation}
(u_0,v_0)(x,y,z,w)= (x+y+z+w,x-y+z-w)
\end{equation}
is anti-holomorphic, considering  $~\CC= \RR^{2}$ and $\RR^{4}$ with the
complex structures $J_0$ and $i =J_0\times J_0$, respectively,
or equivalently,
 $(u_0,v_0)$ satisfies the first eq. of (8.2)
\begin{Pp} If $f$ is anti-$i$-holomorphic and at a point $p$, $r=k=0$
that is, at $p$,~
$\frac{\partial u}{\partial y}\frac{\partial u}{\partial z}-
\frac{\partial u}{\partial x}\frac{\partial u}{\partial w}=
-\frac{\partial s}{\partial y}\frac{\partial s}{\partial z}+
\frac{\partial s}{\partial x}\frac{\partial s}{\partial w},$
 then $\Gamma_f$
is a minimal submanifold with e.k.a. $\theta$ at $p$ given by
\[\cos^2\theta=\frac{ (\frac{\partial s}{\partial x}-
\frac{\partial u}{\partial z})^2+(\frac{\partial s}{\partial y}
-\frac{\partial u}{\partial w})^2}{
(1+\|\frac{\partial f}{\partial x}\|^2)
(1+\|\frac{\partial f}{\partial z}\|^2)-\langle
\frac{\partial f}{\partial x}, \frac{\partial f}{\partial z}\rangle^2}
\]
\end{Pp}
\noindent
\em Proof. \em  We only have to apply the above formulas, and the fact that
since $f$ is anti-$i$-holomorphic
$\frac{\partial f}{\partial x}=
(\frac{\partial u}{\partial x},\frac{\partial u}{\partial y},
\frac{\partial s}{\partial x},\frac{\partial s}{\partial y}), 
$$~\frac{\partial f}{\partial y}=
(\frac{\partial u}{\partial y},-\frac{\partial u}{\partial x},
\frac{\partial s}{\partial y},-\frac{\partial s}{\partial x}), 
$$~\frac{\partial f}{\partial z}=
(\frac{\partial u}{\partial z},\frac{\partial u}{\partial w},
\frac{\partial s}{\partial z},\frac{\partial s}{\partial w}), 
$$~\frac{\partial f}{\partial w}=
(\frac{\partial u}{\partial w},-\frac{\partial u}{\partial z},
\frac{\partial s}{\partial w},-\frac{\partial s}{\partial z})$ 
and so
$\|\frac{\partial f}{\partial y}\|=\|\frac{\partial f}{\partial x}\|$.
\qed
\begin{Cr}
If $f:\RR^{4}\ra \RR^{4}$ satisfies\\[-6mm]
\begin{eqnarray}\begin{array}{ll}
\frac{\partial u}{\partial x}=-\frac{\partial v}{\partial y}=
\frac{\partial u}{\partial z}=-\frac{\partial v}{\partial w}
&~~~~\frac{\partial v}{\partial x}=
\frac{\partial u}{\partial y}=\frac{\partial v}{\partial z}=
\frac{\partial u}{\partial w}\\
\frac{\partial s}{\partial x}=
-\frac{\partial t}{\partial y}=
\frac{\partial s}{\partial z}=-\frac{\partial t}{\partial w}
&~~~~\frac{\partial t}{\partial x}=
\frac{\partial s}{\partial y}=\frac{\partial t}{\partial z}=
\frac{\partial s}{\partial w}\\[-3mm]
\end{array}
\end{eqnarray}
Then $f$ is in the conditions of Proposition 8.1 at every $p\in \RR^{4}$,
 with $h=o=d=n$ and $q=l=\|\frac{\partial f}{\partial x}\|^2
=\|\frac{\partial f}{\partial z}\|^2=
(\frac{\partial u}{\partial x})^2+
(\frac{\partial u}{\partial y})^2+
(\frac{\partial s}{\partial x})^2+
(\frac{\partial s}{\partial y})^2$, and
\[\cos^2\theta =
\frac{(\frac{\partial s}{\partial x}-\frac{\partial u}{\partial z})^2
+(\frac{\partial s}{\partial y}-\frac{\partial u}{\partial w})^2}
{1+2\La{(}(\frac{\partial u}{\partial x})^2+
(\frac{\partial u}{\partial y})^2+
(\frac{\partial s}{\partial x})^2+
(\frac{\partial s}{\partial y})^2\La{)}}\]
$\Gamma_f$ is a complete Cayley submanifold with no $J_0$-complex points. 
Furthermore,
$(\frac{\partial s}{\partial x}-\frac{\partial u}{\partial z})^2
+(\frac{\partial s}{\partial y}-\frac{\partial u}{\partial w})^2$
is bounded iff
 $\cos^2\theta$ is bounded by a constant $\delta <1$.
\end{Cr}
\noindent
\em Proof. \em 
Set $X=-\frac{\partial u}{\partial x}$, $Y=-\frac{\partial u}{\partial y}$,
$a=(B-2X)^2+(C-2Y)^2=(\frac{\partial s}{\partial x}+
\frac{\partial u}{\partial x})^2 +
(\frac{\partial s}{\partial y}+ \frac{\partial u}{\partial y})^2
$, and $\zeta =B^2+ C^2$.
Then $\cos^2\theta =\frac{\zeta}{1+a+\zeta}$
with $\zeta, a\geq 0$. This function is an increasing function on 
$\zeta$, what implies the last  assertion.\qed
\\[5mm]
\em Remark. \em For any constants $\al,\be$, the map
$f=(\al (u_0,v_0), \be (u_0,v_0))$, where $(u_0,v_0)$ is
given by (8.3), satisfies the conditions of
Corol. 8.1.\\[2mm]
The following example of \cite{[D-S]} was announced in 
\cite{[S1]} (we note that in
\cite{[S1]} is missing a sequareroot on the denominator of the
expression of $\cos\theta$). It is an example on the
conditions of Cor.8.1.
\begin{Pp} ([D-S], [S,1]) Let $\phi(t)=\sin(t)$, $\xi(t)=\sinh(t)$, and
\\[-2mm]
\begin{equation}\begin{array}{l}
u(x,y,z,w)=\phi(x+z)\xi'(y+w)\\
v(x,y,z,w)= -\phi'(x+z)\xi(y+w)
\end{array}\\[-2mm]\end{equation}
then:\\
$(a)$~ If $f=(u,v,u,v)$, $\Gamma_f$ is a complete minimal Lagrangian
submanifold.\\
$(b)$~ If $f=(u,v,-u,-v)$, $\Gamma_f$ is a complete Cayley submanifold
 with e.k.a and\\[-3mm]
\begin{equation}
\cos\theta=2\sqrt{\frac{\cos^2(x+z)+\sinh^2(y+w)}
{1 + 4(\cos^2(x+z)+\sinh^2(y+w))}}
\end{equation}
Thus $\Gamma_f$ has no $J_0$-complex points, but
 $\cos\theta$ assume all values of $[0,1[$. The set of
Lagrangian points is an infinite discrete family of
parallel 2-planes
${\cal L}= \sm{\bigcup_{k\in Z}} \RR^{}\cdot (1,0,-1,0) \oplus
\RR^{}\cdot(0,1,0,-1)+(0,0, \frac{\pi}{2} +k\pi, 0)$.
\end{Pp}
\begin{Pp} Let  $f:\RR^{4}\ra\RR^{4}$ be a map.\\[2mm]
$(1)$~  A point $p_0$  is  a $J_0$-complex  point of $\Gamma_f$
  iff $df(p_0):\RR^{4}\ra\RR^{4}$ is
a complex structure of $\RR^{4}$. If that is the case, then it is 
$g_M$-orthogonal. It is $g_0$-orthogonal iff $g_M=2g_0$.\\[2mm]
$(2)$~ If $\Jw$ is a $g_0$-orthogonal complex structure
of $\RR^{}$ and at a point $p_0$,
  $df(X)=a\Jw(X)$ where $a$ is  any nonzero real number, 
then $\Gamma_f$ has e.k.a. at a point $p_0$, with
 $\cos\theta=\frac{2|a|}{1+a^2}$ and
$\Gamma_f^*\omega_0(X,Y)=g_M(\cos\theta\epsilon \Jw(X),Y)$, where
$\Jw$ is also a $g_M$-orthogonal structure on $\RR^{4}$, and $\epsilon
= sign\, a $.
\end{Pp}
\noindent
\em Proof. \em (1) At $p_0$, $\Gamma_f$ is  a $J_0$-complex 
submanifold iff
$\forall X\in \RR^{4}$ $\exists Y\in \RR^{4}$ s.t.
$(Y,df(Y))=J_0(X,df(X))=(-df(X),X)$.
that is $-df(df(X))=X$, but this is equivalent to $df(p_0):\RR^4\ra \RR^4$ 
to be a complex structure $\Jw$. Now
$g_M(X,Y)=g_0(X,Y)+g_0(\Jw(X),\Jw(Y))$ 
and so $g_M$ is $\Jw$-Hermitian, or equivalently
$\Jw$ is $g_M$-orthogonal. Now easily follows that $\Jw$
is $g_0$-orthogonal iff $g_M=2g_0$.\\[2mm]
 (2) the condition e.k.a, $\Gamma_f\omega_0
=\cos\theta\Jw$ (under a $g_M$-musical isomorphism),  means
$-g_0(df(X),Y)+g_0(X,df(Y))=$ $\cos\theta g_0(\Jw(X),Y)
+$ $\cos\theta g_0(df(\Jw X), df(Y))$
for some $J_M$-orthogonal structure $\Jw$ on $\RR^{4}$.
Obviously if $df(p_0)=a\Jw$ with $\Jw$ $g_0$-orthogonal, then
immediately we verify that $\Gamma_f$ has e.k.a. at $p_0$
with $\cos\theta=\frac{2|a|}{1+a^2}$ and $\Gamma_f^* \omega_0
=\cos\theta\epsilon \Jw$.\qed
\begin{Cr} Assume $f$ is anti-$i$-holomorphic and
$\Gamma_f$ is at a point $p_0$
a $J_0$-complex submanifold of $\RR^{8}$, that is $df(p_0)=\Jw$
where $\Jw$ is a complex structure of $\RR^{4}$, $g_M$-orthogonal.
Then $i$ and $\Jw$ anti-commute and are both $g_M$-orthogonal on $\RR^{4}$.
\end{Cr}
\noindent
\em Proof. \em
Let $X,Y\in \RR^{4}$. From $df(p_0)(iX)=-idf(p_0)(X)$
we have $\Jw\circ i =-i \circ \Jw$.
That is, $\Jw $ and $i$ anti-commute. Now, at $p_0$
$g_M(iX,iY) $ $= g_0(iX,iY)+g_0(df(p_0)(iX),df(p_0)(iY))$ $
= g_0(X,Y) +g_0(-idf(p_0)(X),-idf(p_0)(Y))$ $=g_M(X,Y)$.
So $i$ is also $g_M$-orthogonal. \qed\\[4mm]
Now we are ready to get examples of non-$J_0$-holomorphic Cayley
submanifolds of $(\RR^{8}, J_0,g_0)$ with $J_0$-complex points,
or non-linear Cayley graphs with $\cos\theta \leq \delta <1$.\\[2mm]
Consider the complex structure $j$ of $\RR^{4}$
 $j=\frac{1}{\sqrt{2}}(e_1\wedge e_3-e_2\wedge e_4)$. Then
$ji=-ij$, so $j$
 can be seen as a (linear) anti-$i$-holomorphic
map of $\RR^{4}$, $j(x,y,z,w)=(-z,w,x,-y)$  (with $dj(p)=j$
$\forall p=(x,y,z,w))$.
\begin{Pp} Let $\tilde{f}:\RR^{4}\ra\RR^{4}$ be any anti-$i$-holomorphic
 with $d\tilde{f}(0)=0$. Then $f:\RR^{4}\ra\RR^{4}$, $f=j+\tilde{f}$
is s.t. its graph defines a Cayley submanifold of $\RR^{8}$
with a $J_0$-complex point $0$.
\end{Pp}
\noindent
\em Proof. \em Since $df(0)=j$,
 by  proposition 8.3(1) the tangent space of $\Gamma_f$
is at $0$ a $J_0$-complex subspace of $\RR^{8}$.\qed
\begin{Cr} $(1)$~ $f(x,y,z,w)=j(x,y,z,w)+(x^2-y^2, -2xy, z^2-w^2, -2zw)$
defines a non $J_0$-holomorphic Cayley submanifold of
$\RR^{8}$ with only one $J_0$-complex point, namely at $0$.\\[2mm]
$(2)$~$f(x,y,z,w)=j(x,y,z,w)+(x^2-y^2, -2xy,0,0 )$
defines a non $J_0$-holomorphic Cayley submanifold of
$\RR^{8}$ with set of  $J_0$-complex point ${\cal C}=\RR^{2}\times \{
(0, 0)\}$.
\end{Cr}
\noindent
\em Proof. \em  (1)
From previous proposition $0$ is a $J_0$-complex point.
If $p=(x,y,z,w)$ is a $J_0$-complex point of $\Gamma_f$, then
$df(p)= j + \xi$,  with $(j+\xi)^2=-Id$,  where
\[\xi=(2xe_*^1-2ye_*^2,
-2ye_*^1-2xe_*^2,2ze_*^3-2we_*^4, -2we_*^3-2ze_*^4)\]
From $-Id=(j+\xi)^2= j^2+j\xi +\xi j +\xi^2=-Id +j\xi +\xi j+ \xi^2$,
we should have $j\xi +\xi j=- \xi^2$. But
$\xi^2=4(x^2+ y^2)(e_*^1\otimes e_1 +e_*^2\otimes e_2)
+4(z^2+ w^2)(e_*^3\otimes e_3 +e_*^4\otimes e_4)$
and $(j\xi +\xi j)(e_1)=2(x+z)e_3+2(y-w)e_4$,~
$(j\xi +\xi j)(e_3)=-2(x+z)e_1+2(y-w)e_2$,
and so $j\xi +\xi j=- \xi^2$ is only possible for $p=0$.
The case (2) is similar with
$\xi=(2xe_*^1-2ye_*^2,-2ye_*^1-2xe_*^2,0,0)$.
\qed
\begin{Pp}  Let $\tilde{f}=f +(\al (u_0,v_0), \be (u_0,v_0))$
where $f=(u,v,u,v)$ is given by Prop.$8.2(a)$ and $(u_0,v_0)$ by
$(8.3)$ and $\al, \be$ any constants. Then
 $\tilde{f}$ is anti-i-holomorphic satisfying $(8.4)$, and
$\cos\theta\leq\frac{2(\al^2+\be^2)}{1+ 2(\al^2+\be^2)}<1$.
\end{Pp}
\noindent
\em Proof. \em  Use proof of Cor. 8.1  to check the upper bound
of $\cos\theta$.\qed
\subsection{Coassociative graphs}
A coassociative graph is a Cayley graph of a map $f:\RR^4\ra\RR^3
\!\subset\! \RR^4$ (see \cite{[H-L1]}). In this case at each point $p$, 
$df(p):\RR^4\ra\RR^3\!\subset\! \RR^4$ cannot
be an isomorphism. Thus, by Prop. 8.3 (1) we have:
\begin{Cr} If $\Gamma_f$ is a coassociative graph then it has no
$J_0$-complex points.
\end{Cr}
\noindent
An example of a coassociative graph given in \cite{[H-L1]}
 is the graph of $\eta: \RR^4\ra \RR^3$
$\eta(x)=\frac{\sqrt{5}}{2\|x\|}\bar{x}\epsilon x$ where the
product is the quaternionic product and
$\epsilon$ is a unit of $\RR^3=Im \RR^4$. This is the cone of the
Hopf map from $S^3$ to $S^2$.

\end{document}